\documentclass[12pt]{amsart}

\input{0_Paper_Style.sty}

\newcommand{\bbQ}{\mathbb{Q}}
\newcommand{\bbR}{\mathbb{R}}
\newcommand{\bbZ}{\mathbb{Z}}

\newcommand{\calP}{\mathcal{P}}

\newcommand{\prm}{^\prime}
\newcommand{\parent}[1]{\left( #1 \right)}

\newcommand{\inv}{^{-1}}
\newcommand{\abs}[1]{\left|#1\right|}

\newcommand{\floor}[1]{\left\lfloor #1 \right\rfloor}

\newcommand{\primeres}[1]{\calP^*\parent{#1}} %The largest ell such that primorial(ell) <= n %Formerly logvar
\newcommand{\primorial}[1]{\calP\parent{#1}} %The product over p <= ell of p
\newcommand{\vp}[1]{v_{#1}}

\newcommand{\ksigma}[1]{\sigma^{[#1]}}%{\sigma^{\left[{#1}\right]}}

\newcommand{\krho}[2]{\rho^{\left[{#1}\right]}_{#2}}
\newcommand{\kG}[1]{G^{\left[{#1}\right]}}

\newcommand{\ku}[1]{u^{\left[{#1}\right]}}

\newcommand{\kF}[1]{F^{\left[{#1}\right]}}
\newcommand{\kfp}[2]{f^{\left[{#1}\right]}_{#2}}
\newcommand{\kalpha}[2]{\alpha_{#2}^{\left[{#1}\right]}}

\newcommand{\ka}[2]{a_{#2}^{\left[{#1}\right]}}

\newcommand{\kx}[2]{x_{#2}^{\left[{#1}\right]}}

\newcommand{\kEp}[2]{E_{#2}^{\left[{#1}\right]}}
\newcommand{\kE}[1]{E^{\left[{#1}\right]}}
\newcommand{\kepsilon}[2]{\epsilon_{#2}^{\left[{#1}\right]}}
\newcommand{\kN}[2]{N^{\left[{#1}\right]}\parent{#2}}

\newcommand{\kNi}[2]{N^{\left[{#1}\right]}_{#2}}

\begin{document}
\title{A family of analogues to the Robin criterion}
\author{Steve Fan, Mits Kobayashi, and Grant Molnar}
\address{Department of Mathematics, University of Georgia, Athens, GA 30602}
\email{Steve.Fan@uga.edu}
\address{Department of Mathematics, Dartmouth College, Hanover, NH 03755}
\email{Mits.Kobayashi@dartmouth.edu}
\address{}
\email{molnar.grant.5772@gmail.com}

\begin{abstract}
	The Robin criterion states that the Riemann hypothesis is equivalent to the inequality $\sigma(n) < e^\gamma n \log \log n$ for all $n>5040$, where $\sigma(n)$ is the sum of divisors of $n$, and $\gamma$ is the Euler--Mascheroni constant. Define the family of functions
\[ \ksigma k (n)\coloneqq\sum_{[d_1,\dots,d_k]=n}d_1\dots d_k \]
where $[d_1, \dots, d_k]$ is the least common multiple of $d_1, \dots, d_k$. These functions behave asymptotically like $\sigma(n)^k$ as $k\to\infty$.
We prove the following analogue of the Robin criterion: for any $k \geq 2$, the Riemann hypothesis holds if and only if $\ksigma k (n) < \frac{(e^\gamma n \log \log n)^k}{\zeta(k)}$ for all $n > 2162160$, where $\zeta$ is the Riemann zeta function.
\end{abstract}

\maketitle

\section{Introduction}\label{Section: Introduction}

In 1894, von Sterneck \cite{vonSterneck} introduced arithmetic functions $F$ of the form
\[ F(n)\coloneqq\sum_{[d_1,\dots,d_k]=n}f_1(d_1)\cdots f_k(d_k), \]
where $[d_1, \dots, d_k]$ is the least common multiple of $d_1, \dots, d_k$ and $f_1, \dots, f_k$ are arithmetic functions. (See also Lehmer \cite{Lehmer2, Lehmer}.) In particular, von Sterneck considered $f_1=\cdots=f_k=\varphi$, the Euler totient function, in which case $F$ is the Jordan totient function. Note that the definition of $F$ is equivalent to the identity
\[\sum_{n\geq 1}\frac{F(n)}{n^s}=\sum_{d_1,\dots,d_k\geq 1}\frac{f_1(d_1)\cdots f_n(d_n)}{[d_1, \dots, d_k]^s}.\]

Taking $f_i(n)=n$ for all $i=1,\dots,k$, we make the following definition.
\begin{definition}\label{Definition: ksigma (n)}
	For $k \geq 1$ an integer, and $n \in \bbZ_{>0}$, we define $\ksigma 1(n):=n$ and
	\[
	\ksigma k (n) \coloneqq \sum_{[d_1, \dots, d_k] = n} d_1 \dots d_k.
	\]
\end{definition}
The function $\ksigma 2 (n)$ is a special case of \cite[(5.10)]{Buschman}, but the study of the family of functions $\ksigma k (n)$ appears to be new.

Lehmer \cite{Lehmer} proved that for any arithmetic function $f : \bbZ_{>0} \to \bbR$, we have
\begin{equation}\label{Equation: LCM identity}
	\sum_{[d_1, \dots, d_k] = n} f(d_1) \dots f(d_k) = \sum_{d \mid n} \mu(n/d) \parent{\sum_{\delta \mid d} f(\delta)}^k,
\end{equation}
where $\mu$ is the M\"obius function. The expression on the right makes sense if the integer $k$ is replaced by any complex number $\kappa$. We call this expression the \emph{$\kappa$th LCM-power of $f$}. In this paper we will focus on the case real $\kappa > 1$.

\begin{definition}\label{Definition: ksigma (n) (revised)}
	For $\kappa$ a positive real number, and $n \in \bbZ_{>0}$, we define 
	\begin{equation}
	\ksigma \kappa (n) \coloneqq \sum_{d \mid n} \mu(n/d) \sigma(d)^\kappa,\label{Equation: ksigma as a sum of sigmas}
	\end{equation}
	where $\sigma(n)$ is the sum-of-divisors function.
\end{definition}

Note that by \eqref{Equation: LCM identity}, \Cref{Definition: ksigma (n) (revised)} agrees with \Cref{Definition: ksigma (n)} whenever $\kappa = k$ is a positive integer.

The notation $\ksigma\kappa$ is motivated by the fact that for each $n \in \bbZ_{>0}$,
\[
\ksigma \kappa (n) \sim \sigma(n)^{\kappa} \ \text{as} \ \kappa \to \infty.
\]
Moreover, $\ksigma \kappa (n)$ approaches $\sigma(n)^\kappa$ monotonically from below as $\kappa \to \infty$. Indeed, since $\ksigma\kappa(n)$ is multiplicative, it suffices to note that
\[ \lim_{\kappa\to\infty}\frac{\ksigma\kappa(p^\ell)}{\sigma(p^\ell)^\kappa} =  \lim_{\kappa\to\infty}\frac{\sigma(p^\ell)^\kappa-\sigma(p^{\ell-1})^\kappa}{\sigma(p^\ell)^\kappa}=\lim_{\kappa\to\infty}\parent{1-\parent{\frac{\sigma(p^{\ell-1})}{\sigma(p^\ell)}}^\kappa}=1. \]
These relations motivate us to examine other properties of $\ksigma\kappa$ analogous to those of $\sigma^\kappa$. 

We start by estimating the partial sums of $\ksigma \kappa(n)$.

\begin{theorem}\label{Intro Theorem: Mean value for ksigma}
    Let $\kappa>1$. We have
    \[\sum_{n\leq x}\ksigma \kappa (n)= \frac{c(\kappa)}{(\kappa + 1) \zeta(\kappa + 1)} x^{\kappa+1}+O(x^\kappa (\log x)^\kappa)\]
for all $x\ge2$, where $\zeta$ is the Riemann zeta function,
    \begin{equation}\label{Equation: definition of c kappa}
    		c(\kappa) \coloneqq \sum_{n \geq 1} \frac{\sigma_{-1}^{[\kappa]}(n)}{n},
    	\end{equation}
    	and
    \begin{equation}\label{Equation: Definition of ksigma k -1 (n)}
    		\sigma_{-1}^{[\kappa]}(n) \coloneqq \sum_{d \mid n} \mu(n/d) \sigma_{-1}(d)^\kappa
    \end{equation}
    is the $\kappa$th LCM-power of $1/n$.
\end{theorem}

We next derive an upper bound for $\ksigma \kappa(n)$, and thereby prove an analogue of Gr\"{o}nwall's theorem \cite[(25)]{Gronwall}
	\begin{equation} \label{Equation: Gronwall's Theorem}
    	\limsup_{n \to \infty} \frac{\sigma(n)}{e^\gamma n \log \log n} = 1,
  \end{equation}
namely that the maximal order of $\ksigma\kappa$ is
\[\frac{(e^\gamma n \log \log n)^\kappa}{\zeta(\kappa)},\]
where $\gamma$ is the Euler--Mascheroni constant.

\begin{theorem}[$\kappa$-Gr\"{o}nwall's Theorem]\label{Intro Theorem: k-Gronwall's Theorem}
	Let $\kappa > 1$ be a real number. We have
	\[
	\limsup_{n \to \infty} \frac{\zeta(\kappa) \ksigma \kappa(n)}{(e^\gamma n \log \log n)^\kappa} = 1.
	\]
\end{theorem}

Our next result provides an equivalence between the Riemann hypothesis and an elementary inequality for $\ksigma\kappa(n)$, analogous to the Robin criterion. The Robin criterion, which was established in 1984 \cite{Robin}, states that the Riemann hypothesis is equivalent to the inequality 
\begin{equation}\label{Equation: Robin's classical inequality}
\sigma(n) < e^\gamma n \log \log n
\end{equation}
for all $n>5040$. The weaker statement that if the Riemann hypothesis holds then \eqref{Equation: Robin's classical inequality} holds for sufficiently large $n$ was already shown by Ramanjuan in 1915 \cite{Ramanujan1, Ramanujan2}. (See \cite{Broughan} for the history.) We prove the following analogue to the Robin criterion.

\begin{theorem}[$\kappa$-Robin criterion]\label{Intro Theorem: k-Ramanujan-Robin Theorem}
	Let $\kappa > 3/2$ be a real number. The following are equivalent:
	\begin{enumerate}
		\item The Riemann hypothesis holds;\label{Condition: RH holds}
		\item For all $n$ sufficiently large, we have
		\[
		\ksigma \kappa (n) < \frac{\parent{e^\gamma n \log \log n}^\kappa}{\zeta(\kappa)}.
		\]\label{Condition: for n large ksigma (n) is bounded}
	\end{enumerate}
	If $\kappa \geq 2$, we may replace \eqref{Condition: for n large ksigma (n) is bounded} above with the condition
	\begin{enumerate}
		\item[\textnormal{(2$^\prime$)}]
      \makeatletter
      \edef\@currentlabel{2$^\prime$}
      \makeatother
      \label{Condition: for n > 2162160 ksigma(n) is bounded} For all $n > 2162160$, we have 
		\begin{equation}
		\ksigma \kappa (n) < \frac{\parent{e^\gamma n \log \log n}^\kappa}{\zeta(\kappa)}.\label{Equation: k-Robin inequality}
		\end{equation}
	\end{enumerate}
\end{theorem}

\begin{remark}
We obtain the Robin criterion with the larger threshold $n > 2162160$ in place of $n>5040$ by multiplying both sides of \eqref{Equation: k-Robin inequality} by $\zeta(\kappa)$, taking the $\kappa$th roots, and appealing to \Cref{Lemma: counterexamples accrue as kappa shrinks} below. To recover the original version of the Robin criterion, one needs only to verify Robin's inequality \eqref{Equation: Robin's classical inequality} numerically (as Robin himself did) for $5040< n\le 2162160$.
\end{remark}

Let $H_n \coloneqq \sum_{1 \leq m \leq n} 1/m$ denote the $n$th harmonic number. It is well-known that 
\[
H_n = \log n + \gamma +O(1/n)
\]
as $n \to \infty$. 
In 2000, Lagarias provided an alternative formulation of the Robin criterion, establishing the equivalence of the Riemann hypothesis to the inequality
\begin{equation}
\sigma(n) < H_n + e^\gamma e^{H_n} \log H_n\label{Equation: Lagarias criterion}
\end{equation}
for all $n >1$ \cite[Theorem 1.1]{Lagarias}. We prove the following analogue to Lagarias' criterion.

\begin{theorem}[$\kappa$-Lagarias criterion]\label{Intro Theorem: Lagarias criterion}
	Let $\kappa \geq 4$ be a real number. The following are equivalent:
	\begin{enumerate}
		\item The Riemann hypothesis holds;
		\item For all $n > 1$, we have
		\[
		\ksigma \kappa (n) < \frac{\parent{H_n +  e^\gamma e^{H_n} \log H_n}^\kappa}{\zeta(\kappa)}.
		\]
	\end{enumerate}
\end{theorem}

\noindent{\bf Outline and Notation.} The remainder of this paper is organized as follows. In \Cref{Section: Arithmetic estimates}, we recall several classical arithmetic estimates needed for the proofs of our main results. In \Cref{Section: asymptotics}, we determine the mean value and extremal orders of the arithmetic function $\ksigma \kappa (n)$. In \Cref{Section: kappa-colossally abundant numbers}, we define the $\kappa$-colossally abundant numbers in analogy with the colossally abundant numbers of Ramanujan \cite{Ramanujan1, Ramanujan2} and Robin \cite{Robin}, and develop their properties by means of an auxiliary function $\kF \kappa (x, a)$. In \Cref{Section: Robin's Theorem}, we leverage Robin's Theorem (\Cref{Theorem: Classical Robin's Theorem}) to prove $\eqref{Condition: for n large ksigma (n) is bounded} \Rightarrow \eqref{Condition: RH holds}$ in \Cref{Intro Theorem: k-Ramanujan-Robin Theorem}. In \Cref{Section: Ramanujan's Theorem}, we use the theory developed in \Cref{Section: kappa-colossally abundant numbers} to prove $\eqref{Condition: RH holds} \Rightarrow \eqref{Condition: for n large ksigma (n) is bounded}$ and then (when $\kappa \geq 2$) to prove $\eqref{Condition: RH holds} \Rightarrow \eqref{Condition: for n > 2162160 ksigma(n) is bounded}$ in \Cref{Intro Theorem: k-Ramanujan-Robin Theorem}. In \Cref{Section: An analogue to Lagarias criterion}, we establish \Cref{Intro Theorem: Lagarias criterion}. In \Cref{Section: Future work}, we give some possible directions for future work. In our appendix, we give a direct proof of our analogue to Robin's theorem in the case $\kappa = 2$.

For any $x\in\mathbb{R}$, we denote the floor %integer part 
of $x$ by $\lfloor x\rfloor$, which is the largest integer not exceeding $x$, and the ceiling of $x$ by $\lceil x \rceil$, which is the least integer no less than $x$. The letter $p$ always represents a prime number. We write $\pi(x)$ for the prime counting function and $\theta(x) \coloneqq \sum_{p\leq x}\log p$ for the Chebyshev theta function.
We use Landau's big-$O$ notation and Vinogradov's notation $\ll$ interchangeably, and also adopt the standard order notations $o,\gg,\sim$ from analytic number theory.

\section{Arithmetic estimates}\label{Section: Arithmetic estimates}

In this section, we collect some classical arithmetic estimates which we will require later in the paper.

\subsection{Unconditional estimates}

Even without the Riemann hypothesis, we can obtain meaningful bounds for various functions of arithmetic interest.

\begin{lemma}\label{Lemma: Unconditional bound on prod (1 - p^-1)^-1}
	For $x \geq 286$, we have
	\[
	e^\gamma \log x \parent{1 - \frac{1}{2 \parent{\log x}^2}} < \prod_{p \leq x} \parent{1 - p\inv}\inv < e^\gamma \log x \parent{1 + \frac{1}{2 \parent{\log x}^2}}.
	\]
\end{lemma}

%Used to prove \Cref{Theorem: Bound on ksigma k(n)}

\begin{proof}
	Rosser--Schoenfeld \cite[Theorem 8, p. 70]{Rosser-Schoenfeld}. %\cite[Lemma 5.12]{Broughan}.
\end{proof}

\begin{lemma}\label{Lemma: sum of p^-ell}
	Let $\lambda > 1$ and let $x > 1$. We have
	\[
	\sum_{p > x} p^{-\lambda} < \frac{1.01624 \lambda x^{1 - \lambda}}{(\lambda - 1) \log x}.
	\]
\end{lemma}

%Used to prove \Cref{Lemma: Unconditional bound on prod (1 - p^-k)^-1} below.

\begin{proof}
	We follow the proof in Rosser--Schoenfeld \cite[p. 87]{Rosser-Schoenfeld}, using their Theorem 9 without rounding up the value $1.01624$ to $1.02$.
\end{proof}

\begin{lemma}\label{Lemma: Unconditional bound on prod (1 - p^-k)^-1}
	Let $\kappa > 1$ be a real number. For $x > 20000$, we have
	\[
	\frac{1}{\zeta(\kappa)} < \prod_{p \leq x} \parent{1 - p^{-\kappa}} < \frac{1}{\zeta(\kappa)} \exp\parent{\frac{1.01624 \kappa x^{1 - \kappa}}{\log x} \parent{\frac{1}{\kappa - 1}  + 0.000052}}.
	\]
\end{lemma}

%Used to prove \Cref{Theorem: Bound on ksigma k(n)}
%Used to prove Theorem {Theorem: Ineffective Ramanujan's Theorem}

\begin{proof}
	The left inequality is immediate. By \Cref{Lemma: sum of p^-ell}, we have
	\begin{align*}
		\prod_{p \leq x} \parent{1 - p^{-\kappa}} %&= \frac{1}{\zeta(\kappa)} \exp\parent{- \sum_{p > x} \log\parent{1 - p^{-\kappa}}} \\
		&= \frac{1}{\zeta(\kappa)} \exp\parent{\sum_{n \geq 1} \frac 1n \sum_{p > x} p^{-\kappa n}} \\
		&< \frac{1}{\zeta(\kappa)} \exp\parent{\frac{1.01624 \kappa x}{\log x} \sum_{n \geq 1} \frac{x^{- \kappa n}}{\kappa n - 1}} \\
		%&= \frac{1}{\zeta(\kappa)} \exp\parent{\frac{1.01624 \kappa x^{1 - \kappa}}{(\kappa - 1) \log x} + \frac{1.01624 \kappa x}{\log x} \sum_{n \geq 2} \frac{x^{- \kappa n}}{\kappa n - 1}} \\
		%&\leq \frac{1}{\zeta(\kappa)} \exp\parent{\frac{1.01624 \kappa x^{1 - \kappa}}{(\kappa - 1) \log x} + \frac{1.01624 \kappa x^{1-\kappa}}{\log x} \sum_{n \geq 1} \frac{x^{- \kappa n}}{\kappa n + \kappa - 1}} \\
		%&\leq \frac{1}{\zeta(\kappa)} \exp\parent{\frac{1.01624 \kappa x^{1 - \kappa}}{(\kappa - 1) \log x} + \frac{1.01624 x^{1-\kappa}}{\log x} \sum_{n \geq 1} \frac{x^{- \kappa n}}{n}} \\
		%&\leq \frac{1}{\zeta(\kappa)} \exp\parent{\frac{1.01624 \kappa x^{1 - \kappa}}{\log x} \parent{\frac{1}{\kappa - 1} + \sum_{n \geq 1} \frac{x^{- \kappa n}}{n}}} \\
		%&= \frac{1}{\zeta(\kappa)} \exp\parent{\frac{1.01624 \kappa x^{1 - \kappa}}{\log x} \parent{\frac{1}{\kappa - 1} - \log\parent{1 - x^{-\kappa}}}} \\
		&\leq \frac{1}{\zeta(\kappa)} \exp\parent{\frac{1.01624 \kappa x^{1 - \kappa}}{\log x} \parent{\frac{1}{\kappa - 1} - \log\parent{1 - x^{-\kappa}}}} \\
		&< \frac{1}{\zeta(\kappa)} \exp\parent{\frac{1.01624 \kappa x^{1 - \kappa}}{\log x} \parent{\frac{1}{\kappa - 1} + 0.000052}}
		%&= \frac{1}{\zeta(\kappa)} \parent{\frac{e^{1/\parent{\kappa - 1}}}{1 - x^{-\kappa}}}^{1.01624 \kappa x^{1 - \kappa}/\log x}
	\end{align*}
	where the last equality follows because $-\log(1 - x^{-\kappa}) < -\log(1 - 20000^{-1}) < 0.000052$.
\end{proof}

\begin{lemma}\label{Lemma: Unconditional bound on theta(x)}
	For $x \geq 19421$ we have
	\[
	\abs{\theta(x) - x} < \frac{x}{8 \log x}.
	\]
\end{lemma}

%Used to prove Theorem \ref{Theorem: Bound on ksigma k(n)}

\begin{proof}
	Schoenfeld \cite[Corollary 2*, p. 359]{Schoenfeld}.
	%\cite[Theorem 4]{Rosser-Schoenfeld}.
\end{proof}

\begin{lemma}\label{Lemma: bound on (1 - (p+1)^-k)(1 - p^-k)}
	Let $\kappa \geq 3/2$ and $x>1$. We have
	\begin{equation}
	1 \leq \prod_{p > x} \frac{1 - (p+1)^{-\kappa}}{1 - p^{-\kappa}} < \exp\parent{\frac{1.0779 (\kappa + 1)}{\kappa x^{\kappa} \log x}}\label{Equation: bound on (1 - (p+1)^-k)(1 - p^-k)}.
	\end{equation}
\end{lemma}

%Used to prove Theorem \ref{Theorem: Ineffective Ramanujan's Theorem}

\begin{proof}
	The left inequality is immediate. Now recalling Bernoulli's inequality $1+rx\leq (1+x)^r$ for $r\geq 1, x\geq -1$ and that $\log\parent{1 + u} < u$ for $-1<u<0$, we compute
	\begin{align*}
		\prod_{p > x} \frac{1 - (p+1)^{-\kappa}}{1 - p^{-\kappa}} %&= \exp\parent{\sum_{p > x} \log\parent{\frac{1 - (p+1)^{-\kappa}}{1 - p^{-\kappa}}}} \\
		&= \exp\parent{\sum_{p > x} \log\parent{1 + \frac{1 - \parent{1 - \frac{1}{p+1}}^{\kappa}}{p^\kappa - 1}}} \\
		&< \exp\parent{\sum_{p > x} \frac{1 - \parent{1 - \frac{1}{p+1}}^{\kappa}}{p^\kappa - 1}} \\
		&< \exp\parent{\sum_{p > x} \frac{\kappa}{(p + 1) \parent{p^\kappa - 1}}}.
	\end{align*}
	We have
	\[
	\parent{p + 1}\parent{p^\kappa - 1} > \frac{2^{3/2} p^{\kappa + 1}}{3} > \frac{p^{\kappa + 1}}{1.06067},
	\]
	so
	\[
	\prod_{p > x} \frac{1 - (p+1)^{-\kappa}}{1 - p^{-\kappa}} < \exp\parent{\sum_{p > x} \frac{1.06067\kappa}{p^{\kappa + 1}}}.
	\]
	Applying \Cref{Lemma: sum of p^-ell}, we obtain \eqref{Equation: bound on (1 - (p+1)^-k)(1 - p^-k)}.
\end{proof}

\subsection{Conditional estimates}

If the Riemann hypothesis holds, we can obtain stronger bounds in a few cases of interest.

\begin{lemma}\label{Lemma: bound on product (1 - p^-2)}%See Lemma 7.8 of Broughan
	If the Riemann hypothesis holds, then for $x \geq 20000$, we have
	\[
	\prod_{\sqrt{2x} < p \leq x} \parent{1 - p^{-2}} \leq \exp\parent{\frac{-\sqrt{2}}{\sqrt{x} \log x} + \frac{4}{\sqrt{x} \parent{\log x}^2}}.
	\]
\end{lemma}

%Used to prove Theorem \ref{Theorem: Ineffective Ramanujan's Theorem}

\begin{proof}
	Robin \cite[Lemma 6]{Robin}.
\end{proof}

\begin{lemma}\label{Lemma: bound on theta(x)}%Lemma 3.11, Theorem 4.6
	If the Riemann hypothesis holds, then for $x \geq 599$ we have
	\[
	\abs{\theta(x) - x} < \frac{\sqrt{x} \parent{\log x}^2}{8 \pi}.
	\]
	As a consequence, for $x > 0$ we have $\theta(x) < 1.000081x$. For $x \geq 11927$ we have $\theta(x) > 0.985 x$.
\end{lemma}

\begin{proof}
	Schoenfeld \cite[Theorem 10]{Schoenfeld}, Broughan \cite[Lemma 3.11]{Broughan}, and Rosser--Schoenfeld \cite[Corollary to Theorem 6, p. 265]{Rosser-Schoenfeld2}. 
\end{proof}

\begin{remark}
	The bounds $0.985 x < \theta(x)$ and $\theta(x) < 1.000081x$ could be refined using \cite[Proposition 2.1]{LitchmanPomerance}, but we do not require such improvements for our results.
\end{remark}

\begin{lemma}[Nicolas's bound]\label{Lemma: Nicolas's bound}%See Lemma 5.24 and 5.33 of Broughan, or Lemma 5 of Robin, or Nicolas.
	If the Riemann hypothesis holds, then for $x \geq 20000$, we have
	\[
	\prod_{p \leq x} \parent{1 - \frac 1p}\inv \leq e^\gamma \log\theta(x) \exp\parent{\frac{2 + \beta}{\sqrt x \log x} + \frac{\alpha(x)}{\sqrt x \parent{\log x}^2}},
	\]
	where
	\begin{equation}
	\alpha(x) \coloneqq \frac{\parent{\theta(x) - x}^2\parent{1.31 + \log x}}{2 x^{3/2}} + \parent{\beta - 2} + \frac{8 + 4 \beta}{\log x} + \frac{2 \log x}{x^{1/6}} + \frac{(\log 2 \pi) \log x}{x^{1/2}},\label{Equation: definition of alpha(x)}
	\end{equation}
	and $\beta \coloneqq \gamma + 2 - \log 4 \pi$.
\end{lemma}

\begin{proof}
	Robin \cite[Lemma 5]{Robin}. We note in passing that Broughan inadvertently introduces an extra $2$ in front of the term $\log 2 \pi \log x/x^{1/2}$ when he rederives this result \cite[Lemma 7.7]{Broughan}.
\end{proof}

%Used to prove Theorem \ref{Theorem: Ineffective Ramanujan's Theorem}

\section{Asymptotics of \texorpdfstring{$\ksigma \kappa (n)$}{sigma[k](n)}}\label{Section: asymptotics}

In this section, we provide both the mean value and extremal order for the functions $\ksigma \kappa (n)$ when $\kappa > 1$.

\subsection{Mean value of \texorpdfstring{$\ksigma \kappa (n)$}{sigma[k](n)}}\label{Subsection: average order}

We furnish the mean value for the functions $\ksigma \kappa (n)$ when $\kappa > 1$, by means of the following proposition.

\begin{proposition}[{\cite[Corollary 1, p. 66]{Balakrishnan-Petermann}}]\label{Proposition: BP}
For $\kappa>0$ we have
\[\sum_{n\leq x} \sigma^\kappa(n) = \frac{c(\kappa)}{\kappa+1}x^{\kappa+1}+x^\kappa\sum_{r=0}^{\lceil\kappa/3\rceil-1}a_r(\log x)^{\kappa-r}+O\left(x^{\kappa}(\log x)^{2\kappa/3}(\log\log x)^{4\kappa/3}\right) \]
for all $x\ge3$, where $c(\kappa) = \sum_{n \geq 1} \sigma_{-1}^{[\kappa]}(n)/n$ as in \eqref{Equation: definition of c kappa}, and
$a_r=a_r(\kappa)$ are real constants.
\end{proposition}

We now proceed with the main theorem of this subsection.

\begin{theorem}[\Cref{Intro Theorem: Mean value for ksigma}]\label{Theorem: Mean value for ksigma}
    Let $\kappa>1$. We have
    \[\sum_{n\leq x}\sigma^{[\kappa]}(n)=\frac{c(\kappa)}{(\kappa + 1) \zeta(\kappa + 1)} x^{\kappa+1}+x^\kappa\sum_{r=0}^{\lceil\kappa/3\rceil-1}a_r'(\log x)^{\kappa-r}+O\left(x^\kappa(\log x)^{2\kappa/3}(\log\log x)^{4\kappa/3}\right)\]
    for all $x\ge3$, where $a_r'=a_r'(\kappa)$ are real constants.
\end{theorem}

\begin{proof}
We apply Dirichlet's hyperbola method to obtain
    \begin{align*}
    \sum_{n\leq x}\sigma^{[\kappa]}(n)&=\sum_{a\leq y}\mu(a)\sum_{b\leq x/a}\sigma^\kappa(b)+\sum_{b\leq x/y}\sigma^\kappa(b)\sum_{a\leq x/b}\mu(a)-\sum_{a\leq y}\mu(a)\sum_{b\leq x/y}\sigma^\kappa(b),
\end{align*}
with $y \coloneqq x^{1/\kappa}$.

To estimate the first double sum,
\begin{equation}
\sum_{a\leq x^{1/\kappa}}\mu(a)\sum_{b\leq x/a}\sigma^\kappa(b),\label{Equation: the first double sum in mean value result}
\end{equation}
we use \Cref{Proposition: BP} on the inner sum of \eqref{Equation: the first double sum in mean value result}, and find \eqref{Equation: the first double sum in mean value result} equals
\begin{multline*}
    \frac{c(\kappa)}{\kappa+1} \sum_{a\leq x^{1/\kappa}}\mu(a)\left(\frac{x}{a}\right)^{\kappa+1}+\sum_{r=0}^{\lceil\kappa/3\rceil-1}a_r\sum_{a\leq x^{1/\kappa}}\mu(a)\left(\frac{x}{a}\right)^{\kappa}\left(\log\frac{x}{a}\right)^{\kappa-r}+\\
    O\left(\sum_{a\leq x^{1/\kappa}}\left(\frac{x}{a}\right)^{\kappa}\left(\log\frac{x}{a}\right)^{2\kappa/3}\left(\log\log\frac{x}{a}\right)^{4\kappa/3}\right).    
\end{multline*}
The final big-O term may be estimated as being
\begin{align*}
O\left(\sum_{a\leq x^{1/\kappa}}\left(\frac{x}{a}\right)^{\kappa}\left(\log x\right)^{2\kappa/3}\left(\log\log x\right)^{4\kappa/3}\right)
&=O\left(x^\kappa(\log x)^{2\kappa/3}(\log\log x)^{4\kappa/3}\sum_{a\leq x^{1/\kappa}}\frac{1}{a^\kappa}\right)\\
&=O\left(x^\kappa(\log x)^{2\kappa/3}(\log\log x)^{4\kappa/3}\right),
\end{align*}
since $\kappa>1$. Since
\[\sum_{n\leq x}\frac{\mu(n)}{n^{\kappa}}=\frac{1}{\zeta(\kappa)}+O\left(x^{1-\kappa}\right),\]
the main term of \eqref{Equation: the first double sum in mean value result} is
\[\frac{c(\kappa)}{\kappa+1} x^{\kappa+1}\sum_{a\leq x^{1/\kappa}}\frac{\mu(a)}{a^{\kappa+1}}=\frac{c(\kappa)}{(\kappa+1) \zeta(\kappa + 1)} x^{\kappa+1}+O\left(x^{\kappa}\right).\]
To estimate each of the terms
\[a_rx^\kappa\sum_{a\leq x^{1/\kappa}}\frac{\mu(a)}{a^\kappa}\left(\log\frac{x}{a}\right)^{\kappa-r},\]
we use Newton's binomial theorem with error term:
\[\left(\log x - \log a\right)^{\kappa-r} = \sum_{j=0}^J\binom{{\kappa-r}}{j}(-1)^j(\log a)^j(\log x)^{{\kappa-r}-j}+O\left((\log a)^{J+1}(\log x)^{\kappa-r-J-1}\right)\]
where for each $r$ we split the series at $J=J_r$ satisfying $\kappa-r-J_r-1\leq 2\kappa/3 < \kappa-r-J_r$, so that
\[J_r=\lceil\kappa/3\rceil-r-1.\]
Then we have
\[\left(\log\frac{x}{a}\right)^{\kappa-r}=\sum_{j=0}^{J_r}\binom{{\kappa-r}}{j}(-1)^j(\log a)^j(\log x)^{{\kappa-r}-j}+O\left((\log a)^{\lceil\kappa/3\rceil-r}(\log x)^{2\kappa/3}\right),\]
yielding the estimate
\begin{align*}
\sum_{a\leq x^{1/\kappa}}\frac{\mu(a)}{a^\kappa}\left(\log\frac{x}{a}\right)^{\kappa-r}&=\sum_{a\leq x^{1/\kappa}}\frac{\mu(a)}{a^\kappa}\sum_{j=0}^{J_r}\binom{{\kappa-r}}{j}(-1)^j(\log a)^j(\log x)^{{\kappa-r}-j}\\
&+O\left(\sum_{a\leq x^{1/\kappa}}\frac{(\log a)^{\lceil\kappa/3\rceil-r}}{a^\kappa}(\log x)^{2\kappa/3}\right)\\
&=\sum_{a\leq x^{1/\kappa}}\frac{\mu(a)}{a^\kappa}\sum_{j=0}^{J_r}\binom{{\kappa-r}}{j}(-1)^j(\log a)^j(\log x)^{{\kappa-r}-j}+O\left((\log x)^{2\kappa/3}\right).
\end{align*}
Thus we have
\begin{multline*}
a_rx^\kappa\sum_{a\leq x^{1/\kappa}}\frac{\mu(a)}{a^\kappa}\left(\log\frac{x}{a}\right)^{\kappa-r}= \\
a_rx^\kappa\sum_{a\leq x^{1/\kappa}}\frac{\mu(a)}{a^\kappa}\sum_{j=0}^{J_r}\binom{{\kappa-r}}{j}(-1)^j(\log a)^j(\log x)^{{\kappa-r}-j}+O\left(x^\kappa(\log x)^{2\kappa/3}\right).
\end{multline*}

for $0\leq r \leq \lceil\kappa/3\rceil-1$. So far we have that \eqref{Equation: the first double sum in mean value result} is already equal to our desired expression,
\[\frac{c(\kappa)}{(\kappa + 1) \zeta(\kappa + 1)} x^{\kappa+1}+x^\kappa\sum_{r=0}^{\lceil\kappa/3\rceil-1}a_r'(\log x)^{\kappa-r}+O\left(x^\kappa(\log x)^{2\kappa/3}(\log\log x)^{4\kappa/3}\right).\]

Next we show that the two remaining double sums belong to the error term. For the second double sum, we have
\[\sum_{b\leq x^{1-1/\kappa}}\sigma^\kappa(b)\sum_{a\leq x/b}\mu(a)=\sum_{b\leq x^{1-1/\kappa}}\sigma^\kappa(b)\,O\left(\frac{x}{b}\right)=O\left(x\sum_{b\leq x^{1-1/\kappa}}\frac{\sigma^\kappa(b)}{b}\right).\]
Now we use partial summation and by Lemma 1 the estimate
\begin{equation}\label{eq1}\sum_{n\leq x}\sigma^\kappa(n)=O\left(x^{\kappa+1}\right)\end{equation}
to obtain
\begin{align*}
    \sum_{b\leq x^{1-1/\kappa}}\frac{\sigma^\kappa(b)}{b} & =\frac{x^{1/\kappa}}{x}\sum_{b\leq x/x^{1/\kappa}}\sigma^{\kappa}(b)+\int_{1^-}^{x/x^{1/\kappa}}\frac{1}{u^2}\sum_{b\leq u}\sigma^\kappa(b)\,du\\
    & =\frac{x^{1/\kappa}}{x}O\left(\left(\frac{x}{x^{1/\kappa}}\right)^{\kappa+1}\right)+\int_{1^-}^{x/x^{1/\kappa}}\frac{1}{u^2}O\left(u^{\kappa+1}\right)\,du\\
    &= O\left(x^{\kappa-1}\right)
\end{align*}
Thus we have
\[\sum_{b\leq x^{1-1/\kappa}}\sigma^\kappa(b)\sum_{a\leq x/b}\mu(a)=O\left(x^{\kappa}\right).\]
To estimate the final double sum, we use the trivial estimate
\[\sum_{n\leq x}\mu(n)=O(x)\]
and the estimate \eqref{eq1}
so that
\[ \sum_{a\leq x^{1/\kappa}}\mu(a)\sum_{b\leq x^{1-1/\kappa}}\sigma^\kappa(b)=O(x^{1/\kappa}(x/x^{1/\kappa})^{\kappa+1})=O\left(x^{\kappa}\right).\]
This establishes our result.
\end{proof}

\subsection{Mean value of \texorpdfstring{$\ksigma 2 (n)$}{sigma[2](n)}}

In this subsection, we improve on \Cref{Theorem: Mean value for ksigma} in the case $\kappa=2$ by detemining explicit constants. To do so, we require another result from \cite{Balakrishnan-Petermann}.

\begin{proposition}[{\cite[Theorem 1, Theorem 2, and Lemma 3]{Balakrishnan-Petermann}}]\label{Proposition: BP1}

Let $\{a(n)\}_{n=1}^\infty$ be a sequence of real numbers satisfying
\[
\sum_{n=1}^\infty\frac{a(n)}{n^s}=\zeta(s)\zeta^\alpha(s+1)f(s+1)
\]
for $\alpha>0$, and suppose moreover that $f(s)$ has a Dirichlet series expansion which is absolutely convergent in the half plane $\sigma>1-\lambda$ for some $\lambda>0$. Let 
\begin{equation}\label{Equation: definition of v(n)}
\zeta^\alpha(s+1)f(s+1)=\sum_{n=1}^\infty\frac{v(n)}{n^s}.
\end{equation}
Then
\begin{equation}\label{Equation: core equation ffor Balakrishnan-Petermann}
\sum_{n\leq x}a(n) = x\sum_{n\leq x}\frac{v(n)}{n}-\frac{1}{2}\sum_{n\leq x}v(n) + O((\log x)^{2\alpha/3}(\log\log x)^{4\alpha/3})
\end{equation}
for all $x\ge3$.
\end{proposition}

We also require a formulation of the Selberg-Delange method.

\begin{proposition}[{\cite[Theorem II.5.2]{Tenenbaum}}]\label{Proposition: Selberg Delange}
Let $\set{a_n}_{n \geq 1}$ be a positive sequence of real numbers, and suppose $F(s)\coloneqq\sum_{n\geq 1}a_nn^{-s}$ is a Dirichlet series and $\alpha \geq 0$ is a real number such that $\zeta(s)^{-\alpha} F(s)$ may be continued as a holomorphic function of $s = \sigma + i \tau$ for $\sigma \geq 1 - c_0/(1 + \max(0, \log \tau))$, and $\abs{\zeta(s)^{-\alpha} F(s)} \leq M (1 + \abs{\tau})^{1 - \delta}$. For $x \geq 3$, we have
$$ \sum_{n\leq x}a_n=x(\log x)^{\alpha-1}\left(\sum_{0\leq k\leq \alpha-1}\frac{\lambda_k(\alpha)}{(\log x)^k}+O\left(Me^{-c_1\sqrt{\log x}}\right)\right) $$
where $c_1 = c_1(\alpha, M, c_0) > 0$, $\lambda_k(\alpha)\coloneqq \mu_k(\alpha)/\Gamma(\alpha-k)$, and the $\mu_k(\alpha)$ are defined by the Taylor expansion
$$ \frac{s^\alpha F(s+1)}{s+1}=\sum_{k\geq 0}\mu_k(\alpha)s^k. $$
\end{proposition}

The results \cite[Theorem 1, Theorem 2, and Lemma 3]{Balakrishnan-Petermann} and \cite[Theorem II.5.2]{Tenenbaum} are in fact more general than \Cref{Proposition: BP1} and \Cref{Proposition: Selberg Delange}, but these formulations suffice for our purposes.

\begin{lemma}\label{Lemma: sum for n<= x of sigma(n)/n squared}
	We have 
\begin{equation}\label{Equation: partial sum of sigma(n)/n squared}
\sum_{n\leq x}\left(\frac{\sigma(n)}{n}\right)^2=\frac{5}{2}\zeta(3)x-\frac{1}{4}\parent{\log x}^2+O\left((\log x)^{4/3}(\log\log x)^{8/3}\right),
\end{equation}
for all $x\ge3$.
\end{lemma}

\begin{proof}
We use Ramanujan's identity
\[ \sum_{n=1}^\infty\frac{\sigma_a(n)\sigma_b(n)}{n^s}=\frac{\zeta(s)\zeta(s-a)\zeta(s-b)\zeta(s-a-b)}{\zeta(2s-a-b)}, \]
taking $a=b=1$, to get
\begin{equation}
\sum_{n=1}^\infty\frac{\sigma^2(n)}{n^s}=\frac{\zeta(s)\zeta^2(s-1)\zeta(s-2)}{\zeta(2s-2)}.\label{Equation: Ramanujan's identity when a = b = 1}
\end{equation}

Recalling that $\sigma_{-1}(n) = \sigma(n)/n$, and by \eqref{Equation: Definition of ksigma k -1 (n)} and \eqref{Equation: Ramanujan's identity when a = b = 1}, we have
\[
\sum_{n \geq 1} \frac{\sigma_{-1}^{[2]}(n)}{n^s} = \frac{\zeta(s+1)^2 \zeta(s + 2)}{\zeta(2 s + 2)}.
\]
When $s = 1$, this implies
\[ c(2) = \frac{\zeta^2(2)\zeta(3)}{\zeta(4)}=\frac{5}{2}\zeta(3).\]
We apply %Balakrishnan--Petermann 
\Cref{Proposition: BP1} with $a(n) = \left(\frac{\sigma(n)}{n}\right)^2$ and $\alpha = 2$ to obtain 
\[
\sum_{n \leq x} \left(\frac{\sigma(n)}{n}\right)^2 \sim \frac{5}{2}\zeta(3) x
\]
for $x$ sufficiently large.

We now refine this estimate by applying the Selberg--Delange method to the sequence $a_n = n \sigma_{-1}^{[2]}(n)$ with $\alpha = 2$, $c_0 = 1/2$, arbitrary $\delta \in (0, 1)$, and $M = M(\delta)$. Suppose that $x$ is sufficiently large. Here, as before, $\sigma_{-1}^{[2]}(n)$ is defined by \eqref{Equation: Definition of ksigma k -1 (n)}. We thereby determine asymptotics for
\[\sum_{n\leq x}n \sigma_{-1}^{[2]}(n) \]
using
\[ F(s)=\sum_{n=1}^\infty\frac{n\sigma_{-1}^{[2]}(n)}{n^s}=\frac{\zeta(s+1)\zeta^2(s)}{\zeta(2s)}. \]
By \Cref{Proposition: Selberg Delange}, it suffices to determine the leading constant of the power series
\[ \frac{s^2F(s+1)}{s+1}=\frac{\zeta(s+2)s^2\zeta^2(s+1)}{(s+1)\zeta(2s+2)}=\sum_{k\geq 0}\mu_k(2)s^k. \]
Since $s\zeta(s+1)=1 + O(s)$, we have $(s\zeta(s+1))^2=1 + O(s)$. Thus,
\begin{align*}
\mu_0(2)&=\frac{\zeta(2)}{\zeta(2)}=1, \ \text{and} \\
\lambda_0(2) &= \frac{\mu_0(2)}{\Gamma(2 - 0)} = 1.
\end{align*}
We conclude that
\[ \sum_{n\leq x}n\sigma_{-1}^{[2]}(n)=x\log x+O(x). \]
By partial summation, we find
\begin{equation}\label{Equation: Partial sums of sigma -1 kappa (n)}
\sum_{n\leq x}\sigma_{-1}^{[2]}(n)=\frac{1}{2}\parent{\log x}^2 + O(\log x)
\end{equation}
and
\begin{equation}\label{Equation: Partial sums of sigma -1 kappa (n)/n}
\sum_{n\leq x}\frac{\sigma_{-1}^{[2]}(n)}{n}=F(2)-\sum_{n>x}\frac{\sigma_{-1}^{[2]}(n)}{n}=\frac{5}{2}\zeta(3)+O\left(\frac{\log x}{x}\right).
\end{equation}
Substituting \eqref{Equation: Partial sums of sigma -1 kappa (n)} and \eqref{Equation: Partial sums of sigma -1 kappa (n)/n} into \eqref{Equation: core equation ffor Balakrishnan-Petermann} with $v(n) = \sigma_{-1}^{[2]}(n)$, we obtain
\[
\sum_{n\leq x}\left(\frac{\sigma(n)}{n}\right)^2 = \frac{5}{2}\zeta(3)x -\frac{1}{4}\parent{\log x}^2+O\left((\log x)^{4/3}(\log\log x)^{8/3}\right),
\]
as required.
\end{proof}

\begin{theorem}\label{Theorem: Mean value for 2sigma}
We have
\[ \sum_{n\leq x}\sigma^{[2]}(n) = \frac{5}{6}x^3 +\frac{3}{4\zeta(2)}x^2 \parent{\log x}^2 +O\left(x^2 (\log x)^{4/3}(\log\log x)^{8/3}\right) \]
for all $x\ge3$.
\end{theorem}

\begin{proof}
	We apply partial summation to \eqref{Equation: partial sum of sigma(n)/n squared} to get 
\[
\sum_{n\leq x}\sigma^2(n) = \frac{5}{6}\zeta(3)x^3 +\frac{3}{4} x^2 \parent{\log x}^2 + O\left(x^2 (\log x)^{4/3}(\log\log x)^{8/3}\right)
\]
for all $x\ge3$. We then follow the proof of \Cref{Theorem: Mean value for ksigma} for the special case $\kappa=2$ to note that the leading constant must be divided by $\zeta(3)$ and the secondary constant must be divided by $\zeta(2)$ to obtain our result.
\end{proof}

\subsection{An upper bound for \texorpdfstring{$\ksigma \kappa (n)$}{sigma[k](n)}}\label{Section: upper bound}

In this subsection, we provide an upper bound for $\ksigma \kappa (n)$. Although this bound is not sharp, it is a maximal order, meaning that the the limit superior of $\ksigma \kappa (n)$ divided by the upper bound is 1. As a corollary, we obtain \Cref{Intro Theorem: k-Gronwall's Theorem}, our analogue to Gr\"onwall's theorem. Before proving an upper bound for $\ksigma \kappa (n)$, we set some notation.

\begin{definition}\label{Definition: primorial and residual}
	For $x$ a nonnegative real number, we define $\primorial{x}$, the \defi{primorial} of $x$, to be the product of the primes less than or equal to $x$; that is, $\primorial{x} \coloneqq \prod\limits_{p \leq x} p$. For $x$ a nonnegative real number, we define the \defi{primorial residual} $\primeres x$ to be the largest prime $\ell$ such that $\primorial{\ell} \leq x$ if such a prime exists; otherwise, we define $\primeres x = 0$.
\end{definition}

%Note that $\primorial x = e^{\theta(x)}$.

\begin{remark}
The primorial residual takes its name from the following observation: if we restrict the domain of the primorial to the set of primes, then it becomes a residuated mapping, and the primorial residual (with domain restricted to the interval $[2, \infty)$) is its residual \cite[page 11]{BlythJanowitz}. Under \Cref{Definition: primorial and residual}, however, the primorial is merely quasi-residuated \cite[page 9]{BlythJanowitz}.
\end{remark}

Our upper bound for $\ksigma \kappa (n)$ is loosely inspired by \cite[Theorem 2]{Robin} (as corrected by \cite[Theorem 7.13]{Broughan}). Before we prove a global upper bound for $\ksigma \kappa (n)$, we require a local upper bound.

\begin{lemma}\label{Lemma: sigma k(p^ell)/p^k ell < (1-p^-k)/(1 - p^-1)^k}
	Let $\kappa > 1$ be a real number. For $p$ prime and $\ell \geq 1$ an integer, we have
	\[
	\frac{\ksigma \kappa (p^\ell)}{p^{\kappa \ell}} < \frac{1 - p^{-\kappa}}{\parent{1 - p\inv}^{\kappa}}.
	\]
\end{lemma}

%Used to prove Theorem \ref{Theorem: Bound on ksigma k(n)}

\begin{proof}
For $\ell \geq 1$ arbitrary, we compute
	\begin{align}
	\frac{\ksigma \kappa (p^\ell)}{p^{\kappa \ell}} %&= \frac{\parent{p^{\ell+1} - 1}^\kappa - \parent{p^\ell -1}^\kappa}{p^{\kappa \ell} (p-1)^\kappa} \nonumber \\
	%&= \frac{\parent{1 - p^{-\ell-1}}^\kappa}{\parent{1 - p\inv}^\kappa} \parent{1 - \parent{\frac{1 - p^{-\ell}}{p - p^{-\ell}}}^\kappa} \nonumber \\	
\label{Equation: ksigma(p^ell)/p^k ell}	&= \frac{\parent{1 - p^{-\ell-1}}^\kappa -  \parent{p\inv - p^{-\ell - 1}}^\kappa}{\parent{1 - p\inv}^\kappa}. 
	\end{align}
	The derivative of $\eqref{Equation: ksigma(p^ell)/p^k ell}$ with respect to $\ell$ is
	\[
	\frac{\parent{1 - p^{-\ell-1}}^{\kappa-1} -  \parent{p\inv - p^{-\ell - 1}}^{\kappa-1}}{\parent{1 - p\inv}^\kappa} \kappa p^{-\ell-1} \log p,
	\]
	which is positive, so $\frac{\ksigma \kappa (p^\ell)}{p^{\kappa \ell}}$ is strictly increasing in $\ell$. Thus
	\[
	\frac{\ksigma \kappa (p^\ell)}{p^{\kappa \ell}} < \lim_{\ell \to \infty} \frac{\ksigma \kappa (p^\ell)}{p^{\kappa \ell}} = \frac{1 - p^{-\kappa}}{\parent{1 - p\inv}^{\kappa}}%.\label{Equation: inequality for ksigma(p^ell)/p^kell}
	\]
	as desired.
\end{proof}

\begin{thm}\label{TheoremV: Bound on ksigma k(n)}
For any $\kappa>1$ and $n\ge e^{19183}$, we have
\begin{equation}\label{eq:thmVBound on ksigma k(n)1}
\sigma^{[\kappa]}(n) < \left(e^{\gamma}n\log\log n + \frac{0.42n}{\log\log n}\right)^{\kappa}\prod_{p|n}\left(1-p^{-\kappa}\right)
\end{equation}
and
\begin{align}\label{eq:thmVBound on ksigma k(n)2}
\begin{split}
	\sigma^{[\kappa]}(n) &< \frac{1}{\zeta(\kappa)}\left(e^{\gamma}n\log\log n + \frac{0.42n}{\log\log n}\right)^{\kappa} \exp\Biggl(\frac{1.01624 \kappa (\log n)^{1 - \kappa}}{\log\log n} \parent{\frac{1}{\kappa - 1}  + 0.000052}\cdot \\
		&\phantom{< \frac{1}{\zeta(\kappa)}\left(e^{\gamma}n\log\log n + \frac{0.42n}{\log\log n}\right)^{\kappa} }\left(1-\frac{0.14}{\log\log n}\right)^{1-\kappa}\left(1+\frac{1}{7(\log\log n)^2}\right)\Biggl).
	\end{split}
\end{align}
\end{thm}

\begin{proof}
By \Cref{Lemma: sigma k(p^ell)/p^k ell < (1-p^-k)/(1 - p^-1)^k}, we have
\begin{equation}\label{eq:thmVBound on ksigma k(n)6}
\sigma^{[\kappa]}(n) < n^k \prod_{p|n}\left(1-p^{-\kappa}\right)\left(1-p^{-1}\right)^{-\kappa}.
\end{equation}
Now suppose $x=\primeres n$ so $x$ is prime, and let $x_+$ be the next prime so that
$$ \primorial{x}\leq n <\primorial{x_+}. $$
We claim that $\omega(n)\leq \pi(x)$. For if $\omega(n)\geq \pi(x_+)$, each distinct prime factor of $n$ would be at least as big as each of the primes $\leq x_+$, so $n\geq \primorial{x_+}$, a contradiction.

To prove \eqref{eq:thmVBound on ksigma k(n)1}, we use
\[\prod_{p|n}\left(1-p^{-1}\right)^{-1}\leq \prod_{p\le x}\left(1-p^{-1}\right)^{-1},\]
which is true since each factor is greater than 1,  the number of factors on the left is at most the number of factors on the right, and each $p\mid n$ is either $\leq x$, in which case we have matching factors, or for each factor with $p>x$ on the left we can find an unmatched factor on the right that is greater.

From \cite[Theorem 5.9]{Dusart} we have
\begin{equation}\label{eq:thmVBound on ksigma k(n)3}
\prod_{p\le x}\left(1-p^{-1}\right)^{-1}< e^{\gamma}\log x\left(1+\frac{0.10836}{\parent{\log x}^2}\right)
\end{equation}
for $x \ge 2278382$, and straightforward numerical computation shows that this inequality in fact holds for $x\geq 19421$.
We observe that for the prime $19421$ we have $\theta(19421)=19182.3\dots$, so that $n\geq e^{19183}>e^{\theta(19421)}$. Since $x = \primeres n$, we have $n\geq \primorial{x}$ and so $\log n \ge \theta(x)$. Then
\Cref{Lemma: Unconditional bound on theta(x)} gives us
	\[
	\log \log n \geq \log \theta(x) \geq \log x + \log\parent{1 - \frac{1}{8 \log x}} > \log x - \frac 1{7.899 \log x},
	\]
so that our result would follow if we have
\[e^{\gamma}\log x\left(1+\frac{0.10836}{\parent{\log x}^2}\right)< e^{\gamma}\left(\log x - \frac{1}{7.899\log x}\right) +\frac{0.42}{(\log x - 1/(7.899\log x))},\]
which is clearly true by the stronger inequality
\begin{equation}\label{Equation: thmVbound stronger inequality}
e^{\gamma}\log x\left(1+\frac{0.10836}{\parent{\log x}^2}\right) < e^{\gamma}\left(\log x - \frac{1}{7.899\log x}\right)+\frac{0.42}{\log x}.
\end{equation}
This proves \eqref{eq:thmVBound on ksigma k(n)1}.

To prove \eqref{eq:thmVBound on ksigma k(n)2}, we use the inequality $\left(1-p^{-\kappa}\right)\left(1-p^{-1}\right)^{-\kappa}>1$ and argue as before that since the expression on the left decreases with increasing $p$, we have
$$  \prod_{p\mid n}\left(1-p^{-\kappa}\right)\left(1-p^{-1}\right)^{-\kappa} \leq  \prod_{p\le x}\left(1-p^{-\kappa}\right)\left(1-p^{-1}\right)^{-\kappa}. $$
Thus,
\begin{equation}\label{eq:thmVBound on ksigma k(n)7}
\sigma^{[\kappa]}(n) < n^k \prod_{p\le x}\left(1-p^{-\kappa}\right)\left(1-p^{-1}\right)^{-\kappa}<\left(e^{\gamma}n\log\log n + \frac{0.42n}{\log\log n}\right)^{\kappa}\prod_{p\le x}\left(1-p^{-\kappa}\right).
\end{equation}
By \Cref{Lemma: Unconditional bound on prod (1 - p^-k)^-1}, we have
\begin{equation}\label{eq:thmVBound on ksigma k(n)8}
\prod_{p\le x}\left(1-p^{-\kappa}\right)<\frac{1}{\zeta(\kappa)} \exp\parent{\frac{1.01624 \kappa x^{1 - \kappa}}{\log x} \parent{\frac{1}{\kappa - 1}  + 0.000052}}.
\end{equation}
Since $n\leq \primorial{x_+}$, we have $\log n\leq \theta(x)+\log x_+$. We seek an upper bound for $\log x_+$ in terms of $x$. By \Cref{Lemma: Unconditional bound on theta(x)}, we have
$$ x_+-\frac{x_+}{8\log x_+}-\log x_+\leq \theta(x_+)-\log x_+=\theta(x) \leq x+\frac{x}{8\log x}. $$
Then
$$ \log x_++\log\left(1-\frac{1}{8\log x_+}-\frac{\log x_+}{x_+}\right)<\log x +\log\left(1+\frac{1}{8\log x}\right). $$
Since $x\mapsto 1/(8\log x)+(\log x)/x$ is decreasing and $x<x_+$, we have
$$ \log x_++\log\left(1-\frac{1}{8\log x}-\frac{\log x}{x}\right)<\log x +\log\left(1+\frac{1}{8\log x}\right). $$
Thus we have the bound
\begin{align*}
\log x_+&< \log x +\log\left(1+\frac{1}{8\log x}\right)-\log\left(1-\frac{1}{8\log x}-\frac{\log x}{x}\right)\\
&=\log x +\log\left(1+\frac{x+4\log^2x}{4x\log x-x/2-4\log^2x}\right).
\end{align*}
For the remainder of this argument we will resort to numerics where we have monotonicity. For instance, it is easy to show that for $x\geq 19421$ we have
$$ \log x_+<\log x+\frac{0.26}{\log x}. $$

 Invoking \Cref{Lemma: Unconditional bound on theta(x)} again, we obtain
 $$ \log n < x+\frac{x}{8\log x}+\log x +\frac{0.26}{\log x}, $$
 and thus
\[\log\log n<\log x+\log\left(1+\frac{1}{8\log x}+\frac{\log x}{x}+\frac{0.26}{x\log x}\right).\]
This time we have for $x\geq 19421$ that
$$ \log\log n < \log x +\frac{0.13}{\log x}. $$
Solving for $\log x$, we have
\begin{align*}
\log x&>\frac{\log\log n+\sqrt{(\log\log n)^2-0.52}}{2}\\
&>\log\log n-\frac{0.14}{\log\log n}\\
&>\log\log n\left(1+\frac{1}{7(\log\log n)^2}\right)^{-1}
\end{align*}
and that 
\[x>\exp\left(-\frac{0.14}{\log\log n}\right)\log n>\left(1-\frac{0.14}{\log\log n}\right)\log n.\]
Hence, the expression inside the exponential in \eqref{eq:thmVBound on ksigma k(n)8} is less than
\begin{align*}
\frac{1.01624 \kappa (\log n)^{1 - \kappa}}{\log\log n} \parent{\frac{1}{\kappa - 1}  + 0.000052}\left(1-\frac{0.14}{\log\log n}\right)^{1-\kappa}\left(1+\frac{1}{7(\log\log n)^2}\right).
\end{align*}
Combining this with \eqref{eq:thmVBound on ksigma k(n)7} and \eqref{eq:thmVBound on ksigma k(n)8} completes the proof of \eqref{eq:thmVBound on ksigma k(n)2}.
\end{proof}

If $\kappa \geq 2$, a computation lets us strengthen \Cref{TheoremV: Bound on ksigma k(n)} (see \Cref{Theorem: Effective Ramanujan's Theorem} and the remark thereafter).

\begin{corollary}\label{Corollary: strengthened bound on ksigma k(n)}
	Let $\kappa \geq 2$ be a real number. For each integer $n > 2162160$, the inequalities \eqref{eq:thmVBound on ksigma k(n)1} and \eqref{eq:thmVBound on ksigma k(n)2} hold.
\end{corollary}

\Cref{TheoremV: Bound on ksigma k(n)} also yields a family of analogues to Gr\"{o}nwall's Theorem.

\begin{cor}\label{Corollary: Analogues to Gronwall's Theorem}
	Let $\kappa > 1$ be a real number. We have
	\[
		\limsup_{n \to \infty} \frac{\zeta(\kappa) \ksigma \kappa (n)}{\parent{e^\gamma n \log \log n}^\kappa} = 1.
	\]
\end{cor}

\begin{proof}
	The inequality \eqref{eq:thmVBound on ksigma k(n)2} furnished by \Cref{TheoremV: Bound on ksigma k(n)} shows that
	\[
	\limsup_{n \to \infty} \frac{\zeta(\kappa) \ksigma \kappa (n)}{\parent{e^\gamma n \log \log n}^\kappa} \leq 1,
	\]
	so it suffices to show this bound is asymptotically obtained. We let $a(n) \coloneqq {\primorial n}^{t(n)}$, where $(t(n))_{n \geq 1}$ is a sequence of nonnegative integers such that
	\[
	t(n) \to \infty \ \text{but} \ \log t(n) = o(\log n) \ \text{as} \ n \to \infty;
	\]
	for instance, we could take $t(n) = \floor{\log n}$. Taking logarithms twice, we see
	\[
		\log \log a(n) = \log \log \primorial{n} + \log {t(n)} \sim \log \log \primorial{n}.
	\]
	Now by \Cref{Lemma: Unconditional bound on prod (1 - p^-1)^-1} and \Cref{Lemma: Unconditional bound on prod (1 - p^-k)^-1}, we compute
	\begin{align*}
		\ksigma \kappa (a(n)) &= a(n)^\kappa \prod_{p \leq n} \frac{1 - p^{-\kappa}}{\parent{1-p^{-1}}^\kappa} \cdot \frac{\parent{1 - p^{-t(n)-1}}^\kappa - \parent{p\inv - p^{-t(n)-1}}^\kappa}{1 - p^{-\kappa}} \\
		&\sim \frac{\parent{e^\gamma a(n) \log \log a(n)}^\kappa}{\zeta(\kappa)} \prod_{p \leq n} \frac{\parent{1 - p^{-t(n)-1}}^\kappa - \parent{p\inv - p^{-t(n)-1}}^\kappa}{1 - p^{-\kappa}},
	\end{align*}
	where the product is over primes $p$.
	
	For $\kappa > 0$, $p \geq 2$, and $t \geq 0$, the function
	\[
	t \mapsto \frac{\parent{1 - p^{-t}}^\kappa - \parent{p\inv - p^{-t}}^\kappa}{1 - p^{-\kappa}}
	\]
	is clearly positive and increasing, with a limit of 1, so
	\[
	\prod_{p \leq n} \frac{\parent{1 - p^{-t(n)-1}}^\kappa - \parent{p\inv - p^{-t(n)-1}}^\kappa}{1 - p^{-\kappa}} \leq 1.
	\]
	On the other hand, for $n$ sufficiently large we have
	\begin{align*}
		\prod_{p \leq n} \frac{\parent{1 - p^{-t(n)-1}}^\kappa - \parent{p\inv - p^{-t(n)-1}}^\kappa}{1 - p^{-\kappa}} &> \prod_{p} \frac{\parent{1 - p^{-t(n)-1}}^\kappa - \parent{p\inv - p^{-t(n)-1}}^\kappa}{1 - p^{-\kappa}} \\
		&> \prod_{p} \frac{\parent{1 - p^{-{t(n)}-1}}^\kappa}{1 - p^{-\kappa}} \parent{1 - p^{-\kappa} \parent{1 + p^{-{t(n)}}}^\kappa} \\
		&= \frac{\zeta(\kappa)}{\zeta(t(n) + 1)} \prod_{p} \parent{1 - p^{-\kappa} \parent{1 + p^{-{t(n)}}}^\kappa},
	\end{align*}
	and
	\[
	\lim_{n \to \infty} \prod_{p} \parent{1 - p^{-\kappa} \parent{1 + p^{-{t(n)}}}^\kappa} = \frac{1}{\zeta(\kappa)},
	\]
	so 
	\[
	\limsup_{n \to \infty} \frac{\zeta(\kappa) \ksigma \kappa (a(n))}{\parent{e^\gamma a(n) \log \log a(n)}^\kappa} = 1
	\]
	as desired.
\end{proof}

Robin proved the following unconditional upper bound for $\sigma(n)$.

\begin{theorem}[{\cite[Theorem 2]{Robin}, \cite[Theorem 7.13]{Broughan}}]\label{Theorem: Robin's unconditional bound}
	For $n \geq 3$, we have
	\[
	\sigma(n) < e^\gamma n \log \log n + \frac{2  n}{3 \log \log n}.
	\]
\end{theorem}

We show that we can recover \Cref{Theorem: Robin's unconditional bound} from \Cref{TheoremV: Bound on ksigma k(n)}.

\begin{proof}[Proof of \Cref{Theorem: Robin's unconditional bound}]
	Taking $\kappa$th roots of both sides of \eqref{eq:thmVBound on ksigma k(n)1} and leting $\kappa \to \infty$ proves a stronger version of this inequality for $n > e^{19183}$. The result now follows from a short computation.
\end{proof}

We can derive Gr\"{o}nwall's Theorem itself \cite[(25)]{Gronwall} from \Cref{TheoremV: Bound on ksigma k(n)} as well. Indeed, taking $\kappa$th roots of both sides of the inequality in \Cref{TheoremV: Bound on ksigma k(n)} and letting $\kappa \to \infty$ gives us 
\[
	\limsup_{n \to \infty} \frac{\sigma(n)}{e^\gamma n \log \log n} \leq 1,
\]
and the converse inequality may be obtained by considering $\sigma(a_n)$ as in the proof of \Cref{Corollary: Analogues to Gronwall's Theorem}.

\begin{cor}[Gr\"{o}nwall's Theorem]
	We have
	\[
	\limsup_{n \to \infty} \frac{\sigma(n)}{e^\gamma n \log \log n} = 1.
	\]
\end{cor}

% As written, \Cref{Theorem: Bound on ksigma k(n)} only applies when $n$ is quite large, but for $\kappa$ fixed we can find a cleaner bound that works for relatively small $n$ by performing a computation on $\primorial n$ for $n \leq 19421$, and by replacing the bound given by \Cref{Theorem: Bound on ksigma k(n)} by something weaker and cleaner. \gm{Do we want to insert a table of examples of this here?}

\subsection{A lower bound for \texorpdfstring{$\ksigma \kappa (n)$}{sigma[k](n)}}\label{Section: Lower bound}

In this subsection, we provide an elementary but sharp lower bound for $\ksigma \kappa (n)$. This bound is also a minimal order, meaning that the limit inferior of $\ksigma \kappa (n)$ divided by the lower bound is 1.

\begin{proposition}\label{Proposition: Lower bound for ksigma(n)}
	Let $\kappa > 1$ be a real number. For each integer $n > 1$ we have
	\[
	\ksigma \kappa (n) \geq \parent{n + 1}^\kappa - 1,
	\]
	with equality if and only if $n$ is prime.
\end{proposition}

\begin{proof}
  Observe that for any $x>0$, the function $f_x : [0, \infty) \to [0, \infty)$ defined by $f_x : t \mapsto (x+t)^\kappa-t^\kappa$ has derivative
  \[f_x'(t)=\kappa\left[(x+t)^{\kappa-1}-t^{\kappa-1}\right]>0\]
so $f_x$ is strictly increasing in $t\geq 0$.
Taking $x=p^\ell$ and comparing the values of $f_x$ at $t=\sigma(p^{\ell-1})$ and $t=1$, we have
	\[
	\ksigma \kappa(p^\ell) = \parent{p^\ell + \sigma(p^{\ell-1})}^\kappa - \sigma(p^{\ell-1})^\kappa \geq \parent{p^\ell + 1}^\kappa - 1,
	\]
	with equality if and only if $\ell = 1$. This establishes the proposition for $n=p^\ell$.

	Now we show that
	\[
	\left((m + 1)^\kappa -1\right)\left((n + 1)^\kappa -1\right) > (mn + 1)^\kappa -1
	\]
  for $\kappa>1$ and $m,n\geq 1$. We again make use of $f_x$. Taking $x=m+n$ and comparing the values of $f_x$ at $t=mn+1$ and $t=1$, we have
\[(m + n + mn + 1)^\kappa  - (mn + 1)^\kappa  > (m + n + 1)^\kappa  - 1.\]
Next, taking $x=m$ and comparing the values of $f_x$ at $t=n+1$ and $t=1$, we obtain
\[(m + n + 1)^\kappa  - (n + 1)^\kappa  > (m + 1)^\kappa  - 1.\] 
Adding up these two inequalities and rearranging the terms yield the desired inequality.

\Cref{Proposition: Lower bound for ksigma(n)} now follows by the multiplicativity of $\ksigma \kappa (n)$ and an easy inductive argument.
\end{proof}

\section{\texorpdfstring{$\kappa$}{k}-colossally abundant numbers}\label{Section: kappa-colossally abundant numbers}

In this section, we develop the theory of $\kappa$-colossally abundant numbers in analogy with the classical theory of colossally abundant numbers. The following material is inspired by \cite[Chapter 6]{Broughan}. 

\begin{definition}\label{Definition: k-colossally abundant numbers}
	Let $\epsilon > 0$ be a real number, and define 
	\begin{equation}
	\krho \kappa \epsilon (n) \coloneqq \frac{\ksigma \kappa (n)}{n^{\kappa(1 + \epsilon)}}.\label{Equaton: Defining rho kappa}
	\end{equation}
	We say that a positive integer $N$ is \defi{$\kappa$-colossally abundant for $\epsilon$} if we have
	\begin{equation}
	\krho \kappa \epsilon (N) \geq \krho \kappa \epsilon (n)\label{Equation: Defining inequality for k-colossally abundant numbers}
	\end{equation}
	for all positive integers $n$. If $N$ is $\kappa$-colossally abundant for $\epsilon$ for some $\epsilon > 0$, we say that $N$ is \defi{$\kappa$-colossally abundant}.
\end{definition}

Under this definition, the usual colossally abundant numbers should be thought of as ``$\infty$-colossally abundant'', because $N$ is colossally abundant if and only if for some $\epsilon > 0$ we have
\[
	\frac{\sigma(N)}{N^{1 + \epsilon}} = \lim_{\kappa \to \infty} \sqrt[\kappa]{\frac{\ksigma \kappa (N)}{N^{\kappa(1 + \epsilon)}}} \geq \lim_{\kappa \to \infty} \sqrt[\kappa]{\frac{\ksigma \kappa (n)}{n^{\kappa(1 + \epsilon)}}} = \frac{\sigma(n)}{n^{1 + \epsilon}}
\]
for all positive integers $n$.

We require the following definition.

\begin{definition}\label{Definition: kF and kfp}
	If $\kappa > 1$ and $\epsilon \in \bbR$, and $x > 1$ are given, then for $a \geq 1$ we define
	\[
	\kfp \kappa \epsilon (x, a) \coloneqq x^{-\kappa(1+\epsilon)} \frac{(x^{a+1}-1)^\kappa-(x^{a}-1)^\kappa}{(x^{a}-1)^\kappa-(x^{a-1}-1)^\kappa},
	\]
  and $f^{[\kappa]}(x,a) \coloneqq \kfp \kappa 0 (x,a).$
	We also define $\kF \kappa : (1, \infty) \times [1, \infty) \to (0, \infty)$ to be
	\[
		\kF \kappa (x, a) \coloneqq \frac{1}{\kappa \log x} \log f^{[\kappa]} (x, a) =  \frac{1}{\kappa \log x}\log\parent{\frac{(x^{a+1}-1)^\kappa-(x^{a}-1)^\kappa}{(x^{a+1}-x)^\kappa-(x^{a}-x)^\kappa}}.
	\]
	We adopt the convention $\kF \kappa(x, 0) \coloneqq \infty$ for $x > 1$, but do not consider $(x, 0)$ to be properly within the domain of $\kF \kappa$.
\end{definition}

By construction, we have
\begin{align}
\kfp \kappa \epsilon (p, a) &= \frac{\ksigma \kappa (p^a)}{p^{\kappa (1 + \epsilon)} \ksigma \kappa (p^{a-1})} = \frac{\krho \kappa \epsilon (p^a)}{\krho \kappa \epsilon (p^{a-1})}, \ \text{and} \label{Equation: kfp versus ksigma} \\
\kF \kappa (p, a) &= \frac{1}{\kappa \log p} \log \frac{\ksigma \kappa (p^a)}{p^{\kappa} \ksigma \kappa (p^{a-1})} = \frac{1}{\kappa \log p}  \log \frac{\krho \kappa \epsilon (p^a)}{\krho \kappa \epsilon (p^{a-1})} - \epsilon \label{Equation: kF versus ksigma} 
\end{align}
for $p$ prime and $a \in \bbZ_{> 0}$.

Equations \eqref{Equation: kfp versus ksigma} and \eqref{Equation: kF versus ksigma} suggest that we may extract information about the $\kappa$-colossally abundant numbers by understanding $\kfp \kappa \epsilon$ and $\kF \kappa$. We spend the next subsection developing our understanding $\kF \kappa$.

\subsection{The function \texorpdfstring{$\kF \kappa (x, a)$}{F[k](x, a)}}

In this subsection, we develop our understanding of the function $\kF \kappa (x, a)$, given in \Cref{Definition: kF and kfp}. We demonstrate that $\kF \kappa (x, a)$ is monotonic in its arguments, %compute an asymptotic for it, 
and deduce some information about its partial inverses.

\begin{theorem}\label{Theorem: behavior of kF}
	The function $\kF \kappa (x, a)$ is continuous and strictly decreasing in $x$ and $a$, and continuous and strictly increasing in $\kappa$. Moreover, for $a \geq 1$ and $x > 1$ respectively, we have
	\[
	\lim_{x \to \infty} \kF \kappa (x, a) = 0 \ \text{and} \ \lim_{a \to \infty} \kF \kappa (x, a) = 0,
	\]
	as well as
	\[
	\lim_{x \to 1^+} \kF \kappa(x, a) = \infty.
	\]
	We also have
	\begin{equation}
	\lim_{\kappa \to \infty} \kF \kappa (x, a) = \frac{1}{\log x} \log\parent{\frac{x^{a+1}-1}{x^{a+1}-x}} \eqqcolon F(x, a)\label{Equation: Definition of F(x, a)}
	\end{equation}
	and
	\[
	\lim_{\kappa \to 1^+} \kF \kappa (x, a) = 0.
	\]
\end{theorem}
We remark that $F(x, a)$, defined in \eqref{Equation: Definition of F(x, a)} plays an important role in the study of colossally abundant numbers \cite{ErdosNicolas}. The monotonicity of $F$ in its arguments is essential to that program.

\begin{proposition}\label{Proposition: F is decreasing in x and a}
	Let $F(x, a)$ be as in \eqref{Equation: Definition of F(x, a)}. Then $F$ is decreasing in $x$ and $a$ on $(1, \infty) \times (1, \infty)$.
\end{proposition}

\begin{proof}
Write 
\[F(x,a)=\frac{1}{\log x} \log\parent{1+\frac{1}{xG(x,a)}},\]
where 
\[
G(x,a)\coloneqq\frac{x^{a}-1}{x-1}.
\]

Note that
\[\frac{\partial G}{\partial x}(x,a)=\frac{ax^{a-1}(x-1)-(x^a-1)}{(x-1)^2}.\]
Applying Lagrange's mean value theorem to the function $t \mapsto t^a - 1$, we see that for some $y \in (1, x)$ we have
\[
x^a-1=ay^{a-1}(x-1) < ax^{a-1}(x-1),
\]
and therefore $\frac{\partial G}{\partial x} > 0$ for $x>1$. This shows that $G(x,a)$ is increasing in $x>1$, whence $F(x,a)$ is decreasing in $x>1$.

It is immediate that $G(x, a)$ is increasing in $a$, whence $F(x, a)$ is decreasing in $a > 1$ as well.
\end{proof}

\begin{remark}
The monontonicity of $F(x,a)$ was known to Robin \cite{Robin}, at least when $a$ is an integer.
\end{remark}

Analogously, our function $\kF \kappa (x, a)$ plays an important role in our study of $\kappa$-colossally abundant numbers. So it is not surprising that we need to investigate the monontonicity of $\kF \kappa (x, a)$, but this time in each of the three variables $x,a,\kappa$. Despite the simplicity of the statement of \Cref{Theorem: behavior of kF}, its proof is rather involved. We start by proving that for any fixed $a\ge1$ and $\kappa>1$, $\kF \kappa (x, a)$ is strictly decreasing in $x \in (1,\infty)$. It suffices to show 
\begin{equation}\label{Equ:f^k(x,a)}
f^{[\kappa]}(x,a)=\frac{(x^{a+1}-1)^{\kappa}-(x^{a}-1)^{\kappa}}{(x^{a+1}-x)^{\kappa}-(x^{a}-x)^{\kappa}}
\end{equation}
is strictly decreasing in $x\in(1,\infty)$.
\begin{prop}\label{prop:f^k(x a) in x}
	Let $a\geq1$ and $\kappa>1$ be positive real numbers. Then 
	$f^{[\kappa]}(x,a)$ is strictly decreasing in $x\in(1,\infty)$.
\end{prop}
\begin{proof}
Fix $a\geq1$ and $\kappa>1$, and let $q^{[\kappa]}(x,a) \coloneqq (x^a-1)^{\kappa-1}$ for $x>1$. Then 
\begin{align}
	(x^{a+1}-1)^{\kappa}-(x^{a}-1)^{\kappa}&=(x^{a+1}-1)q^{[\kappa]}(x,a+1)-(x^{a}-1)q^{[\kappa]}(x,a), \label{Equation: technical term 1} \\
	(x^{a+1}-x)^{\kappa}-(x^{a}-x)^{\kappa}&=x^k\left[(x^{a}-1)q^{[\kappa]}(x,a)-(x^{a-1}-1)q^{[\kappa]}(x,a-1)\right].\label{Equation: technical term 2} 
\end{align} 
We compute the derivatives of \eqref{Equation: technical term 1} and \eqref{Equation: technical term 2}:
\begin{align*}
	\frac{\partial}{\partial x}\left[(x^{a+1}-1)^{\kappa}-(x^{a}-1)^{\kappa}\right]&=\kappa x^{-1}\left[(a+1)x^{a+1}q^{[\kappa]}(x,a+1)-ax^{a}q^{[\kappa]}(x,a)\right],\\
	\frac{\partial}{\partial x}\left[(x^{a+1}-x)^{\kappa}-(x^{a}-x)^{\kappa}\right]&=\kappa x^{\kappa-1}\left\{[(a+1)x^a-1]q^{[\kappa]}(x,a)-(ax^{a-1}-1)q^{[\kappa]}(x,a-1)\right\}.
\end{align*}
Thus we have
\begin{align*}
	L^{[\kappa]}(x,a)&\coloneqq\frac{\partial}{\partial x}\left[(x^{a+1}-1)^{\kappa}-(x^{a}-1)^{\kappa}\right]\cdot\left[(x^{a+1}-x)^{\kappa}-(x^{a}-x)^{\kappa}\right]\\
	&=\kappa x^{\kappa-1}\left[(a+1)x^{a+1}q^{[\kappa]}(x,a+1)-ax^{a}q^{[\kappa]}(x,a)\right]\cdot\big[(x^{a}-1)q^{[\kappa]}(x,a)-\\
	&~(x^{a-1}-1)q^{[\kappa]}(x,a-1)\big]\\
	&=\kappa x^{\kappa-1}\big[(a+1)(x^{2a+1}-x^{a+1})q^{[\kappa]}(x,a+1)q^{[\kappa]}(x,a)-(a+1)(x^{2a}-x^{a+1})\cdot\\
	&~q^{[\kappa]}(x,a+1)q^{[\kappa]}(x,a-1)-a(x^{2a}-x^a)q^{[\kappa]}(x,a)^2+a(x^{2a-1}-x^a)\cdot\\
	&~q^{[\kappa]}(x,a)q^{[\kappa]}(x,a-1)\big]
\end{align*}
and
\begin{align*}
	R^{[\kappa]}(x,a)&\coloneqq\left[(x^{a+1}-1)^{\kappa}-(x^{a}-1)^{\kappa}\right]\cdot\frac{\partial}{\partial x}\left[(x^{a+1}-x)^{\kappa}-(x^{a}-x)^{\kappa}\right]\\
	&=\kappa x^{\kappa-1}\left[(x^{a+1}-1)q^{[\kappa]}(x,a+1)-(x^{a}-1)q^{[\kappa]}(x,a)\right]\cdot\big\{[(a+1)x^a-1]q^{[\kappa]}(x,a)\\
	&~-(ax^{a-1}-1)q^{[\kappa]}(x,a-1)\big\}\\
	&=\kappa x^{\kappa-1}\big\{[(a+1)x^{2a+1}-x^{a+1}-(a+1)x^a+1]q^{[\kappa]}(x,a+1)q^{[\kappa]}(x,a)\\
	&~-(ax^{2a}-x^{a+1}-ax^{a-1}+1)q^{[\kappa]}(x,a+1)q^{[\kappa]}(x,a-1)-[(a+1)x^{2a}-\\
	&~(a+2)x^a+1]q^{[\kappa]}(x,a)^2+(ax^{2a-1}-x^a-ax^{a-1}+1)q^{[\kappa]}(x,a)q^{[\kappa]}(x,a-1)\big\}.
\end{align*}
To show that $f^{[\kappa]}(x,a)$ is strictly decreasing in $x\in(1,\infty)$, it suffices to prove that $R^{[\kappa]}(x,a)-L^{[\kappa]}(x,a)>0$ for all $x>1$, namely,
\begin{align}
&A(x,a)q^{[\kappa]}(x,a+1)q^{[\kappa]}(x,a)+B(x,a)q^{[\kappa]}(x,a+1)q^{[\kappa]}(x,a-1)+\nonumber\\
&~~~C(x,a)q^{[\kappa]}(x,a)q^{[\kappa]}(x,a-1)>(x^a-1)^2q^{[\kappa]}(x,a)^2\label{Equ:ABC}
\end{align}
for all $x>1$, where
\begin{align*}
	A(x,a)&\coloneqq ax^{a+1}-(a+1)x^a+1,\\
	B(x,a)&\coloneqq  x^{2a}-ax^{a+1}+ax^{a-1}-1,\\
	C(x,a)&\coloneqq  (a-1)x^a-ax^{a-1}+1.
\end{align*}
Now we show that $A(x,a), B(x,a),C(x,a)>0$ for all $x, a >1$. For $A(x,a)$ we have $A_x(x,a)=a(a+1)x^{a-1}(x-1)>0$ for all $x>1$, which implies that $A(x,a)$ is strictly increasing in $x\in(1,\infty)$. Thus $A(x,a)>A(1,a)=0$ for all $x>1$. Similarly, we have $C(x,a)>0$ for all $x>1$ whenever $a>1$. For $B(x,a)$ we find  $$B_x(x,a)=2ax^{2a-1}-a(a+1)x^a+a(a-1)x^{a-2}=ax^{a-2}[2x^{a+1}-(a+1)x^2+(a-1)].$$
Since $a>1$ implies that
$$\frac{\partial}{\partial x}[2x^{a+1}-(a+1)x^2+(a-1)]=2(a+1)x(x^{a-1}-1)>0$$
for all $x>1$, so the function $2x^{a+1}-(a+1)x^2+(a-1)$ is strictly increasing in $x\in(1,\infty)$. It follows that $B_x(x,a)>0$ for all $x>1$. Hence $B(x,a)$ is strictly increasing in $x\in(1,\infty)$ and $B(x,a)>B(1,a)=0$ for all $x>1$. 
\par Now we prove \eqref{Equ:ABC}, which can be rewritten as
\begin{multline*}
A(x,a)[q(x,a+1)q(x,a)]^{\kappa-1}+B(x,a)[q(x,a+1)q(x,a-1)]^{\kappa-1}+\\C(x,a)[q(x,a)q(x,a-1)]^{\kappa-1}> q(x,a)^{2\kappa},
\end{multline*}
where $q(x,a) \coloneqq q^{[1]}(x, a) = x^a-1$. When $a=1$, the inequality above becomes 
\[A(x,1)[q(x,2)q(x,1)]^{\kappa-1}>q(x,1)^{2\kappa},\] 
which is true for $x>1$ due to $A(x,1)=q(x,1)^2$ and $q(x,2)=x^2-1>q(x,1)$. So it remains to consider the case $a>1$. A straightforward computation shows that
$$\frac{ A(x,a)}{q(x,a+1)q(x,a)}+\frac{B(x,a)}{q(x,a+1)q(x,a-1)}+\frac{C(x,a)}{q(x,a)q(x,a-1)}=1$$
for all $x>1$. Since $k>1$ and 
$$A(x,a)+B(x,a)+C(x,a)=(x^a-1)^2=q(x,a)^2,$$
it follows by the power mean inequality that
\begin{multline*}
A(x,a)[q(x,a+1)q(x,a)]^{\kappa-1}+B(x,a)[q(x,a+1)q(x,a-1)]^{\kappa-1}+\\C(x,a)[q(x,a)q(x,a-1)]^{\kappa-1}\geq\left(A(x,a)+B(x,a)+C(x,a)\right)^{\kappa}=q(x,a)^{2\kappa}
\end{multline*}
for all $x>1$ with equality if and only if 
$$q(x,a+1)q(x,a)=q(x,a+1)q(x,a-1)=q(x,a)q(x,a-1),$$
or equivalently if and only if $q(x,a+1)=q(x,a)=q(x,a-1)$, which is clearly impossible for $x>1$. This completes the proof of \eqref{Equ:ABC} and hence that of the proposition.
\end{proof}

Next, we prove that $\kF \kappa (x, a)$ is strictly decreasing in $a\in[1,\infty)$ for $x > 1$. Again, it suffices to consider $f^{[\kappa]}(x,a)$.
\begin{prop}\label{prop:f^k(x a) in a}
Let $x>1$ and $\kappa>1$. Then $f^{[\kappa]}(x,a)$ is a strictly decreasing function of $a\in[1,\infty)$.
\end{prop}
\begin{proof}
As in the proof of \Cref{prop:f^k(x a) in x}, we put $q(x,a)=x^a-1>0$. We compute the partial derivative $\partial f^{[\kappa]}/\partial a$ to obtain
	\[
	x^{\kappa}\cdot\frac{\partial f^{[\kappa]}}{\partial a}(x,a)=\frac{\kappa x^{a-1}\log x}{[q(x,a)^{\kappa}-q(x,a-1)^{\kappa}]^2}v^{[\kappa]}(x,a),\]
	where
	\begin{align*}
		v^{[\kappa]}(x,a)&=\left[q(x,a+1)^{\kappa-1}x^2-q(x,a)^{\kappa-1}x\right]\left[q(x,a)^{\kappa}-q(x,a-1)^{\kappa}\right]\\
		&~~~-\left[q(x,a)^{\kappa-1}x-q(x,a-1)^{\kappa-1}\right]\left[q(x,a+1)^{\kappa}-q(x,a)^{\kappa}\right]\\
		&=\left[q(x,a+1)^{\kappa-1}x^2-q(x,a)^{\kappa-1}x\right]\left[q(x,a)^{\kappa-1}(x^a-1)-q(x,a-1)^{\kappa-1}(x^{a-1}-1)\right]\\
		&~~~-\left[q(x,a)^{\kappa-1}x-q(x,a-1)^{\kappa-1}\right]\left[q(x,a+1)^{\kappa-1}(x^{a+1}-1)-q(x,a)^{\kappa-1}(x^a-1)\right]\\
		&=x(1-x)[q(x,a+1)q(x,a)]^{\kappa-1}+(x^2-1)[q(x,a+1)q(x,a-1)]^{\kappa-1}\\
		&~~~+(1-x)[q(x,a)q(x,a-1)]^{\kappa-1}.
	\end{align*}
	Thus it suffices to show $v^{[\kappa]}(x,a)<0$ for all $a>1$. Equivalently, we must show that
	\[x[q(x,a+1)q(x,a)]^{\kappa-1}+[q(x,a)q(x,a-1)]^{\kappa-1}>(1+x)[q(x,a+1)q(x,a-1)]^{\kappa-1},\]
	which can be rewritten as
	\begin{equation}\label{Equprop:f^k(x,a) in a1}
		\frac{x}{1+x}\left(\frac{q(x,a)}{q(x,a-1)}\right)^{\kappa-1}+\frac{1}{1+x}\left(\frac{q(x,a)}{q(x,a+1)}\right)^{\kappa-1}>1.
	\end{equation}
	Simple computation shows that
	\[\frac{x}{1+x}\left(\frac{q(x,a)}{q(x,a-1)}\right)^{-1}+\frac{1}{1+x}\left(\frac{q(x,a)}{q(x,a+1)}\right)^{-1}=1.\]
	Since $\kappa-1>0$, it follows from the power mean inequality that the left side of \eqref{Equprop:f^k(x,a) in a1} is greater than or equal to
	\[\left(\frac{x}{1+x}\left(\frac{q(x,a)}{q(x,a-1)}\right)^{-1}+\frac{1}{1+x}\left(\frac{q(x,a)}{q(x,a+1)}\right)^{-1}\right)^{1-\kappa}=1\]
	with equality if and only if 
	\[\frac{q(x,a)}{q(x,a-1)}=\frac{q(x,a)}{q(x,a+1)}.\]
	Since the equality above does not hold for $x>1$ and $a>1$, we finish the proof of \eqref{Equprop:f^k(x,a) in a1}.
\end{proof}
\vspace*{2mm}
\par For $x,\kappa>1$ and $a\geq1$, we have shown that
\[F^{[\kappa]}(x,a)=\frac{\log f^{[\kappa]}(x,a)}{\kappa\log x}=\frac{1}{\kappa\log x}\log\frac{(x^{a+1}-1)^{\kappa}-(x^{a}-1)^{\kappa}}{(x^{a+1}-x)^{\kappa}-(x^{a}-x)^{\kappa}}\]
is strictly decreasing as $x$ or $a$ increases.

The arithmetic mean--geometric mean inequality implies that $x^{a+1}+x^{a-1}>2x^a$, which is equivalent to 
\begin{equation}\label{eq:AC<B^2}
(x^{a+1}-1)^{\kappa}(x^{a-1}-1)^{\kappa}<(x^{a}-1)^{2\kappa}.
\end{equation}  
Thus we have
\[\frac{(x^{a+1}-1)^{\kappa}-(x^{a}-1)^{\kappa}}{(x^{a+1}-x)^{\kappa}-(x^{a}-x)^{\kappa}}<\frac{(x^{a+1}-1)^{\kappa}}{(x^{a+1}-x)^{\kappa}},\]
which implies that $F^{[\kappa]}(x,a)<F(x,a)$, with $F$ as in \eqref{Equation: Definition of F(x, a)}. In fact, we shall show that as a function of $\kappa$, $F^{[\kappa]}(x,a)$ is strictly increasing on $(1,\infty)$. To this end, let us rewrite
\[\kF \kappa (x,a)=\frac{1}{\kappa\log x}\log\frac{q(x,a+1)^{\kappa}-q(x,a)^{\kappa}}{q(x,a)^{\kappa}-q(x,a-1)^{\kappa}}-1.\]
We need the following elementary lemma.
\begin{lemma}\label{lem:H(t)}
	The function
	\[H(t)=\frac{t\parent{\log t}^2}{(t-1)^2}\]
	is strictly decreasing on $(1,\infty)$.
\end{lemma}
\begin{proof}
	We compute
	\[H'(t)=\frac{(\parent{\log t}^2+2\log t)(t-1)^2-2(t-1)t\parent{\log t}^2}{(t-1)^4}=-\frac{[(t+1)\log t-2(t-1)]\log t}{(t-1)^3}\]
	for all $t>1$. To show $H(t)$ is strictly decreasing on $(1,\infty)$, it is thus sufficient to prove that
	\[\log t>\frac{2(t-1)}{t+1}\]
	for all $t>1$. This follows directly from the fact that
	\[\frac{d}{dt}\left(\log t-\frac{2(t-1)}{t+1}\right)=\frac{(t-1)^2}{t(t+1)^2}>0\]
	for all $t>1$. This proves our lemma.
\end{proof}
Before proving the monotonicity of $F^{[\kappa]}(x,a)$ in $\kappa$, we need to study 
\begin{equation}\label{Equ:h^k}
h^{[\kappa]}(x,a) \coloneqq x^\kappa f^{[\kappa]} (x, a) =\frac{q(x,a+1)^{\kappa}-q(x,a)^{\kappa}}{q(x,a)^{\kappa}-q(x,a-1)^{\kappa}}
\end{equation}
as a function of $\kappa$.
\begin{prop}\label{prop:h^k in k}
	Given any positive real numbers $x>1$ and $a\geq1$, $h^{[\kappa]}(x,a)$ is strictly increasing and strictly log-concave as a function of $\kappa\in(1,\infty)$.
\end{prop}
\begin{proof}
	For $a=1$ we have $h^{[\kappa]}(x,a)=(x+1)^{\kappa}-1$. It is clear that $h^{[\kappa]}(x,1)$ is strictly increasing in $\kappa\in(1,\infty)$. Since
	\[\frac{\partial}{\partial \kappa}(\log h^{[\kappa]}(x,1))=\frac{(x+1)^{\kappa}\log(x+1)}{(x+1)^{\kappa}-1}=\left(1+\frac{1}{(x+1)^{\kappa}-1}\right)\log(x+1)\]
	is strictly decreasing in $\kappa\in(1,\infty)$, we see that $h^{[\kappa]}(x,1)$ is strictly log-concave in $\kappa\in(1,\infty)$.
	\par In what follows, we shall suppose  that $a>1$. For simplicity, let us write $A \coloneqq q(x,a+1)$, $B \coloneqq q(x,a)$ and $C \coloneqq q(x,a-1)$. Then $A>B>C>0$ and $AC<B^2$, the latter of which is equivalent to \eqref{eq:AC<B^2}. We calculate the partial derivative $\partial h^{[\kappa]}/\partial \kappa$ to obtain
	\begin{align*}
		\frac{\partial h^{[\kappa]}}{\partial \kappa} (x,a)&=\frac{(A^{\kappa}\log A-B^{\kappa}\log B)(B^{\kappa}-C^{\kappa})-(A^{\kappa}-B^{\kappa})(B^{\kappa}\log B-C^{\kappa}\log C)}{(B^{\kappa}-C^{\kappa})^2}\\
		&=\frac{1}{(B^{\kappa}-C^{\kappa})^2}\left[(AB)^{\kappa}\log\frac{A}{B} +(BC)^{\kappa}\log\frac{B}{C}-(AC)^{\kappa}\log\frac{A}{C}\right]
	\end{align*}
	for all $\kappa>1$. By the weighted arithmetic mean--geometric mean inequality we have
	\[\frac{\log(A/B)}{\log(A/C)}(B/C)^{\kappa} +\frac{\log(B/C)}{\log(A/C)}(B/A)^{\kappa}>\left[(B/C)^{\log(A/B)}(B/A)^{\log(B/C)}\right]^{\kappa/\log(A/C)}=1.\]
	This implies that $\partial h^{[\kappa]}/\partial \kappa>0$ for all $\kappa>1$. Hence $h^{[\kappa]}(x,a)$ is is strictly increasing in $\kappa\in(1,\infty)$.
	\par Now we show that $h^{[\kappa]}(x,a)$ is strictly log-concave in $\kappa\in(1,\infty)$. Note that
	\begin{align}
		\frac{\partial}{\partial \kappa}(\log h^{[\kappa]}(x,a))&=\frac{1}{(A^{\kappa}-B^{\kappa})(B^{\kappa}-C^{\kappa})}\left[(AB)^{\kappa}\log\frac{A}{B} +(BC)^{\kappa}\log\frac{B}{C}-(AC)^{\kappa}\log\frac{A}{C}\right]\nonumber\\
		&=\frac{A^{\kappa}\log A-B^{\kappa}\log B}{A^{\kappa}-B^{\kappa}}-\frac{B^{\kappa}\log B-C^{v}\log C}{B^{\kappa}-C^{\kappa}}.\label{Equprop:h^k in k}
	\end{align}
	Since
	\begin{equation*}
		\frac{\partial}{\partial \kappa}\left(\frac{A^{\kappa}\log A-B^{\kappa}\log B}{A^{\kappa}-B^{\kappa}}\right) %&=\frac{(A^{\kappa}(\log A)^2-B^{\kappa}(\log B)^2)(A^{\kappa}-B^{\kappa})-(A^{\kappa}\log A-B^{\kappa}\log B)^2 }{(A^{\kappa}-B^{\kappa})^2}\\
		=-\frac{A^{\kappa}B^{\kappa}(\log(A/B))^2}{(A^{\kappa}-B^{\kappa})^2},
	\end{equation*}
	and since 
	\[\frac{\partial}{\partial \kappa}\left(\frac{B^{\kappa}\log B-C^{\kappa}\log C}{B^{\kappa}-C^{\kappa}}\right)=-\frac{B^{\kappa}C^{\kappa}(\log(B/C))^2}{(B^{\kappa}-C^{\kappa})^2}\]
	by symmetry, we have
	\begin{align*}
	\frac{\partial^2}{\partial \kappa^2}(\log h^{[\kappa]}(x,a))&=-\frac{A^{\kappa}B^{\kappa}(\log(A/B))^2}{(A^{\kappa}-B^{\kappa})^2}+\frac{B^{\kappa}C^{\kappa}(\log(B/C)^2)}{(B^{\kappa}-C^{\kappa})^2}
	\end{align*}
	for all $\kappa>1$. Set $r_{\kappa}= (A/B)^{\kappa}$ and $s_{\kappa}= (B/C)^{\kappa}$. Then $s_{\kappa}>r_{\kappa}>1$, since $AC<B^2$. By \Cref{lem:H(t)} we have
	\[\frac{\partial^2}{\partial \kappa^2}(\log h^{[\kappa]}(x,a))=\frac{-H(r_{\kappa})+H(s_{\kappa})}{\kappa^2}<0\]
	for all $\kappa>1$. This proves that $h^{[\kappa]}(x,a)$ is strictly log-concave in $\kappa\in(1,\infty)$ as required.
\end{proof}
We are now ready to show that $F^{[\kappa]}(x,a)$ is a strictly increasing function of $\kappa\in(1,\infty)$.
\begin{prop}\label{prop:F^k in k}
	Given any positive real numbers $x>1$ and $a\geq1$, $F^{[\kappa]}(x,a)$ is strictly increasing in $\kappa\in(1,\infty)$.
\end{prop}
\begin{proof}
	Fixing $x>1$ and $a\geq1$, we have 
	\[F^{[\kappa]}(x,a)=\frac{\log h^{[\kappa]}(x,a)}{\kappa\log x}-1,\]
	where $h^{[\kappa]}(x,a)$ is defined by \eqref{Equ:h^k}. For $a=1$ we have
	\[F^{[\kappa]}(x,1)=\frac{\log((x+1)^{\kappa}-1)}{\kappa\log x}-1.\]
	Since
	\begin{align*}
		\frac{\partial F^{[\kappa]}}{\partial \kappa}(x,1)&=\frac{1}{\kappa^2\log x}\left(\frac{\kappa(x+1)^{\kappa}\log(x+1)}{(x+1)^{\kappa}-1}-\log((x+1)^{\kappa}-1)\right)\\
		&=\frac{1}{\kappa^2\log x}\left(\frac{\kappa\log(x+1)}{(x+1)^{\kappa}-1}-\log(1-(x+1)^{-\kappa})\right)>0
	\end{align*}
	for all $\kappa>1$, it follows that $F^{[\kappa]}(x,1)$ is a strictly increasing function of $\kappa\in(1,\infty)$.
	\par Suppose now that $a>1$. Note that
	\[\frac{\partial F^{[\kappa]}}{\partial \kappa}(x,a)=\frac{1}{\kappa^2\log x}\left(\kappa\cdot\frac{\partial}{\partial\kappa}(\log h^{[\kappa]}(x,a))-\log h^{[\kappa]}(x,a)\right)\]
	and
	\[\frac{\partial}{\partial \kappa}\left(\kappa\cdot\frac{\partial}{\partial\kappa}(\log h^{[\kappa]}(x,a))-\log h^{[\kappa]}(x,a)\right)=\kappa\cdot\frac{\partial^2}{\partial\kappa^2}(\log h^{[\kappa]}(x,a))<0\]
	for all $\kappa>1$, since $h^{[\kappa]}(x,a)$ is strictly log-concave in $\kappa\in(1,\infty)$ by \Cref{prop:h^k in k}. Thus 
	\[\kappa\cdot\frac{\partial}{\partial\kappa}(\log h^{[\kappa]}(x,a))-\log h^{[\kappa]}(x,a)\] 
	is strictly decreasing in $\kappa\in(1,\infty)$. It is clear that
	\[\log h^{[\kappa]}(x,a)=\log\frac{q(x,a+1)^{\kappa}-q(x,a)^{\kappa}}{q(x,a)^{\kappa}-q(x,a-1)^{\kappa}}<{\kappa}\log\frac{q(x,a+1)}{q(x,a)},\]
	since $q(x,a+1)q(x,a-1)<q(x,a)^2$. By \eqref{Equprop:h^k in k} we have
	\[\lim\limits_{k\to\infty}\frac{\partial}{\partial\kappa}(\log h^{[\kappa]}(x,a))=\log\frac{q(x,a+1)}{q(x,a)}.\]
	Since $h^{[\kappa]}(x,a)$ is strictly log-concave in $\kappa\in(1,\infty)$ by \Cref{prop:h^k in k}, we know that $\partial h^{[\kappa]}/\partial\kappa$ is strictly decreasing in $\kappa\in(1,\infty)$. Hence 
	\[\frac{\partial}{\partial\kappa}(\log h^{[\kappa]}(x,a))>\log\frac{q(x,a+1)}{q(x,a)}\]
	for all $\kappa\in(1,\infty)$. It follows that
	\[\frac{\partial h^{[\kappa]}}{\partial \kappa}(x,a)=\frac{1}{\kappa^2\log x}\left(\kappa\cdot\frac{\partial}{\partial\kappa}(\log h^{[\kappa]}(x,a))-\log h^{[\kappa]}(x,a)\right)>0\]
	for all $\kappa\in(1,\infty)$. Hence $F^{[\kappa]}(x,a)$ is a strictly increasing function of $\kappa\in(1,\infty)$.
\end{proof}

\begin{proof}[Proof of \Cref{Theorem: behavior of kF}]
The monotonicity part follows from \Cref{prop:f^k(x a) in x}, \Cref{prop:f^k(x a) in a} and \Cref{prop:F^k in k}. To proceed, we observe that $F^{[\kappa]}(x,a)\to0$ as $x\to\infty$ as a consequence of \Cref{prop:f^k(x a) in x}. Since $F^{[\kappa]}(x,a)$ is a strictly decreasing function of $x\in(1,\infty)$, we have $F^{[\kappa]}(x,a)>0$ for all $x>1$. Note that 
\begin{align*}
	f^{[\kappa]}(x,a)=\frac{(x^{a+1}-1)^{\kappa}-(x^{a}-1)^{\kappa}}{(x^{a+1}-x)^{\kappa}-(x^{a}-x)^{\kappa}}&>\frac{(x^{a+1}-1)^{\kappa}-(x^{a}-1)^{\kappa}+(x^{a}-x)^{\kappa}}{(x^{a+1}-x)^{\kappa}}\\
	&=\left(\frac{x^{a+1}-1}{x(x^{a}-1)}\right)^{\kappa}+\left(\frac{x^{a-1}-1}{x^{a}-1}\right)^{\kappa}-x^{-{\kappa}}\\
	&\to\left(1+\frac{1}{a}\right)^{\kappa}+\left(1-\frac{1}{a}\right)^{\kappa}-1
\end{align*}
as $x\to1^+$. By Bernoulli's inequality we have
\[\left(1+\frac{1}{a}\right)^{\kappa}+\left(1-\frac{1}{a}\right)^{\kappa}-1>\left(1+\frac{\kappa}{a}\right)+\left(1-\frac{\kappa}{a}\right)-1=1.\]
Hence $F^{[\kappa]}(x,a)\to\infty$ as $x\to1^+$, as desired.
\end{proof}

We now define two partial inverses to $\kF \kappa (x, a)$.
\begin{definition}\label{Definition: kx and ka}
	For any $x > 1$, $\kF \kappa(x, \cdot) : [1, \infty) \to \left(0, \frac{\log\parent{(x+1)^\kappa - 1)}}{\kappa \log x} - 1\right]$ is strictly decreasing, and so it has a (strictly decreasing) inverse
	\[
	\ka \kappa x : \left(0, \frac{\log\parent{(x+1)^\kappa - 1)}}{\kappa \log x} - 1\right] \to [1, \infty).
	\]
	Likewise, for any $a \geq 1$, $\kF \kappa (\cdot, a) : (1, \infty) \to (0, \infty)$ is strictly decreasing, and so it has a (strictly decreasing) inverse
	\[
	\kx \kappa a : (0, \infty) \to (1, \infty).
	\]
\end{definition}

By \Cref{Theorem: behavior of kF}, $\ka \kappa x$ is decreasing as a function of $x$, $\kx \kappa a$ is decreasing as a function of $a$, and both $\ka \kappa x$ and $\kx \kappa a$ are increasing as functions of $\kappa$. We are interested in the equality
\[
\kF \kappa (x_a, a) = \epsilon,
\]
so we write $\epsilon$ for the argument of $\ka \kappa x$ and of $\kx \kappa a$.

We now study $\kx \kappa a$ more carefully. To ease notation, we write $x_a$ for $\kx \kappa a$ where no confusion arises. We typically view $x_a$ as a function of $\epsilon$, as noted above, but we also may consider it as a function of $x$ via the implict function theorem and the following equality:
\[
\kF \kappa (x, 1) = \kF \kappa ({x_a}, a).
\]
Note that $x_a$ is \emph{increasing} as a function of $x$, where it is \emph{decreasing} as a function of $\epsilon$.

\begin{lemma}\label{Proposition: asymptotic for xell}%See Lemma 6.18
For any fixed $\kappa>1$ and $a\geq1$, we have $x_a\sim(ax)^{1/a}$ as $\epsilon\to 0^+$ (or equivalently, as $x \to \infty$).
\end{lemma}

\begin{proof}
Note that $x_a\to\infty$ as $\epsilon\to0^+$. If $a\geq1$, then
\[F^{[\kappa]}(x,a)=\frac{\log f^{[\kappa]}(x,a)}{\kappa\log x}=\frac{1}{\kappa\log x}\left(f^{[\kappa]}(x,a)-1+O\left((f^{[\kappa]}(x,a)-1)^2\right)\right)\]
for sufficiently large $x$, where $f^{[\kappa]}(x,a)$ is defined by \eqref{Equ:f^k(x,a)}. Since $(1-t)^{\kappa}=1-\kappa t+O(t^2)$ for $t\in[0,1]$, it follows that
\begin{align*}
	f^{[\kappa]}(x,a)-1&=\frac{(1-x^{-a-1})^\kappa-(1-x^{-a})^{\kappa}x^{-\kappa}}{(1-x^{-a})^{\kappa}-(1-x^{-a+1})^{\kappa}x^{-\kappa}}-1\\
	&=\frac{(1-\kappa x^{-a-1}+O(x^{-2a-2}))-(1-\kappa x^{-a}+O(x^{-2a}))x^{-\kappa}}{(1-
		\kappa x^{-a}+O(x^{-2a}))-(1-\kappa x^{-a+1}+O(x^{-2a+2}))x^{-\kappa}}-1\\
	&=\frac{\kappa x^{-a}(1-x^{-1})-\kappa x^{-a-\kappa+1}(1-x^{-1})+O(x^{-2a}+x^{-2a+2-\kappa})}{(1-\kappa x^{-a}+O(x^{-2a}))-(1-\kappa x^{-a+1}+O(x^{-2a+2}))x^{-\kappa}}\\
	&=\frac{\kappa x^{-a}(1+O(x^{1-\min(\kappa,2)}))}{1+O(x^{-\min(\kappa,a)})}\\
	&=\kappa x^{-a}\left(1+O\left(x^{1-\min(\kappa,2)}+x^{-\min(\kappa,a)}\right)\right)\\
	&=\kappa x^{-a}\left(1+O\left(x^{1-\min(\kappa,2)}\right)\right),
\end{align*}
where we have used the assumption that $\kappa>1$ and $a\geq1$. Thus we have
\begin{equation}\label{L3Equ1}
	F^{[\kappa]}(x,a)=\frac{1+O\left(x^{1-\min(\kappa,2)}+x^{-a}\right)}{x^a\log x}=\frac{1+O\left(x^{1-\min(\kappa,2)}\right)}{x^a\log x}\sim\frac{1}{x^a\log x}
\end{equation}
as $x\to\infty$, where all the implicit constants depend only on $\kappa$. In particular, we obtain
\[F^{[\kappa]}(x_a,a)\sim\frac{1}{x_a^a\log x_a}.\]
as $\epsilon\to0^+$. Since $F^{[\kappa]}(x_a,a)=F^{[\kappa]}(x,1)=\epsilon$, we see that $x_a^a\log x_a\sim x\log x$ as $\epsilon\to0^+$. The lemma follows.
\end{proof}

Now we specialize in the case $a=2$. The following lemma provides a more precise asymptotic for $x_2$.

\begin{lemma}\label{Lemma: asymptotic for x2 with error}
For any $\kappa>1$ and $\epsilon>0$, we have
\begin{equation}
x_2=\sqrt{2x}\left(1-\frac{\log 2}{2\log x}+O\left(\frac{1}{\parent{\log x}^2}\right)\right)\label{Equation: asymptotic for x2 with error}
\end{equation}
for sufficiently large $x$, where the implicit constant in the error term depends only on $\kappa$.
\end{lemma}

\begin{proof}
Let $\kappa>1$ and $\epsilon>0$. Denote by $\xi=\xi(\kappa,\epsilon)>1$ the unique solution to the equation $\xi^2\log\xi=x\log x$. Since $x_2^2\log x_2\sim x\log x$ as $\epsilon\to0^+$, as shown in the proof of \Cref{Proposition: asymptotic for xell}, one may think of $\xi$ as a proxy to $x_2$. Then $\xi=\sqrt{2x}(1-\eta)$, where $\eta\to0$ as $\epsilon\to0^+$ by \Cref{Proposition: asymptotic for xell}. Carrying this back into the equation $\xi^2\log\xi=x\log x$, we obtain 
\[(1-2\eta+O(\eta^2))\left(1+\frac{\log 2}{\log x}-\frac{2\eta+O(\eta^2)}{\log x}\right)=1,\]
which implies that $\eta=(\log 2+o(1))/(2\log x)$. From this it follows that
\begin{equation}\label{L5Equ1}
	\xi=\sqrt{\frac{x\log x}{\log \xi}}=\sqrt{\frac{2x\log x}{\log x+\log2+O(1/\log x)}}=\sqrt{2x}\left(1-\frac{\log2}{2\log x}+O\left(\frac{1}{\parent{\log x}^2}\right)\right).
\end{equation}
By \eqref{L3Equ1} we have 
\[F^{[\kappa]}(x_2,2)=F^{[\kappa]}(x,1)=\frac{1}{x\log x}+O\left(\frac{1}{x^{\delta_{\kappa}}\log x}\right)\]
and
\[F^{[\kappa]}(\xi,2)=\frac{1}{\xi^2\log\xi}+O\left(\frac{1}{\xi^{1+\delta_{\kappa}}\log\xi}\right)=\frac{1}{x\log x}+O\left(\frac{1}{x^{(1+\delta_{\kappa})/2}\log x}\right),\]
where $\delta_{\kappa}=\min(\kappa,2)\in(1,2]$. It follows that 
\begin{equation}\label{L5Equ2}
	F^{[\kappa]}(\xi,2)-F^{[\kappa]}(x_2,2)=O\left(\frac{1}{x^{(1+\delta_{\kappa})/2}\log x}\right).
\end{equation}
On the other hand, we have
\[\frac{\partial F^{[\kappa]}}{\partial x}(x,2)=\left(-\frac{F^{[\kappa]}(x,2)}{x}+\frac{1}{\kappa f^{[\kappa]}(x,2)}\cdot\frac{\partial f^{[\kappa]}}{\partial x}(x,2)\right)\frac{1}{\log x}.\]
Note that
\[\frac{1}{\kappa f^{[\kappa]}(x,2)}\cdot\frac{\partial f^{[\kappa]}}{\partial x}(x,2)=\frac{g^{[\kappa]}(x)}{[(x^3-1)^{\kappa}-(x^2-1)^{\kappa}][(x^3-x)^{\kappa}-(x^2-x)^{\kappa}]}\sim\frac{g^{[\kappa]}(x)}{x^{6\kappa}}\]
as $x\to\infty$, where 
\begin{align*}
	g^{[\kappa]}(x)&=\left[3(x^3-1)^{\kappa-1}x^2-2(x^2-1)^{\kappa-1}x\right]\left[(x^3-x)^{\kappa}-(x^2-x)^{\kappa}\right]\\
	&~~~-\left[(x^3-1)^{\kappa}-(x^2-1)^{\kappa}\right]\left[(x^3-x)^{\kappa-1}(3x^2-1)-(x^2-x)^{\kappa-1}(2x-1)\right].
\end{align*}
But for sufficiently large $x$, we find that
\[3(x^3-1)^{\kappa-1}x^2-2(x^2-1)^{\kappa-1}x=3x^{3\kappa-1}-(3\kappa-3)x^{3\kappa-4}-2x^{2\kappa-1}+O\left(x^{3\kappa-7}+x^{2\kappa-3}\right),\]
\[	(x^3-x)^{\kappa}-(x^2-x)^{\kappa}=x^{3\kappa}-\kappa x^{3\kappa-2}-x^{2\kappa}+O\left(x^{3\kappa-4}+x^{2\kappa-1}\right),\]
\[(x^3-1)^{\kappa}-(x^2-1)^{\kappa}=x^{3\kappa}-\kappa x^{3\kappa-3}-x^{2\kappa}+O\left(x^{3\kappa-6}+x^{2\kappa-2}\right),\]
\[(x^3-x)^{\kappa-1}(3x^2-1)-(x^2-x)^{\kappa-1}(2x-1)=3x^{3\kappa-1}-(3\kappa-2)x^{3\kappa-3}-2x^{2\kappa-1}+O\left(x^{3\kappa-5}+x^{2\kappa-2}\right).\]
It follows that
\begin{align*}
	g^{[\kappa]}(x)&=\left(3x^{6\kappa-1}-3\kappa x^{6\kappa-3}-5x^{5\kappa-1}+O\left(x^{5\kappa-2}+x^{6\kappa-4}\right)\right)\\
	&~~~-\left(3x^{6\kappa-1}-(3\kappa-2)x^{6\kappa-3}-5x^{5\kappa-1}+O\left(x^{5\kappa-2}+x^{6\kappa-4}\right)\right)\\
	&=-2x^{6\kappa-3}+O\left(x^{5\kappa-2}+x^{6\kappa-4}\right)
\end{align*}
for sufficiently large $x$. Hence
\[\frac{1}{\kappa f^{[\kappa]}(x,2)}\cdot\frac{\partial f^{[\kappa]}}{\partial x}(x,2)\sim-2x^{-3}\]
as $x\to\infty$. Combining this with \eqref{L3Equ1} we obtain
\[\frac{\partial F^{[\kappa]}}{\partial x}(x,2)=-\frac{2+o(1)}{x^3\log x}\]
as $x\to\infty$. Now the mean value theorem implies that there exists $\eta=\eta(\kappa,\epsilon)>1$ which lies between $\xi$ and $x_2$, such that
\[|F^{[\kappa]}(\xi,2)-F^{[\kappa]}(x_2,2)|=\left|\frac{\partial F^{[\kappa]}}{\partial x}(\eta,2)(\xi-x_2)\right|.\]
Since $\xi\sim x_2\sim\sqrt{2x}$ as $\epsilon\to0^+$, we must have
\[\left|\frac{\partial F^{[\kappa]}}{\partial x}(\eta,2)\right|\gg\frac{1}{x^{3/2}\log x}.\]
when $\epsilon$ is sufficiently small. Together with \eqref{L5Equ1} and \eqref{L5Equ2} this implies that
\[x_2=\xi+O\left(x^{1-\delta_{\kappa}/2}\right)=\sqrt{2x}\left(1-\frac{\log2}{2\log x}+O\left(\frac{1}{\parent{\log x}^2}\right)\right)\]
for sufficiently small $\epsilon$, since $1-\delta_{\kappa}/2<1/2$. This completes the proof of the lemma.
\end{proof}

\begin{remark}
	The quantity $\xi$ defined by $\xi^2 \log \xi = x \log x$ may be written explicitly in terms of $x$ by recourse to (the principal branch of) the Lambert $W$-function, which is defined to be the inverse of $x \mapsto x e^x$. Indeed, we have 
	\[
	\xi = e^{W(2 x \log x)/2} = \sqrt{\frac{2 x \log x}{W(2 x \log x)}}.
	\]
	The asymptotic \eqref{L5Equ1} then follows from the asymptotic series representation for $W$ around $\infty$, given by
	\begin{equation}\label{Equation: Asymptotic for W}
	W(x) = \log x - \log \log x + \sum_{m \geq 1} \sum_{n \geq 0} \frac{(-1)^n}{m!} \begin{bmatrix} n + m \\ n + 1 \end{bmatrix} \frac{\parent{\log \log x}^m}{\parent{\log x}^{m+n}},
	\end{equation}
	where $\begin{bmatrix} n \\ m \end{bmatrix}$ is an unsigned Stirling number of the first kind. The expression given in \eqref{Equation: Asymptotic for W} converges for real $x > e$. See \cite[Section 4.1.4]{Mezo} for more details on the asymptotics of the Lambert $W$-function. 
\end{remark}

\Cref{Lemma: asymptotic for x2 with error} implies that $x_2<\sqrt{2x}$ for sufficiently large $x$. Now we prove the following stronger result in the case $\kappa\ge3/2$ which holds for all $x>1$.

\begin{lemma}\label{Lemma:upper bounds for x_2}
For any $\kappa\geq3/2$ and any $\epsilon>0$, we have $x_2<\sqrt{2x}$.
\end{lemma}
\begin{proof}
	If $x\leq 2$, then $x_2<x\leq\sqrt{2x}$. Suppose now that $x>2$. Put $y=\sqrt{2x}>2$. The inequality $x_2<\sqrt{2x}$ is equivalent to $F^{[\kappa]}(y^2/2,1)=F^{[\kappa]}(x_2,2)>F^{[\kappa]}(y,2)$. That is,
	\[\frac{1}{\kappa\log(y^2/2)}\log\frac{(y^2+2)^{\kappa}-2^{\kappa}}{y^{2\kappa}}>\frac{1}{\kappa\log y}\log\frac{(y^2+y+1)^{\kappa}-(y+1)^{\kappa}}{((y+1)^{\kappa}-1)y^{\kappa}}.\]
	Since $\log(y^2/2)<2\log y$, it suffices to show
	\[(y^2+2)^{\kappa}-2^{\kappa}>\left(\frac{(y^2+y+1)^{\kappa}-(y+1)^{\kappa}}{(y+1)^{\kappa}-1}\right)^2,\]
	which can be rewritten as
	\begin{equation}\label{L6Equ1}
		\frac{\left(1-\left(\frac{2}{y^2+2}\right)^{\kappa}\right)\left(1-\left(\frac{1}{y+1}\right)^{\kappa}\right)^2}{\left(1-\left(\frac{y+1}{y^2+y+1}\right)^{\kappa}\right)^{2}\left(1-\frac{2y+1}{(y^2+2)(y+1)^2}\right)^{\kappa}}>1.
	\end{equation}
	By Lagrange's mean value theorem, there exists 
	\[\xi\in\left(\frac{1}{y+1},\frac{y+1}{y^2+y+1}\right)\]
	such  that
	\[2\log\left(1-\left(\frac{1}{y+1}\right)^{\kappa}\right)-2\log\left(1-\left(\frac{y+1}{y^2+y+1}\right)^{\kappa}\right)=\frac{2\kappa\xi^{\kappa-1}y}{(1-\xi^{\kappa})(y+1)(y^2+y+1)}.\]
	Since
	\[\frac{\xi^{\kappa-1}}{1-\xi^{\kappa}}>\frac{1}{(y+1)^{\kappa-1}},\]
	we obtain
	\[2\log\left(1-\left(\frac{1}{y+1}\right)^k\right)-2\log\left(1-\left(\frac{y+1}{y^2+y+1}\right)^{\kappa}\right)>\frac{2\kappa y}{(y+1)^{\kappa}(y^2+y+1)}.\]
	It is clear from the inequality $\log(1-z)<-z$ for $z\in(0,1)$ that
	\[-\kappa\log\left(1-\frac{2y+1}{(y^2+2)(y+1)^2}\right)>\frac{\kappa(2y+1)}{(y^2+2)(y+1)^2}.\]
	From the inequality $\log(1-z)>-z-3z^2/4$ for all $z\in(0,1/3]$ it follows that
	\[\log\left(1-\left(\frac{2}{y^2+2}\right)^{\kappa}\right)>-\left(\frac{2}{y^2+2}\right)^{\kappa}-\frac{3}{4}\left(\frac{2}{y^2+2}\right)^{2\kappa}.\]
	Thus the natural logarithm of the left-hand side of \eqref{L6Equ1} is greater than
	\[\frac{2\kappa y}{(y+1)^{\kappa}(y^2+y+1)}+\frac{\kappa(2y+1)}{(y^2+2)(y+1)^2}-\left(\frac{2}{y^2+2}\right)^{\kappa}-\frac{3}{4}\left(\frac{2}{y^2+2}\right)^{2\kappa}.\]
	We have to show that the above expression is greater than 0, or equivalently,
	\[\frac{2\kappa y}{y^2+y+1}+\frac{\kappa(2y+1)(y+1)^{\kappa-2}}{y^2+2}-\left(\frac{2(y+1)}{y^2+2}\right)^{\kappa}-\frac{3}{4}\left(\frac{2\sqrt{y+1}}{y^2+2}\right)^{2\kappa}>0.\]
 Since $y^2+2>2(y+1)$ and $y^2+2>2\sqrt{y+1}$ for any $y>2$, each term on the left-hand side with the sign attached is an increasing function of $\kappa\in[3/2,\infty)$. Hence it suffices to prove
	\begin{equation}\label{L6Equ2}
		\frac{3y}{y^2+y+1}+\frac{3(2y+1)}{2(y^2+2)\sqrt{y+1}}-\left(\frac{2(y+1)}{y^2+2}\right)^{3/2}-\frac{3}{4}\left(\frac{2\sqrt{y+1}}{y^2+2}\right)^{3}>0
	\end{equation}
	for all $y>2$. But
	\begin{align*}
		\frac{2y+1}{(y^2+2)\sqrt{y+1}}-\left(\frac{2\sqrt{y+1}}{y^2+2}\right)^{3}&=\frac{(2y+1)(y^2+2)^2-8(y+1)^2}{(y^2+2)^3\sqrt{y+1}}\\
		&>4\cdot\frac{(y^2+2)^2-(2y+2)^2}{(y^2+2)^3\sqrt{y+1}}>0,
	\end{align*}
	where the last inequality holds because $2y+1>4$. Using that facts that $3y>2(y+1)$ and that $2y+1>\sqrt{y+1}$, the left-hand side of \eqref{L6Equ2} is 
	\begin{align}
   &>\frac{3y}{y^2+y+1}+\frac{3(2y+1)}{4(y^2+2)\sqrt{y+1}}-\left(\frac{2(y+1)}{y^2+2}\right)^{3/2}\nonumber\\
   &>\frac{2(y+1)}{y^2+y+1}+\frac{3\sqrt{y+1}}{4(y^2+2)}-\left(\frac{2(y+1)}{y^2+2}\right)^{3/2}\nonumber\\
   &=\frac{(y+1)^{3/2}}{4(y^2+2)^{3/2}}\left(\frac{8(y^2+2)^{3/2}}{(y^2+y+1)\sqrt{y+1}}+\frac{3\sqrt{y^2+2}}{y+1}-8\sqrt{2}\right).\label{L6Equ3}
	\end{align}
	To prove \eqref{L6Equ2}, it suffices to show that the factor in the parentheses in \eqref{L6Equ3} is positive. Since
	\[\frac{d}{dy}\left(\frac{y^2+2}{(y+1)^2}\right)=\frac{2(y-2)}{(y+1)^3}>0\]
	for all $y>2$, it follows that
	\[\frac{3\sqrt{y^2+2}}{y+1}=3\sqrt{\frac{y^2+2}{(y+1)^2}}\]
	is strictly increasing on $[2,\infty)$. Similarly, we have
	\[\frac{d}{dy}\left(\frac{2(y^2+2)^{3/2}}{(y^2+y+1)\sqrt{y+1}}\right)=\frac{\sqrt{y^2+2}}{(y^2+y+1)^2(y+1)^{3/2}}\left(y^4+5y^3-y^2-8y-6\right)>0,\]
	since 
	\[y^4+5y^3-y^2-8y-6=(y^2(y^2-1)-6)+y(5y^2-8)>0\]
	for all $y>2$. Thus 
	\[\frac{8(y^2+2)^{3/2}}{(y^2+y+1)\sqrt{y+1}}\]
	is also strictly increasing on $[2,\infty)$. It follows that the factor in the parentheses in \eqref{L6Equ3} is greater than
	\[\frac{48\sqrt{2}}{7}+\sqrt{6}-8\sqrt{2}>\frac{48\sqrt{2}}{7}+\frac{3\sqrt{2}}{2}-8\sqrt{2}=\frac{5\sqrt{2}}{14}>0\]
	for all $y>2$. This completes the proof that \eqref{L6Equ3} is positive and hence that of the lemma.
\end{proof}

\begin{remark}\label{rmk:Lemma:upper bounds for x_2}
	It seems that \Cref{Lemma:upper bounds for x_2} holds for all $\kappa>1$. This would follow if one could show that the function
	\[[(x^2+2)^{\kappa}-2^{\kappa}][(x+1)^{\kappa}-1]^2-[(x^2+x+1)^{\kappa}-(x+1)^{\kappa}]^2\]
	is strictly increasing in $x\in[0,\infty)$ and $\kappa\in[1,\infty)$. In fact, monotonicity in either of the two variable suffices.
\end{remark}

If $\kappa = 2$, we can strengthen \Cref{Lemma: asymptotic for x2 with error} substantially.

\begin{lemma}\label{Lemma: lower bound on x2 when kappa = 2}
	We have
	\begin{equation}
	\sqrt{2x} \parent{1 - \frac{\log 2}{2 \log x}} < \kx 2 2\label{Equation: lower bound on kx 2 2}
	\end{equation}
	whenever $x \geq 2^{15}$.
\end{lemma}

\begin{proof}
	We follow the proof of \cite[Lemma 6.17]{Broughan}. Let
	\begin{equation}
	y \coloneqq \sqrt{2x} \parent{1 - \frac {\log 2}{2 \log x}}, \label{Equation: y in terms of x}
	\end{equation}
	so \eqref{Equation: lower bound on kx 2 2} is equivalent to
	\begin{equation}
	\frac{\log\parent{1 + \frac{2}{x}}}{2 \log x} = \kF 2(x, 1) = \kF 2(x_2, 2) < \kF 2(y, 2) = \frac{\log\parent{1 + \frac{2}{y^2 + 2 y}}}{2 \log y}.\label{Equation: First restatement x2 > sqrt(2x)(1 - u)}
	\end{equation}
	
	Recall that
	\begin{equation}
	\frac 1{t+1} < \log\parent{1 + \frac 1t} < \frac 1 t\label{Equation: elementary inequality for log(1 + t)}
	\end{equation}
	whenever $t > 1$. Using \eqref{Equation: elementary inequality for log(1 + t)}, we find \eqref{Equation: First restatement x2 > sqrt(2x)(1 - u)} follows from
	\begin{equation}
	2\parent{y^2 + 2 y + 2}\frac{\log y}{\log x} < 2x.\label{Equation: Second restatement x2 > sqrt(2x)(1 - u)}
	\end{equation}
	For ease of notation, we write $u \coloneqq \frac{\log 2}{\log x}$, so $y = \sqrt{2 x} \parent{1 - u/2}$. We have
	\begin{equation}
	2\frac{\log y}{\log x} = \frac{\log 2 + \log x + 2\log\parent{1 - \frac {\log 2}{2 \log x}}}{\log x} < %\frac{\log 2 + \log x - \frac {\log 2}{\log x}}{\log x} = 
	1 + \frac{\log 2}{\log x} - \frac{\log 2}{\parent{\log x}^2} = 1 + u - \frac{u^2}{\log 2},\label{Equation: first weakening of x2 > sqrt(2x)(1 - u)}
	\end{equation}
	by \eqref{Equation: y in terms of x} and because $\log\parent{1 - t} < -t$ whenever $0 < t < 1$.
	
	Note that $u \leq 1/15$ by the hypothesis $x \geq 2^{15}$. We now estimate the left-hand side of \eqref{Equation: Second restatement x2 > sqrt(2x)(1 - u)} in two pieces. First, we have
	\[
	2y^2 \frac{\log y}{\log x} < 2 x \parent{1 - \frac {u}{2}}^2 \parent{1 + u - \frac{u^2}{\log 2}} < 2 x \parent{1 - 2.0804 u^2},
	\]
	where the first inequality follows from \eqref{Equation: first weakening of x2 > sqrt(2x)(1 - u)} and the last inequality follows because
	\begin{align*}
		\parent{1 - \frac {u}{2}}^2 \parent{1 + u - \frac{u^2}{\log 2}} &= 1 - \parent{3 + \frac{4}{\log 2} - \parent{1 + \frac{4}{\log 2}} u + \frac{u^2}{\log 2}}\frac{u^2}{4} \\
		&\leq 1 - \parent{3 + \frac{4}{\log 2} - \parent{1 + \frac{4}{\log 2}} u + \frac{u^2}{\log 2}}\bigg|_{u=\frac{1}{15}} \frac{u^2}{4} \\
		&< 1 - 2.0804 u^2.
	\end{align*}
	Second, we have
	\begin{align*}
	4\parent{y + 1}\frac{\log y}{\log x} &< 2 \parent{y + 1}\parent{1 + u - \frac{u^2}{\log 2}} \\
	% &= 2 \parent{\sqrt{2 x} \parent{1 - \frac u 2} + 1}\parent{1 + u - \frac{u^2}{\log 2}} \\
	&\leq \parent{2 \parent{1 - \frac u 2}\parent{1 + u - \frac{u^2}{\log 2}} \sqrt{2 x} + 2 \parent{1 + u - \frac{u^2}{\log 2}}}\bigg|_{u=\frac{1}{15}} \\
	&\leq 2.059 \sqrt{2x},
	\end{align*}
	where the first inequality again follows from \eqref{Equation: first weakening of x2 > sqrt(2x)(1 - u)}
	Therefore,
	\[
	2\parent{y^2 + 2 y + 2}\frac{\log y}{\log x} \leq 2 x \parent{1 - 2.0804 u^2} + 2.059 \sqrt{2x} < 2x,
	\]
	where the final inequality follows because
	\[
	-4.1608 u^2 x + 2.059 \sqrt{2x} < -67 < 0
	\]
	for $x \geq 2^{15}$. This completes the proof.
\end{proof}

\subsection{Understanding the \texorpdfstring{$\kappa$}{k}-colossally abundant numbers}

With the notation and results we have developed, we can now understand $\kappa$-colossally abundant numbers using $\kF \kappa$.

\begin{proposition}\label{Proposition: structure of k-colossally abundant numbers} %See Lemma 6.9
	If $N$ is $\kappa$-colossally abundant for $\epsilon$, then for every prime $p$, we have
	\begin{equation}
	\kF \kappa (p, \vp p (N) + 1) \leq \epsilon \leq \kF \kappa (p, \vp p (N)),\label{Equation: epsilon bounded by kF on both sides}
	\end{equation}
	where $\vp p$ is the $p$-adic valuation on $\bbQ$. Conversely, if $N$ is a positive integer and $\epsilon > 0$ is a real number, and \eqref{Equation: epsilon bounded by kF on both sides} holds for every prime $p$, then $N$ is $\kappa$-colossally abundant for $\epsilon$.
\end{proposition}

\begin{proof}
	Let $\kappa > 1$ and $\epsilon > 0$ be given, and let $p$ be a prime. From \Cref{Definition: kF and kfp} and \eqref{Equation: kfp versus ksigma}, we conclude
	\[
	\kfp \kappa \epsilon (p, 1) \kfp \kappa \epsilon (p, 2) \cdots \kfp \kappa \epsilon (p, a) = \krho{\kappa}{\epsilon}(p^a)=\frac{\ksigma \kappa (p^a)}{p^{a \kappa(1 + \epsilon)}}.
	\]
	Now $\frac{\ksigma \kappa (p^a)}{p^{a \kappa (1 + \epsilon)}} \geq \frac{\ksigma \kappa (p^{a-1})}{p^{(a-1)\kappa(1 + \epsilon)}}$ if and only if $\kfp \kappa \epsilon (p, a) \geq 1$. By \Cref{prop:f^k(x a) in a}, $\kfp \kappa \epsilon(p, a)$ is a strictly decreasing function of $a\ge1$, and $\lim_{a \to \infty} \kfp \kappa \epsilon (p, a) = p^{-\kappa \epsilon} < 1$. Write $X_p \coloneqq \set{p^a \ : \ a \in \bbZ_{\geq 0}}$ for ease of notation. If $\kfp \kappa \epsilon (p, a) < 1$ for all $a \geq 1$, then $\krho \kappa \epsilon |_{X_p}$ is maximized at $p^a = 1$; otherwise, it is maximized when $\kfp \kappa \epsilon (p, a) \geq 1$ and $\kfp \kappa \epsilon (p, a + 1) \leq 1$. In particular, if $\kfp \kappa \epsilon (p, a) > 1$ for some $a$ and $\kfp \kappa \epsilon (p, a) \neq 1$ for all $a$, then $\krho \kappa \epsilon |_{X_p}$ is maximized when $a$ is chosen so $\kfp \kappa \epsilon (p, a) > 1$ and $\kfp \kappa \epsilon (p, a + 1) < 1$; if $\kfp \kappa \epsilon (p, a) = 1$ for some $a \geq 1$, then $p^{a-1}$ and $p^a$ are the unique global maxima of $\krho \kappa \epsilon |_{X_p}$, and of course $\krho \kappa \epsilon |_{X_p} (p^{a-1}) = \krho \kappa \epsilon |_{X_p} (p^a)$. With the convention $\kF \kappa (p, 0) \coloneqq \infty$, we have shown $\krho \kappa \epsilon |_{X_p}$ is maximized at $p^a$ if
	\[
	\kF \kappa (p, a + 1) \leq \epsilon \leq \kF \kappa (p, a).
	\]
	Observe also that $\lim_{p \to \infty} \kfp \kappa \epsilon (p, 1) = 0$, so $\krho \kappa \epsilon |_{X_p}$ is maximized at $p^a = 1$ for all $p$ sufficiently large. 
	
	Now as $\krho \kappa \epsilon$ is multiplicative, every $\kappa$-colossally abundant number $N$ for $\epsilon$ is of the form $N = \prod_{p \ \text{prime}} p^{a_p}$, where $p^{a_p}$ maximizes $\krho \kappa \epsilon |_{X_p}$, and every number obtained in this fashion is $\kappa$-colossally abundant number for $\epsilon$.
\end{proof}

The proof of \Cref{Proposition: structure of k-colossally abundant numbers} above shows that the triples $(p, a, \epsilon)$ for which $\kfp \kappa \epsilon (p, a)= 1$ are especially important to our understanding of $\kappa$-colossally abundant numbers, so we make the following definition.

\begin{definition}\label{Definition: kE}
	Let
	\[
	\kEp \kappa p \coloneqq \set{\kF \kappa (p, j) \ : \ j \in \bbZ_{> 0}}
  \]
  for $p$ prime, and define
  \[
	\kE \kappa \coloneqq \bigcup_{p \ \text{prime}} \kEp \kappa p \cup \{\infty\}.
	\]
	Write $\kE \kappa = \set{\kepsilon \kappa i}_{i = 0}^\infty$, where $\infty = \kepsilon \kappa 0 > \kepsilon \kappa 1 > \kepsilon \kappa 2 > \dots > 0$, and $\lim_{i \to \infty} \kepsilon \kappa i = 0$. For $0 < \epsilon < \infty$, if $\epsilon \in \kE \kappa$, we say $\epsilon$ is a \defi{$\kappa$-critical value}. Otherwise, $\epsilon$ is a \defi{$\kappa$-noncritical value}.
\end{definition}

Thus, a real number $\epsilon$ is a critical value precisely if $\kfp \kappa \epsilon (p, a) = 1$ for some prime $p$ and some integer $a$.

\begin{definition}\label{Definition: N_epsilon}
	For each $\epsilon > 0$, we define $\kN \kappa \epsilon$ to be the largest $\kappa$-colossally abundant number for $\epsilon$. 
\end{definition}

\begin{example}\label{Example: Ep and Eq intersect}
	Let $\kappa = 2.956801214357021\ldots$, and note that 
	\[
	\kF \kappa(2, 5) = \kF \kappa (5, 2) = \kepsilon \kappa {13} = 0.019785233524272305\ldots.
	\]
	We have 
	\[
	\kN \kappa {\kepsilon{\kappa}{13}} = 2^5 \cdot 3^3 \cdot 5^2 \cdot 7 \cdot 11 \cdot 13 \cdot 17 = 2 \cdot 5 \cdot \kN \kappa {\kepsilon{\kappa}{12}}.
	\]
	However, both $2 \cdot \kN \kappa {\kepsilon{\kappa}{12}}$ and $5 \cdot \kN \kappa {\kepsilon{\kappa}{12}}$ are also $\kappa$-colossally abundant for $\kepsilon \kappa {13}$.
\end{example}

More broadly, if $\kEp \kappa p \cap \kEp \kappa q \neq \emptyset$ and $\epsilon = \epsilon_i \in \kEp \kappa p \cap \kEp \kappa q$, then $\kN \kappa {\epsilon_i}$ is a (likely trivial) multiple of $pq \kN \kappa {\epsilon_{i-1}}$, but $p \kN \kappa {\epsilon_{i-1}}$ and $q \kN \kappa {\epsilon_{i-1}}$ are also $\kappa$-colossally abundant numbers for $\epsilon_i$. Thus $\set{\kN \kappa \epsilon}_{\epsilon > 0}$ may be a proper subset of 
\[
\set{N \in \bbZ_{>0} \ : \ N \ \text{is} \ \kappa\text{-colossally abundant}}.
\]

\begin{theorem}\label{Theorem: Largest kappa-colossally abundant numbers}
	For all $\epsilon > 0$, we have \[
	\kN \kappa \epsilon = \prod_{p \ \text{prime}} p^{\kalpha \kappa p (\epsilon)},
	\]
	where
	\[
	\kalpha \kappa p (\epsilon) \coloneqq \floor{\ka \kappa p (\epsilon)}.
	\]
	If $\epsilon$ is $\kappa$-noncritical, then $\kN \kappa \epsilon$ is the only $\kappa$-colossally abundant number for $\epsilon$.
\end{theorem}

\begin{proof}
	Fix a colossally abundant number $N$ for $\epsilon$, and let $p$ be prime. Let $a$ be the $p$-adic valuation of $N$. By \Cref{Proposition: structure of k-colossally abundant numbers}, we see
	\[
	\kF \kappa (p, a+1) \leq \epsilon \leq \kF \kappa (p, a).
	\]
	On the other hand, we know from \Cref{Theorem: behavior of kF} that $\kF \kappa (p, a)$ is strictly decreasing in $a$, so by the definition of $\ka \kappa p$, the above inequality is equivalent to 
	\begin{equation*}\label{Equation 1 for kappa-colossally abundant numbers}
	a \leq \ka \kappa p (\epsilon) \leq a + 1, 
	\end{equation*}
	or equivalently,
	 \begin{equation}\label{Equation 2 for kappa-colossally abundant numbers}
	\ka \kappa p(\epsilon) - 1 \leq a \leq \ka \kappa p (\epsilon).
	\end{equation}
	Now if $\epsilon \not\in \kEp \kappa p$, then both inequalities in \eqref{Equation 2 for kappa-colossally abundant numbers} are strict, and so $a = \big\lfloor\ka \kappa p (\epsilon)\big\rfloor$. In particular, if $\epsilon$ is $\kappa$-noncritical then $\kN \kappa \epsilon$ is the only $\kappa$-colossally abundant number for $\epsilon$. Also, if $\epsilon \in \kEp \kappa p$, then $\ka \kappa p (\epsilon) \in \bbZ$, and \eqref{Equation 2 for kappa-colossally abundant numbers} implies $a \in \set{\ka \kappa p (\epsilon)-1, \ka \kappa p (\epsilon)}$. But as $\kN \kappa \epsilon$ is maximal by assumption, we again see that $a = \big\lfloor\ka \kappa p (\epsilon)\big\rfloor= \ka \kappa p (\epsilon)$. 
\end{proof}

\begin{cor}\label{Corollary: Nk divides Nk'}
	If $\kappa < \kappa\prm$, then $\kN \kappa \epsilon$ divides $\kN {\kappa\prm} \epsilon$.
\end{cor}

\begin{proof}
	\Cref{Theorem: behavior of kF} tells us $\kF \kappa$ is increasing in $\kappa$, so likewise $\ka \kappa p$ is increasing in $\kappa$ for each prime $p$. The claim is now immediate.
\end{proof}

We emphasize that \Cref{Corollary: Nk divides Nk'} does \emph{not} imply that $\kN \kappa {\kepsilon \kappa i}$ divides $\kN {\kappa\prm} {\kepsilon {\kappa\prm} i}$. 

\begin{example}
	We have $\kN 2 {\kepsilon 2 6} = 2^3 \cdot 3 \cdot 5 \cdot 7$, but $\kN {3} {\kepsilon {3} 6} = 2^3 \cdot 3^2 \cdot 5$. In this case, we have $\kepsilon 2 6 = \frac{\log(51/5)}{2\log 3} - 1 \approx 0.05696$, and $\kepsilon 3 6 = \frac{\log(511)}{3\log 7} - 1 \approx 0.06829$.
\end{example}
	
If $\kEp \kappa p \cap \kEp \kappa q = \kEp {\kappa\prm} p \cap \kEp {\kappa\prm} q = \emptyset$ for all primes $p \neq q$, then $\kN \kappa {\kepsilon \kappa i}$ and $\kN {\kappa\prm} {\kepsilon {\kappa\prm} i}$ must both have $i$ prime factors, counting multiplicity, so one divides the other exactly if one equals the other.

\Cref{Theorem: Largest kappa-colossally abundant numbers} can be reformulated in terms of $\kx \kappa a$ instead of $\ka \kappa x$ as follows:
\begin{equation}
	\kN \kappa \epsilon = \prod_{\ell \geq 1} \prod_{x_{\ell+1} < p \leq x_\ell} p^\ell,\label{Equation: k-colossally abundant numbers in terms of x_ell}
\end{equation}
where $p$ varies over the set of primes $\set{2, 3, 5, 7, \dots}$, and $x_\ell \coloneqq \kx \kappa \ell (\epsilon)$.

\begin{theorem}\label{Theorem: kN is increasing and constant on intervals and and its endpoints are colossally abundant for both ends}
	The function $\kN \kappa \epsilon$ is constant on each half-open interval $[\kepsilon \kappa i, \kepsilon \kappa {i-1})$. These constant values are distinct and increasing. The number $\kN \kappa {\kepsilon \kappa i}$ is also $\kappa$-colossally abundant for $\kepsilon \kappa {i+1}$.
\end{theorem}

\begin{proof}
	Note 
	\[
	\set{\kepsilon \kappa i}_{i = 1}^\infty = \set{\epsilon \in (0, \infty) \ : \ \ka \kappa p (\epsilon) \in \bbZ \ \text{for some prime} \ p}.
	\]
	Thus if $\epsilon, \epsilon\prm \in (\kepsilon \kappa i, \kepsilon \kappa {i-1})$, then for each prime $p$ and each integer $a$, we have the inequality $\ka \kappa p(\epsilon) - 1 \leq a \leq \ka \kappa p (\epsilon)$ if and only if $\ka \kappa p(\epsilon\prm) - 1 \leq a \leq \ka \kappa p (\epsilon\prm)$. We have shown that $\kN \kappa \epsilon$ is constant for $\epsilon \in (\kepsilon \kappa i, \kepsilon \kappa {i-1})$. But the function $\kalpha \kappa p (\cdot) = \floor{\ka \kappa p (\cdot)}$ is a composition of right-continuous functions, and so is right-continuous; thus
	\[
	\kalpha \kappa p (\kepsilon \kappa i) = \lim_{\epsilon \to (\kepsilon \kappa i)^+} \kalpha \kappa p(\epsilon);
	\]
	so in fact $\kN \kappa \epsilon$ is constant for $\epsilon \in [\kepsilon \kappa i, \kepsilon \kappa {i-1})$. Moreover, if $\epsilon < \kepsilon \kappa i < \epsilon\prm$, then for some prime $p$ (not necessarily unique) and some integer $a > 0$ we have
	\[
	\ka \kappa p (\epsilon) < a < \ka \kappa p (\epsilon\prm),
	\]
	so $\kalpha \kappa p (\epsilon) = a - 1 < a = \kalpha \kappa p (\epsilon\prm)$. This proves the second claim.
	
	It remains to show that $\kN \kappa {\kepsilon \kappa i}$ is $\kappa$-colossally abundant for $\kepsilon \kappa {i+1}$. For $p$ prime, if $\kepsilon \kappa {i+1} \not\in \kEp \kappa p$ then $\kalpha \kappa p (\kepsilon \kappa i) = \kalpha \kappa p (\kepsilon \kappa {i+1})$. If $\kepsilon \kappa {i+1} \in \kEp \kappa p$ then $\kalpha \kappa p (\kepsilon \kappa i) + 1 = \kalpha \kappa p (\kepsilon \kappa {i-1})$, but by \Cref{Proposition: structure of k-colossally abundant numbers}, replacing $\kalpha \kappa p (\kepsilon \kappa {i-1})$ with $\kalpha \kappa p (\kepsilon \kappa {i-1}) - 1 = \kalpha \kappa p (\kepsilon \kappa i)$ also yields a $\kappa$-colossally abundant number for $\kepsilon \kappa {i+1}$. Thus $\kN \kappa {\kepsilon \kappa i}$ is $\kappa$-colossally abundant for $\kepsilon \kappa {i+1}$ as desired.
\end{proof}

\Cref{Theorem: kN is increasing and constant on intervals and and its endpoints are colossally abundant for both ends} shows that a $\kappa$-colossally abundant number of the form $\kN \kappa \epsilon$ is in fact of the form $\kN \kappa {\kepsilon \kappa i}$ for some $i \geq 1$. It is therefore natural to index these $\kappa$-colossally abundant numbers by an integral parameter instead of a real parameter.

\begin{definition}\label{Definition: kNi}
	For $\kappa > 1$ given and for $i \geq 1$, we define
	\[
	\kNi \kappa i \coloneqq \kN \kappa {\kepsilon \kappa i}.
	\]
\end{definition}

\begin{cor}\label{Corollary: Infinitely many k-colossally abundant numbers}
	For any $\kappa > 1$, there are infinitely many $\kappa$-colossally abundant numbers.
\end{cor}

\begin{proof}
	We observe $\set{\kNi \kappa i}_{i = 1}^\infty$ is a strictly increasing sequence of $\kappa$-colossally abundant numbers.
\end{proof}

The following definition is motivated by \Cref{Intro Theorem: k-Gronwall's Theorem}.

\begin{definition}\label{Definition: kG}
	For $\kappa$ a positive real number, we define 
	\[
	\kG \kappa (n) \coloneqq \frac{\zeta(\kappa) \ksigma \kappa (n)}{(e^\gamma n \log \log n)^\kappa}.
	\]
\end{definition}

\begin{theorem}\label{Theorem: k-colossally abundant numbers yield local maxima of kG} %See Theorem 6.16
	If $\kNi \kappa i \leq n \leq \kNi \kappa {i+1}$, then 
	\[
	\kG \kappa (n) \leq \max\parent{\kG \kappa \parent{\kNi \kappa i}, \kG \kappa \parent{\kNi \kappa {i+1}}}.
	\]
\end{theorem}

\begin{proof}
	Write $\epsilon \coloneqq \epsilon_{i+1}$, $N_i \coloneqq \kNi \kappa {i}$, and $N_{i+1} \coloneqq \kNi \kappa {i+1}$ to clean up notation. By \Cref{Theorem: kN is increasing and constant on intervals and and its endpoints are colossally abundant for both ends}, $N_i$ and $N_{i+1}$ are both critical points of the function
	\[
	n \mapsto \frac{\ksigma \kappa (n)}{n^{\kappa(1 + \epsilon)}}
	\]
	that yield the same (global) maximum. The function $f(x) \coloneqq \epsilon x - \log\log x$ is concave upwards for $x > 1$, and so for $1 < a \leq \xi \leq b$, we have $f(\xi) \leq \max(f(a), f(b))$. In particular, we see
	\[
	\epsilon \log n - \log\log\log n \leq \max\parent{\epsilon \log N_i - \log \log \log N_i, \epsilon \log N_{i+1} - \log \log \log N_{i+1}};
	\]
	multiplying through by $\kappa > 0$, taking exponentials, and multiplying by $\zeta(\kappa) e^{-\kappa \gamma}$, we see
	\[
	\frac{\zeta(\kappa) n^{\kappa \epsilon}}{(e^\gamma \log \log n)^\kappa} \leq \max\parent{\frac{\zeta(\kappa) N_i^{\kappa \epsilon}}{(e^\gamma \log \log N_i)^\kappa}, \frac{\zeta(\kappa) N_{i+1}^{\kappa \epsilon}}{(e^\gamma \log \log N_{i+1})^\kappa}}.
	\]
	On the other hand, 
	\[
	\frac{\ksigma \kappa (n)}{n^{\kappa(1 + \epsilon)}} \leq \frac{\ksigma \kappa (N_i)}{N_i^{\kappa(1 + \epsilon)}} = \frac{\ksigma \kappa (N_{i+1})}{N_{i+1}^{\kappa(1 + \epsilon)}},
	\]
	and multiplying these two inequalities together yields the desired result.
\end{proof}

We close this section with a lemma that relates $\kx \kappa 1 (\epsilon)$ to $\kN\kappa \epsilon$, conditional on the Riemann hypothesis.

\begin{lemma}\label{Lemma: log log N versus log theta(x) with RH} %See Lemma 7.10
	Fix $\kappa\ge3/2$, and assume the Riemann hypothesis holds. If $x = \kx \kappa 1(\epsilon)$ and $N = \kN \kappa \epsilon$, then we have
	\begin{equation}
	\log \log N > \log \theta(x) \exp\parent{\frac{0.977\sqrt{2}}{\sqrt{x} \log x} + O\parent{\frac{1}{\sqrt{x} \parent{\log x}^2}}}\label{Equation: Ineffective log log N versus log theta x}
	\end{equation}
	for sufficiently large $x$. Moreover, if $\kappa = 2$ then
	\begin{equation}
	\log \log N > \log \theta(x) \exp\parent{\frac{0.977\sqrt{2}}{\sqrt{x} \log x}\parent{1 - \frac{\log 2}{2 \log x}}}\label{Equation: Effective log log N versus log theta x when k = 2}
	\end{equation}
	for $x \geq 2^{15}$.
\end{lemma}

\begin{proof}
We follow the proof of \cite[Lemma 7.10]{Broughan}. 
	
Suppose $x\ge2^{15}$. As usual, we write $x_\ell$ for $\kx \kappa \ell (\epsilon)$. By \eqref{Equation: k-colossally abundant numbers in terms of x_ell}, we see $\log N = \sum_{\ell \geq 1} \theta(x_\ell) \geq \theta(x) + \theta(x_2)$. Now by \Cref{Lemma: bound on theta(x)}, we see
	\[
	\log N > \theta(x) \parent{1 + \frac{0.985 x_2}{1.000081 x}}.
	\]
	Taking logarithms again, recalling that $\log(1+x)>\frac{x}{x+1}$ for $x > 0$, and using Lemma \Cref{Lemma: bound on theta(x)} and \Cref{Lemma:upper bounds for x_2}, we see
	\begin{align*}
	\log \log N &>\log \theta(x) + \log\parent{1 + \frac{0.985 x_2}{1.000081 x}} \\
	&>\log \theta(x) + \frac{0.985 x_2}{1.000081 x + 0.985 x_2} \\
	&> \log \theta(x) + \frac{0.985 x_2}{1.000081 x + 0.985 \sqrt{2 x}} \\
	&>\log \theta(x) + \frac{0.977399 x_2}{x} \\
	&>\log \theta (x) \parent{1 + \frac{0.977399 x_2}{x \log(1.000081x)}}\\
	&>\log \theta (x) \parent{1 + \frac{0.977391x_2}{x \log x}}.
	\end{align*}

	Now for real $m, b > 0$, and $t \in [0, b]$, if $m\leq \frac{\log\parent{1 + b}}{b}$ then $1 + t \geq e^{mt}$. But by \Cref{Lemma:upper bounds for x_2}, $x_2 < \sqrt{2 x}$, so for $x \geq 2^{15}$ we have
	\[
	\frac{0.977391x_2}{x \log x} <\frac{0.977391\sqrt{2}}{\sqrt{x} \log x} \leq \frac{0.977391}{2^7\log 2^{15}}<b:=0.000735.
	\]
	Taking $m\coloneqq \frac{\log\parent{1 + b}}{b}>0.9996$ and observing $0.977391m> 0.977$, we see
	\begin{equation}
	\log \log N >\log \theta (x) \exp\parent{\frac{0.977 x_2}{x \log x}}.\label{Equation: penultimate lower bound on log log N}
	\end{equation}
	Substituting \eqref{Equation: asymptotic for x2 with error} into \eqref{Equation: penultimate lower bound on log log N} yields \eqref{Equation: Ineffective log log N versus log theta x}, and substituting \eqref{Equation: lower bound on kx 2 2} into \eqref{Equation: penultimate lower bound on log log N} yields \eqref{Equation: Effective log log N versus log theta x when k = 2}.
\end{proof}

\section{An analogue to Robin's theorem}\label{Section: Robin's Theorem}

In 1984, Robin proved the following converse of Ramanujan's result.

\begin{theorem}[{\cite[p. 204]{Robin}}, Robin's Theorem]\label{Theorem: Classical Robin's Theorem}
	If the Riemann hypothesis fails, let $B$ be the supremum of the real parts of the nontrivial zeros of $\zeta(s)$, and let $b \in (1 - B, 1/2)$. There is a constant $c = c(b) > 0$ such that the inequality
	\[
	\sigma(n) > e^\gamma n \log \log n \parent{1 + \frac{c}{(\log n)^b}}
	\]
	holds for infinitely many $n$.
\end{theorem}

Our goal in this section is to prove the following result, which is a corollary to \Cref{Theorem: Classical Robin's Theorem}.

\begin{theorem}\label{Theorem: Robin's theorem}
	Let $\kappa > 1$ be a real number. If the Riemann hypothesis fails, let $B$ be the supremum of the real parts of the nontrivial zeros of $\zeta(s)$, and let $b \in (1 - B, 1/2)$. 	There is a constant $c = c(b) > 0$ independent of $\kappa$ such that the inequality
	\[
	\ksigma \kappa(n) > \frac{(e^\gamma n \log \log n)^\kappa}{\zeta(\kappa)} \parent{1 + \frac{c}{(\log n)^b}}^\kappa
	\]
	holds for infinitely many $n$.
\end{theorem}

\Cref{Theorem: Robin's theorem} gives the implication $\eqref{Condition: for n large ksigma (n) is bounded} \Longrightarrow \eqref{Condition: RH holds}$ of \Cref{Intro Theorem: k-Ramanujan-Robin Theorem}. Our proof of \Cref{Theorem: Robin's theorem} has two key components: the first is \Cref{Theorem: Classical Robin's Theorem} itself, and the second is the following lemma.

\begin{lemma}\label{Lemma: counterexamples accrue as kappa shrinks}
	Fix $n \in \bbZ_{>0}$. The function
	\[
	\kappa \mapsto \parent{\zeta(\kappa) \ksigma \kappa (n)}^{1/\kappa}.
	\]
	is smooth and monotonically decreasing in $\kappa \in (1, \infty)$; if $n > 1$, this function is strictly decreasing in $\kappa$.
\end{lemma}

The proof of \Cref{Lemma: counterexamples accrue as kappa shrinks} itself requires the following simple result.

\begin{lemma}\label{Lemma: (q/p)^x and (q^x - 1)/(p^x - 1) and and 1/x suffice}
	If $v > u > 1$, then 
	\[
	H(x;u,v)\coloneqq \parent{\frac{v^x - 1}{u^x - 1}}^{1/x}
	\]
	is a strictly decreasing function of $x\in(0,\infty)$.
\end{lemma}

\begin{proof}
We show that the monotonicity of $H(x;u,v)$ follows quickly from that of $F(x,a)$, which we proved in \Cref{Proposition: F is decreasing in x and a}. To this end, we consider
\[\log H(x;u,v)= \frac{1}{x}\log\parent{\frac{v^x - 1}{u^x - 1}}.\]
Suppose first that $u<v\le u^2$. Observe 
\[\log H(x;u,v) = (F((v/u)^x,(\log u)/\log(v/u)) + 1)\log(v/u).\]
Since for fixed $a > 1$, $F(x,a)$ is strictly decreasing in $x>1$, we see that $\log H(x;u,v)$ is strictly decreasing in $x > 0$. In the general case, we may assume $u^{2^k} < v\le u^{2^{k+1}}$, where $k \ge 0$ is some integer. Then 
\[\log H(x;u,v) = \sum_{i=0}^{k-1}\log H(x;u^{2^i},u^{2^{i+1}}) + \log H(x;u^{2^k},v). \]
By the special case that we just handled, each summand is strictly decreasing in $x>0$. Therefore, $\log H(x;u,v)$ is strictly decreasing in $x>0$, and so is $H(x;u,v)$.
\end{proof}

\begin{proof}[Proof of \Cref{Lemma: counterexamples accrue as kappa shrinks}]
	Suppose $\kappa\prm > \kappa > 1$. We will prove the function
	\[
	\parent{\frac{\ksigma \kappa (p^a)}{1 - p^{-\kappa}}}^{1/\kappa}
	\]
	is decreasing in $\kappa$ for each prime $p$ and each integer $a \geq 0$. Taking a product over all primes $p$, with $a = a(p) = 0$ for almost all $p$, will prove the result.
	
	Suppose first that $a = 0$. In this case, 
	\[
	\parent{\frac{\ksigma \kappa (p^a)}{1 - p^{-\kappa}}}^{1/\kappa} = \frac{1}{\parent{1 - p^{-\kappa}}^{1/\kappa}},
	\]
	so it suffices to prove that $\parent{1 - p^{-\kappa}}^{1/\kappa}$ is increasing. But if $\kappa\prm > \kappa$ then
	\[
	\parent{1 - p^{-\kappa\prm}}^{1/\kappa\prm} > \parent{1 - p^{-\kappa}}^{1/\kappa\prm} > \parent{1 - p^{-\kappa}}^{1/\kappa},
	\]
	so $\parent{1 - p^{-\kappa}}^{1/\kappa}$ is increasing as desired.
	
	Suppose now that $a > 0$. We recall that $\sigma(p^a) = p \sigma(p^{a - 1}) + 1$, so 
	\begin{align*}
	\parent{\frac{\ksigma \kappa (p^a)}{1 - p^{-\kappa}}}^{1/\kappa} &= \parent{\frac{\sigma(p^a)^\kappa - \sigma(p^{a-1})^\kappa}{1 - p^{-\kappa}}}^{1/\kappa} \\
	&= p \sigma(p^{a-1}) \parent{\frac{\parent{p + \sigma(p^{a-1})\inv}^\kappa - 1}{p^\kappa - 1}}^{1/\kappa}.
	\end{align*}
	Taking $u=p$ and $v = p + \sigma(p^{a-1})\inv $ in \Cref{Lemma: (q/p)^x and (q^x - 1)/(p^x - 1) and and 1/x suffice}, %and applying \Cref{Lemma: f^h < g^h}, 
	we see 
	\[
		\parent{\frac{\parent{p + \sigma(p^{a-1})\inv}^\kappa - 1}{p^\kappa - 1}}^{1/\kappa}	
	\]
	is a decreasing function of $\kappa$, and the result follows.
\end{proof}

The proof of \Cref{Theorem: Robin's theorem} is now straightforward.

\begin{proof}[Proof of \Cref{Theorem: Robin's theorem}]
	For any $n > 1$, we have
	\[
	\lim_{\kappa \to \infty} \parent{\zeta(\kappa) \ksigma \kappa (n)}^{1/\kappa} = \sigma(n).
	\]
	Thus by \Cref{Lemma: counterexamples accrue as kappa shrinks}, for fixed $\kappa > 1$ we have
	\[
	\parent{\zeta(\kappa) \ksigma \kappa (n)}^{1/\kappa} > \sigma(n).
	\]
	Suppose now that the Riemann hypothesis fails, and let $b$ and $c$ be as in  \Cref{Theorem: Classical Robin's Theorem}. There are infinitely many $n$ for which
	\[
	\parent{\zeta(\kappa) \ksigma \kappa (n)}^{1/\kappa} > \sigma(n) > e^\gamma n \log \log n \parent{1 + \frac{c}{(\log n)^b}}.
	\]
	Rearranging, we obtain
	\[
	\ksigma \kappa (n) > \frac{\sigma(n)^\kappa}{\zeta(\kappa)} > \frac{\parent{e^\gamma n \log \log n}^\kappa}{\zeta(\kappa)} \parent{1 + \frac{c}{(\log n)^b}}^\kappa,
	\]
	which is what we desired to show.
\end{proof}

\begin{remark}
	Robin's proof of \Cref{Theorem: Classical Robin's Theorem} depends on the inequality 
	\begin{equation}
	x_\ell > x^{1/\ell}.\label{Inequality: xell versus x^1/ell}
	\end{equation}
	For fixed $\kappa>1$ and $\ell\ge2$, \Cref{Proposition: asymptotic for xell} implies that \eqref{Inequality: xell versus x^1/ell} holds for sufficiently large $x$, but we have been unable to establish this inequality uniformly in $\ell$, nor have we been able to adapt Robin's original argument to route around this claim. The difficulty is due to $\kF \kappa (x, a)$ being more complicated than Robin's $F(x, a) = \lim_{\kappa \to \infty} \kF \kappa (x, a)$. We prove \Cref{Theorem: Robin's theorem} in the case $\kappa = 2$ in the appendix by exploiting the relative simplicity of $\kF 2 (x, a)$ to prove \eqref{Inequality: xell versus x^1/ell} uniformly in $\ell$. \end{remark}

\section{An analogue to Ramanujan's theorem}\label{Section: Ramanujan's Theorem}

In this section, we prove an ineffective and then an effective version of Ramanujan's theorem for $\ksigma \kappa$. The following material is inspired by \cite[Chapter 7]{Broughan}. 

\subsection{An ineffective theorem}

Our goal in this subsection is to prove the following analogue to Ramanujan's theorem, which gives the implication $\eqref{Condition: RH holds} \Longrightarrow \eqref{Condition: for n large ksigma (n) is bounded}$ in \Cref{Intro Theorem: k-Ramanujan-Robin Theorem}. 

\begin{theorem}\label{Theorem: Ineffective Ramanujan's Theorem} %See Theorem 7.2
	Let $\kappa > 3/2$ be a real number. If the Riemann hypothesis holds, then 
	\begin{equation}
	\ksigma \kappa(n) < \frac{(e^\gamma n \log \log n)^\kappa}{\zeta(\kappa)}\label{Equation: Robin's Inequality}
	\end{equation}
	for all $n$ sufficiently large.
\end{theorem}

We first prove that we can restrict our attention to $\kappa$-colossally abundant numbers.

\begin{lemma}\label{Lemma: Counterexamples can be chosen to be colossally abundant} %See Lemma 7.1
	Let $i_0$ be a positive integer. If \eqref{Equation: Robin's Inequality} holds for $n = \kNi \kappa i$ for all $i \geq i_0$, then \eqref{Equation: Robin's Inequality} holds all $n \geq \kNi \kappa {i_0}$.
\end{lemma}

\begin{proof}
	Let $n \geq \kN \kappa {i_0}$. By \Cref{Corollary: Infinitely many k-colossally abundant numbers}, there are $\kappa$-colossally abundant numbers $\kNi \kappa {i}$ and $\kNi \kappa {i+1}$ with 
	\[
	\kNi \kappa {i_0} \leq \kNi \kappa i \leq n \leq \kNi \kappa {i+1}.
	\]
	Suppose $n$ violates \eqref{Equation: Robin's Inequality}. In the language of \Cref{Definition: kG}, this means $1 \leq \kG \kappa (n)$. Then by \Cref{Theorem: k-colossally abundant numbers yield local maxima of kG}, we see
	\[
	 1 \leq \kG \kappa (n) \leq \max\parent{\kG \kappa \parent{\kNi \kappa i}, \kG \kappa \parent{\kNi \kappa {i+1}}},
	\]
	so at least one of $\kNi \kappa i$ and $\kNi \kappa {i+1}$ also violates Robin's inequality. Iterating this argument with $n$ larger than $\kNi \kappa {i+2}$, our claim follows.
\end{proof}

We are now in a position to prove \Cref{Theorem: Ineffective Ramanujan's Theorem}.

\begin{proof}[Proof of \Cref{Theorem: Ineffective Ramanujan's Theorem}]
	If there is some $n=N$ violating \eqref{Equation: Robin's Inequality}, then by \Cref{Lemma: Counterexamples can be chosen to be colossally abundant}, we may assume $N = \kN \kappa \epsilon$ for some $\epsilon$. Writing $x_\ell$ for $\kx \kappa \ell (\epsilon)$, \eqref{Equation: k-colossally abundant numbers in terms of x_ell} tells us
	\[
	\frac{\ksigma \kappa (N)}{N^{\kappa}} = \prod_{\ell \geq 1} \prod_{x_{\ell + 1} < p \leq x_\ell} \frac{\ksigma \kappa (p^\ell)}{p^{\kappa \ell}}.
	\] 
	Splitting off the first factor from the rest, and applying \Cref{Lemma: sigma k(p^ell)/p^k ell < (1-p^-k)/(1 - p^-1)^k}, we have
	\[
	\frac{\ksigma \kappa (N)}{N^{\kappa}} < \prod_{x_2 < p \leq x} \frac{\ksigma \kappa (p)}{p^{\kappa}} \prod_{p \leq x_2} \frac{\parent{1 - p^{-\kappa}}}{\parent{1 - p\inv}^{\kappa}}.
	\]
	Let $x \coloneqq x_1$, and assume $x \geq 2^{15}$. Expanding $\ksigma \kappa (p)/ p^\kappa$ and rearranging, we obtain
	\begin{equation}
%	\frac{\ksigma \kappa (N)}{N^{\kappa}} &= \prod_{x_2 < p \leq x} \parent{1 - p\inv}^{-\kappa}\parent{1 - p^{-2}}^{\kappa}\parent{1 - (p+1)^{-\kappa}} \prod_{\ell \geq 2} \prod_{x_{\ell + 1} < p \leq x_\ell} \frac{\ksigma \kappa (p^\ell)}{p^{\kappa \ell}} \nonumber \\
%	&< \prod_{x_2 < p \leq x} \parent{1 - p\inv}^{-\kappa}\parent{1 - p^{-2}}^{\kappa}\parent{1 - (p+1)^{-\kappa}} \prod_{p \leq x_2} \parent{1 - p\inv}^{-\kappa} \parent{1 - p^{-\kappa}} \nonumber \\
%	&= \prod_{p \leq x} \parent{1 - p\inv}^{-\kappa} \prod_{x_2 < p \leq x} \parent{1 - p^{-2}}^{\kappa} \prod_{p \leq x_2} \parent{1 - p^{-\kappa}} \prod_{x_2 < p \leq x} \parent{1 - (p+1)^{-\kappa}} \nonumber \\
	\frac{\ksigma \kappa (N)}{N^\kappa} < \prod_{p \leq x} \parent{1 - p\inv}^{-\kappa} \prod_{x_2 < p \leq x} \parent{1 - p^{-2}}^{\kappa} \prod_{p \le x} \parent{1 - p^{-\kappa}} \prod_{x_2 < p \leq x} \frac{\parent{1 - (p+1)^{-\kappa}}}{1 - p^{-\kappa}}. \label{Equality: bound on ksigma(N)/N^kappa}
	\end{equation}

	Let us consider each of the above products in turn. By \Cref{Lemma: Nicolas's bound} and \eqref{Equation: Ineffective log log N versus log theta x}, we have
	\begin{equation}
	\prod_{p \leq x} \parent{1 - p\inv}\inv %&\leq e^\gamma \log\theta(x) \exp\parent{\frac{2 + \beta}{\sqrt x \log x} + \frac {\alpha(x)}{\sqrt x \parent{\log x}^2}} \\
	\le e^{\gamma} \log\log N \exp\parent{\frac{2 + \beta - 0.977\sqrt{2}}{\sqrt x \log x} + O\parent{\frac {\alpha(x)}{\sqrt x \parent{\log x}^2}}},\label{Equation: Ineffective product over 1-p inv}
	\end{equation}
	where $\beta = \gamma + 2 - \log 4 \pi =0.04619...$ and $\alpha(x)$ is as in \eqref{Equation: definition of alpha(x)}.

	By \Cref{Lemma: asymptotic for x2 with error} and \Cref{Lemma: bound on product (1 - p^-2)}, we see
	\begin{align}
		\prod_{x_2 < p \leq x} \parent{1 - p^{-2}} &\leq \prod_{\sqrt{2x} < p \leq x} \parent{1 - p^{-2}} \leq \exp\parent{- \frac{\sqrt{2}}{\sqrt{x} \log x} + \frac{4}{\sqrt{x} \parent{\log x}^2}}.\label{Equation: asymptotic for (1 - p^-2)}
	\end{align}
		
	By \Cref{Lemma: Unconditional bound on prod (1 - p^-k)^-1}, we get
	\begin{equation}
	\prod_{p \leq x} \parent{1 - p^{-\kappa}} < \frac{1}{\zeta(\kappa)} \exp\parent{\frac{1.01624 \kappa x^{1 - \kappa}}{\log x} \parent{\frac{1}{\kappa - 1} + 0.000052}}.\label{Equation: product of 1-p^-k in Ramanujan's Theorem}
	\end{equation}

	Finally, since $\sqrt x < x_2$ according to \Cref{Lemma: asymptotic for x2 with error}, we find by  \Cref{Lemma: bound on (1 - (p+1)^-k)(1 - p^-k)} that
		\begin{equation}
		\prod_{x_2 < p \leq x} \frac{1 - (p+1)^{-\kappa} }{1 - p^{-\kappa}} < \prod_{p>\sqrt{x}} \frac{1 - (p+1)^{-\kappa} }{1 - p^{-\kappa}}<\exp\parent{\frac{2.1558 (\kappa + 1)}{\kappa x^{\kappa/2} \log x}}.\label{Equation: asymptotic for (1 - (p+1)^-kappa)/(1 - p^-kappa)}
	\end{equation}

	Substituting \eqref{Equation: Ineffective product over 1-p inv}, \eqref{Equation: asymptotic for (1 - p^-2)}, \eqref{Equation: product of 1-p^-k in Ramanujan's Theorem}, and \eqref{Equation: asymptotic for (1 - (p+1)^-kappa)/(1 - p^-kappa)} into \eqref{Equality: bound on ksigma(N)/N^kappa}, we obtain
	\begin{align*}
	\frac{\ksigma \kappa (N)}{N^{\kappa} (\log \log N)^\kappa} &< \frac{e^{\kappa \gamma} }{\zeta(\kappa)} \exp\parent{\frac{(2 + \beta - 1.977\sqrt{2}) \kappa}{\sqrt{x} \log x} + O\parent{\frac{1}{\sqrt{x} \parent{\log x}^2}}} \\ %\frac{B(x)}{\sqrt{x} \parent{\log x}^2}} \\
	&< \frac{e^{\kappa \gamma} }{\zeta(\kappa)} \exp\parent{\frac{-0.7497 \kappa}{\sqrt{x} \log x} + O\parent{\frac{1}{\sqrt{x} \parent{\log x}^2}}}.%+ \frac{B(x)}{x^{1/2} \parent{\log x}^2}},
	\end{align*}
	For $x$ sufficiently large, this expression is always less than $e^{\kappa \gamma}/\zeta(\kappa)$.
	\end{proof}

\subsection{An effective theorem}

Although \Cref{Theorem: Ineffective Ramanujan's Theorem} is quite general, it lacks the effectiveness asserted in \Cref{Intro Theorem: k-Ramanujan-Robin Theorem}. We adapt the proof of \Cref{Theorem: Ineffective Ramanujan's Theorem} and leverage \Cref{Lemma: lower bound on x2 when kappa = 2} to prove the following special case of \Cref{Intro Theorem: k-Ramanujan-Robin Theorem}. Combined with \Cref{Lemma: counterexamples accrue as kappa shrinks}, this will yield our main theorem.

\begin{theorem}\label{Theorem: Effective Ramanujan's Theorem} %See Theorem 7.2
	If the Riemann hypothesis holds, then 
	\begin{equation}
	\ksigma 2(n) < \frac{(e^\gamma n \log \log n)^2}{\zeta(2)}\label{Equation: Robin's Inequality when k = 2}
	\end{equation}
	for all $n > 2162160$.
\end{theorem}

\newcommand{\upperbound}{10^{10^{12.1408}}}

\begin{remark}
	There are a total of 79 counterexamples to \eqref{Equation: Robin's Inequality when k = 2} among $1 < n < \upperbound$. We list them now: 2, 3, 4, 5, 6, 7, 8, 9, 10, 12, 14, 15, 16, 18, 20, 24, 28, 30, 36, 40, 42, 48, 54, 60, 72, 84, 90, 96, 108, 120, 126, 132, 144, 168, 180, 210, 240, 252, 300, 336, 360, 420, 480, 504, 540, 600, 630, 660, 720, 840, 1080, 1260, 1320, 1440, 1680, 2520, 3360, 3780, 4200, 4620, 5040, 7560, 9240, 10080, 12600, 13860, 15120, 18480, 27720, 32760, 55440, 65520, 83160, 110880, 166320, 360360, 720720, 1441440, 2162160.

To verify that no counterexamples exist past 2162160 up to $\upperbound$, we used the code generously supplied by Platt to calculate the result of \cite[Theorem 5]{MorrillPlatt}, replacing the functions related to $\sigma$ with the analogous ones for $\sigma^{[2]}$.
\end{remark}

\begin{proof}[Proof of \Cref{Theorem: Effective Ramanujan's Theorem}]
	 Let $\kappa = 2$. We follow the proof of \Cref{Theorem: Ineffective Ramanujan's Theorem} exactly, save that we replace the ineffective inequality \eqref{Equation: Ineffective product over 1-p inv} with the effective inequality
	 \begin{align}
	\prod_{p \leq x} \parent{1 - p\inv}\inv %&\leq e^\gamma \log\theta(x) \exp\parent{\frac{2 + \beta}{\sqrt x \log x} + \frac {\alpha(x)}{\sqrt x \parent{\log x}^2}} \\
	%&= e^{\gamma} \log\log N \exp\parent{\frac{2 + \beta - 0.9749 \sqrt{2}}{\sqrt x \log x} + O\parent{\frac {\alpha(x)}{\sqrt x \parent{\log x}^2}}}, \\\
	&< e^{\gamma} \log\log N \exp\parent{\frac{2 + \beta - 0.977 \sqrt{2}}{\sqrt x \log x} + \frac {2 \alpha(x) +0.977 \sqrt{2} \log 2}{2 \sqrt x \parent{\log x}^2} }, \label{Equation: Effective product over 1-p inv}
	\end{align}
	using  \Cref{Lemma: Nicolas's bound} and \eqref{Equation: Effective log log N versus log theta x when k = 2}. Substituting \eqref{Equation: asymptotic for (1 - p^-2)}, \eqref{Equation: product of 1-p^-k in Ramanujan's Theorem}, \eqref{Equation: asymptotic for (1 - (p+1)^-kappa)/(1 - p^-kappa)}, and \eqref{Equation: Effective product over 1-p inv} into \eqref{Equality: bound on ksigma(N)/N^kappa}, we obtain an expression of the form
	\begin{align*}
	\frac{\ksigma \kappa (N)}{N^{\kappa} (\log \log N)^\kappa} &< \frac{e^{\kappa \gamma} }{\zeta(\kappa)} \exp\parent{\frac{-1.4994 + \epsilon(x)}{\sqrt{x} \log x}},
	\end{align*}
	where
	\[
	\epsilon(x) \coloneqq \frac {2\alpha(x) +0.977\sqrt{2} \log 2 + 4}{ \log x} + \frac{5.3}{\sqrt x}
	\]
	is decreasing in $x$ on $[2^{15}, \infty)$. But $\epsilon(2^{15}) < 1.1242< 1.4994$ by \Cref{Lemma: bound on theta(x)}, so \eqref{Equation: Robin's Inequality when k = 2} holds when $x \geq 2^{15}$. A short computation finishes the proof.
	\end{proof}

\begin{corollary}\label{Corollary: Effective Ramanujan's Theorem} %See Theorem 7.2
	Let $\kappa \geq 2$. If the Riemann hypothesis holds, then 
	\begin{equation}
	\ksigma \kappa(n) < \frac{(e^\gamma n \log \log n)^\kappa}{\zeta(\kappa)}\label{Equation: Robin's Inequality reprisal}
	\end{equation}
	for all $n > 2162160$.
\end{corollary}

\begin{proof}
	Let $\kappa \geq 2$. By \Cref{Lemma: counterexamples accrue as kappa shrinks}, we have
	\[
	\parent{\zeta(\kappa) \ksigma \kappa(n)}^{1/\kappa} \leq \parent{\zeta(2) \sigma^{[2]}(n)}^{1/2}. 
	\]
	But by \Cref{Theorem: Effective Ramanujan's Theorem}, 
	\[
	\parent{\zeta(2) \ksigma \kappa(n)}^{1/2} < e^\gamma n \log \log n
	\]
	whenever $n > 2162160$. The claim follows.
\end{proof}

When $\kappa \geq 2$, \Cref{Corollary: Effective Ramanujan's Theorem} furnishes the implication $\eqref{Condition: RH holds} \Longrightarrow\eqref{Condition: for n > 2162160 ksigma(n) is bounded}$ in \Cref{Intro Theorem: k-Ramanujan-Robin Theorem}, concluding its proof.

\section{An analogue to Lagarias's criterion}\label{Section: An analogue to Lagarias criterion}

Lagarias proved that the Riemann hypothesis is equivalent to the claim that the inequality \eqref{Equation: Lagarias criterion} holds for $n>1$ \cite[Theorem 1.1]{Lagarias}. He verified his criterion by means of the following lemma.

\begin{lemma}[Lagarias' Lemma]\label{Lemma: Lagarias' Lemma}
	For $n \geq 20$, we have
	\[
	e^\gamma n \log \log n + H_n \leq H_n + e^{H_n} \log H_n \leq e^\gamma n \log \log n + \frac{7 n}{\log n},
	\]
	where $H_n \coloneqq \sum_{1 \leq m \leq n} 1/m$ is the $n$th harmonic number.
\end{lemma}

\begin{proof}
	\cite[Lemmas 3.1 and 3.2]{Lagarias} (see also \cite[Lemma 7.17]{Broughan}).
\end{proof}

We use \Cref{Lemma: Lagarias' Lemma} to formulate an analogue to Lagarias' criterion.

\begin{theorem}\label{Theorem: Lagarias criterion} %See Theorem 7.18
	Let $\kappa \geq 2$ be given. The following are equivalent:
	\begin{enumerate}
		\item The Riemann hypothesis holds;\label{Item: RH}
		\item For $n > 55440$, we have\label{Item: Strong Lagarias inequality}
		\begin{equation}
		\ksigma \kappa(n) < \frac{\parent{e^{H_n} \log H_n}^\kappa}{\zeta(\kappa)};\label{Equation: Strong Lagarias inequality}
		\end{equation}
		\item For $n > 55440$, we have\label{Item: Weak Lagarias inequality}
		\begin{equation}
		\ksigma \kappa(n) < \frac{\parent{H_n + e^{H_n} \log H_n}^\kappa}{\zeta(\kappa)}.\label{Equation: Weak Lagarias inequality}
		\end{equation}
	\end{enumerate}
\end{theorem}

\begin{proof}
	$\eqref{Item: RH} \Longrightarrow \eqref{Item: Strong Lagarias inequality}$. By \Cref{Theorem: Effective Ramanujan's Theorem} and \Cref{Lemma: Lagarias' Lemma}, if the Riemann hypothesis holds, we have
	\[
	\ksigma \kappa (n) <\frac{\parent{e^\gamma n \log \log n}^\kappa}{\zeta(\kappa)} \leq \frac{\parent{e^{H_n} \log H_n}^\kappa}{\zeta(\kappa)}.
	\]
	which proves the first implication for $n > 2162160$. A short computation finishes this implication.
	
	$\eqref{Item: Strong Lagarias inequality} \Longrightarrow \eqref{Item: Weak Lagarias inequality}$. This is immediate.	
	
	$\eqref{Item: Weak Lagarias inequality} \Longrightarrow\eqref{Item: RH}$. Suppose by contradiction that for all $n$ sufficiently large, we have
	\[
		\ksigma \kappa(n) < \frac{\parent{H_n + e^{H_n} \log H_n}^\kappa}{\zeta(\kappa)},
	\]
	but the Riemann hypothesis fails. By \Cref{Theorem: Robin's theorem} and the inequality $e^{x/2}<1+x$ for all $x\in(0,2)$, there is a real number $b$ with $0 < b < \frac 12$ and a constant $c > 0$ such that
	\[
	\frac{(e^\gamma n \log \log n)^\kappa}{\zeta(\kappa)}\exp\parent{\frac{c}{(\log n)^b}} < \ksigma \kappa(n)	\]
	for infinitely many $n$. On the other hand, by \Cref{Lemma: Lagarias' Lemma} and by assumption, the following inequalities hold for all $n$ sufficiently large:
	\begin{align*}
		\ksigma \kappa(n) &< \frac{\parent{H_n + e^{H_n} \log H_n}^\kappa}{\zeta(\kappa)} \\
		&< \frac{(e^\gamma n \log \log n)^\kappa}{\zeta(\kappa)} \parent{1 + \frac{7}{e^\gamma \log n \log \log n}}^\kappa \\
		&< \frac{(e^\gamma n \log \log n)^\kappa}{\zeta(\kappa)} \exp\parent{\frac{7 \kappa}{e^\gamma \log n \log \log n}}.
	\end{align*}
	But for all $n$ sufficiently large, we have
	\[
	\frac{c}{(\log n)^b} > \frac{7 \kappa}{e^\gamma \log n \log \log n},
	\]
	so this is a contradiction.
\end{proof}

\begin{remark}\label{Remark: Lagarias counterexamples}
When $\kappa = 2$, there are a total of 33 counterexamples to \eqref{Equation: Strong Lagarias inequality} among $1 < n < \upperbound$. We list them now: 2, 3, 4, 6, 8, 10, 12, 18, 20, 24, 30, 36, 42, 48, 60, 72, 84, 90, 120, 168, 180, 240, 360, 420, 720, 840, 1260, 1680, 2520, 5040, 27720, 55440.
	
Similarly, when $\kappa = 2$, there are a total of 25 counterexamples to \eqref{Equation: Weak Lagarias inequality} among $1 < n < \upperbound$. We list them now: 2, 4, 6, 12, 18, 24, 30, 36, 48, 60, 72, 84, 120, 180, 240, 360, 420, 720, 840, 1260, 1680, 2520, 5040, 27720, 55440.
\end{remark}

\Cref{Theorem: Lagarias criterion} demands that $n > 55440$, but we need enlarge $\kappa$ only a little bit before the counterexamples we have enumerated disappear.

\begin{corollary}
Given $\kappa \geq 3.89$,
%$\kappa \geq 3.883305332$ 
the following are equivalent:
	\begin{itemize}
		\item The Riemann hypothesis holds;
		\item For $n > 1$, we have
		\begin{equation}\label{Equation: Lagarias inequality}
		\ksigma \kappa(n) < \frac{\parent{H_n + e^{H_n} \log H_n}^\kappa}{\zeta(\kappa)}.
		\end{equation}
	\end{itemize}
\end{corollary}

\begin{proof}
When $\kappa=3.89$, direct computation shows that \eqref{Equation: Lagarias inequality} holds for each of the counterexamples listed in the remark above. By \Cref{Lemma: counterexamples accrue as kappa shrinks}, \eqref{Equation: Lagarias inequality} also holds for every $\kappa\ge 3.89$. The corollary now follows from \Cref{Theorem: Lagarias criterion}.
\end{proof}

Of course, \Cref{Theorem: Ineffective Ramanujan's Theorem} also furnishes us with an ineffective version of the Lagarias criterion for $\kappa > 3/2$.

\begin{theorem}\label{Theorem: Ineffective Lagarias criterion} %See Theorem 7.18
	Let $\kappa > 3/2$ be given. The following are equivalent:
	\begin{enumerate}
		\item The Riemann hypothesis holds;
		\item For $n$ sufficiently large, we have
		\begin{equation*}
		\ksigma \kappa(n) < \frac{\parent{e^{H_n} \log H_n}^\kappa}{\zeta(\kappa)};
		\end{equation*}
		\item For $n$ sufficiently large, we have
		\begin{equation*}
		\ksigma \kappa(n) < \frac{\parent{H_n + e^{H_n} \log H_n}^\kappa}{\zeta(\kappa)}.
		\end{equation*}
	\end{enumerate}
\end{theorem}

\section{Future work}\label{Section: Future work}

We believe that our work makes a promising start on the study of $\ksigma \kappa (n)$, but much work remains to be done. A natural question is: can we reduce the threshold for $\kappa$ in \Cref{Corollary: Effective Ramanujan's Theorem} to $3/2$?

We also wish to highlight the discrepancy between the ranges of $\kappa$ in \Cref{Theorem: Ineffective Ramanujan's Theorem} and \Cref{Theorem: Robin's theorem}. Our proof of \Cref{Theorem: Ineffective Ramanujan's Theorem} fails for $1 < \kappa \leq 3/2$ because the right-hand side of \eqref{Equation: product of 1-p^-k in Ramanujan's Theorem} dominates the other terms in \eqref{Equality: bound on ksigma(N)/N^kappa}. Nonetheless, it would be interesting to know whether \Cref{Theorem: Ineffective Ramanujan's Theorem} still holds for smaller values of $\kappa$.
%\[
%\ksigma \kappa (N) > \frac{\parent{e^\gamma N \log \log N}^\kappa}{\zeta(\kappa)}.
%\]
%However, our efforts to prove this claim have run aground on the putative inequality \eqref{Inequality: xell versus x^1/ell}, just like our efforts to furnish a direct proof of \Cref{Theorem: Robin's theorem} for general $\kappa$.

All the work in this paper presumes that $\kappa > 1$, but it would be both natural and interesting to study $\ksigma \kappa (n)$ for $\kappa \leq 1$. By \eqref{Equation: ksigma(p^ell)/p^k ell}, we see
\[
\ksigma \kappa (p^\ell) \ksigma {-\kappa} (p^\ell) = \parent{1 - \ku \kappa(p, \ell)}\parent{1 - \ku {-\kappa} (p, \ell)},
\]
where
\[
\ku \kappa (p, \ell) \coloneqq \parent{\frac{1 - p^{-\ell}}{p - p^{-\ell}}}^\kappa = \frac{1}{\ku {-\kappa} (p, \ell)}.
\]
Thus, we expect some sort of duality between the behavior of $\ksigma \kappa (n)$ and $\ksigma {-\kappa}(n)$. The details of this correspondence, however, require further study.

Our analysis of $\ksigma \kappa (n)$ depended on our study of $\kappa$-colossally abundant numbers in \Cref{Section: kappa-colossally abundant numbers}. We noted after \Cref{Definition: N_epsilon} that if 
\begin{equation}
\kEp \kappa p \cap \kEp \kappa q \neq \emptyset\label{Equation: Ep cap Eq nonempty}
\end{equation}
then
\[
\set{\kN \kappa \epsilon}_{\epsilon > 0} \subsetneq \set{N \in \bbZ_{>0} \ : \ N \ \text{is} \ \kappa\text{-colossally abundant}}
\]
(see \Cref{Example: Ep and Eq intersect}). For $p \neq q$ fixed primes, can we find $\kappa$ for which \eqref{Equation: Ep cap Eq nonempty} holds? Does \eqref{Equation: Ep cap Eq nonempty} hold for any large $\kappa$? Does \eqref{Equation: Ep cap Eq nonempty} hold for any integral $\kappa$?

The inequality \eqref{Equation: Robin's classical inequality} is equivalent to the statement
\[
\sigma_{-1}(n) < e^\gamma \log \log n,
\]
which suggests an entirely different line of inquiry: rather than studying $\ksigma \kappa (n)$, we could study
\[
\sigma_{-1}^{[\kappa]}(n) = \sum_{d \mid n} \mu(n/d) \sigma_{-1}^\kappa(d),
\]
which satisfies
\[
\sigma_{-1}^{[k]}(n) = \sum_{[d_1, \dots, d_k] = n} \frac{1}{d_1 \dots d_k}
\]
whenever $\kappa = k \geq 1$ is an integer. In \eqref{Equation: Partial sums of sigma -1 kappa (n)} we determined an asymptotic for the partial sums of $\sigma_{-1}^{[2]}(n)$; we wish to generalize this asymptotic to all $\kappa > 1$ . It is not hard to show that for any $\kappa > 1$ we have
\[
\limsup_{n \to \infty} \frac{n \sigma_{-1}^{[\kappa]} (n)}{\kappa^{\omega(n)}(e^\gamma \log \log n)^{\kappa-1}} = 1,
\]
where $\omega(n)$ is the number of distinct factors of $n$. This analogue to \eqref{Equation: Gronwall's Theorem} and \Cref{Intro Theorem: k-Gronwall's Theorem} invites investigation, and we are interested in developing an analogue to the Ramanujan-Robin criterion for $\sigma_{-1}^{[\kappa]}(n)$.

More generally, one can investigate $\kappa$-analogues of Robin's criterion \eqref{Equation: Robin's classical inequality} for other arithmetic functions of interest. For example, Nicolas \cite{Nicolas1,Nicolas2} showed that the Riemann hypothesis implies that $\varphi(n)<e^{-\gamma}n/\log\log n$ for all primorials $n$ and that the falsity of the Riemann hypothesis would imply the reversed inequality for infinitely many primorials $n$. It would be interesting to prove analogues of \Cref{Intro Theorem: k-Ramanujan-Robin Theorem} for 
\[f^{[\kappa]}(n):=\sum_{d\mid n}\mu\left(\frac{n}{d}\right)f(d)^{\kappa}\]
with $f(d)=\varphi(d)$, and perhaps also with its cousin $f(d)=\psi(d)$, where 
\[\psi(d):=n\prod_{p\mid d}\left(1+\frac{1}{p}\right)\]
is the Dedekind totitent function.

In addition, despite the fact that Robin's inequality \eqref{Equation: Robin's classical inequality} is out of reach, Luca, Pomerance and Sol\'{e} \cite{LuPoSo} have considered the exceptional set of positive integers $n$ for \eqref{Equation: Robin's classical inequality}. They showed that the number of $n\le x$ for which \eqref{Equation: Robin's classical inequality} fails is at most $x^{O(1/\log \log x)}$ for $x\ge3$. It may be of interest to consider $\kappa$-analogues of their results for \eqref{Equation: k-Robin inequality}.

Finally, Washington and Yang \cite{Yang-Washington} and Vega \cite{vega} have published variants of the Ramanujan-Robin criterion where the domain of $\sigma(n)$ is restricted to special prime factorizations. Applying their methods to our functions furnishes another natural line of inquiry.

\section*{Appendix: Robin's theorem when \texorpdfstring{$\kappa = 2$}{k = 2}}

%See for instance 2022.04.14.sigma2Alone for a sketch of the arguement. E1 and E2 are both negligible; E1 I am pretty sure is negligible for the same reasons as in Robin's original proof, CONDITIONAL ON the inequality below, that xa > 1^(1/a). E2 is new, but easy to bound.

Robin's proof \cite{Robin} of \Cref{Theorem: Classical Robin's Theorem} could be adapted to prove \Cref{Theorem: Robin's theorem} without recourse to \Cref{Lemma: counterexamples accrue as kappa shrinks}, if we could somehow prove the inequality
\begin{equation}
x_a > x^{1/a}\label{Equation: x a > x^1/a}
\end{equation}
for $a \geq 2$.

Indeed, for colossally abundant $N$, Robin \cite[page 205]{Robin} writes 
\begin{equation}
\frac{\sigma(N)}{N} = \prod_{p \leq x} \parent{1 - p\inv}\inv \prod_{x_2 < p \leq x} \parent{1 - p^{-2}} \prod_{\ell \geq 2} \prod_{x_{\ell + 1} < p \leq x_\ell} \parent{1 - p^{-\ell - 1}}.\label{Equation: Robin's decomposition}
\end{equation}
Similarly, for $\kappa$-colossally abundant $N$, we have
\begin{equation}
\frac{\ksigma \kappa(N)}{N^\kappa} = \frac{1}{\zeta(\kappa)} \parent{\prod_{p \leq x} \parent{1 - p\inv}\inv\prod_{x_2 < p \leq x} \parent{1 - p^{-2}} \prod_{\ell \geq 2} \prod_{x_{\ell + 1} < p \leq x_\ell} \parent{1 - p^{-\ell - 1}}}^\kappa E_\kappa(x),\label{Equation: k-Robin's decomposition}
\end{equation}
where
\begin{equation}
E_\kappa(x) \coloneqq \zeta(\kappa) \prod_{\ell \geq 1} \prod_{x_\ell < p \leq x_{\ell + 1}} \parent{1 - \parent{\frac{1 - p^-\ell}{p - p^{-\ell}}}^\kappa}.\label{Equation: Definition of Ek}
\end{equation}
For $\kappa \geq 2$, a straightforward application of the prime number theorem to \eqref{Equation: Definition of Ek} shows unconditionally that 
\[
E_\kappa(x) = \exp\parent{O\parent{\frac{1}{x \log x}}}
\]
for $x$ sufficiently large, so the contribution of $E_\kappa(x)$ to \eqref{Equation: k-Robin's decomposition} is negligible. The remaining terms in \eqref{Equation: k-Robin's decomposition} are in visible correspondence with the terms in \eqref{Equation: Robin's decomposition}. Robin's handling of $\prod_{p \leq x} \parent{1 - p\inv}\inv$ goes through exactly, and by \Cref{Lemma:upper bounds for x_2} his handling of $\prod_{x_2 < p \leq x} \parent{1 - p^{-2}}$ does as well. However, without \eqref{Equation: x a > x^1/a}, we cannot follow Robin in handling 
\[\prod_{\ell \geq 2} \prod_{x_{\ell + 1} < p \leq x_\ell} \parent{1 - p^{-\ell - 1}}.\]

In the remainder of this appendix, we prove \eqref{Equation: x a > x^1/a} holds when $\kappa = 2$. By the argument sketched above, this gives us a direct proof of the 2-analogue of Robin's theorem.

By the monotonicity of $F_\kappa(x,a)$ in $x$, \eqref{Equation: x a > x^1/a} is equivalent to showing
\[F_\kappa(x_a,a)\leq F_\kappa(x^{1/a},a),\]
which is equivalent to
\[F_\kappa(x_1,1)\leq F_\kappa(x^{1/a},a).\]
This simplifies to
\begin{equation}
(x+1)^\kappa-1\leq \left[\frac{(x^{1+1/a}-1)^\kappa-(x-1)^\kappa}{(x-1)^\kappa-(x^{1-1/a}-1)^\kappa}\right]^a.\label{Equation: reframing of x a > x^1/a}
\end{equation}

Specializing now to the case $\kappa=2$, \eqref{Equation: reframing of x a > x^1/a} becomes
\[(x+1)^2-1\leq \left[\frac{(x^{1+1/a}-1)^2-(x-1)^2}{(x-1)^2-(x^{1-1/a}-1)^2}\right]^a.\]
Letting $y^a=x$, we have
\[(y^a+1)^2-1\leq \left[\frac{(y^{a+1}-1)^2-(y^a-1)^2}{(y^a-1)^2-(y^{a-1}-1)^2}\right]^a.\]
We then factor each of the differences of two squares:
\[y^a(y^a+2)\leq \left[\frac{(y^{a+1}+y^a-2)(y^{a+1}-y^a)}{(y^a+y^{a-1}-2)(y^a-y^{a-1})}\right]^a=y^a\left(\frac{y^{a+1}+y^a-2}{y^a+y^{a-1}-2}\right)^a.\]
Thus it remains to show that
\[y^a+2\leq\left(\frac{y^{a+1}+y^a-2}{y^a+y^{a-1}-2}\right)^a\]
for $y^a\geq 2$. This inequality rearranges to
\[ p(y)\coloneqq(y^{a+1}+y^a-2)^a-(y^a+2)(y^a+y^{a-1}-2)^a\geq 0, \]
so we show that the polynomial $p(y)\geq 0$ for $y^a\geq 2$. Observe that this polynomial has root 1 with multiplicity $a$ and root $2^{1/a}$ with multiplicity at least 1. We will show that in fact these are all of the positive roots by using Descartes' rule of signs and showing that the coefficients of this polynomial has exactly $a+1$ sign changes. Then it follows that $p(y)\geq 0$ for $y^a\geq 2$.

Using the binomial theorem, we find that the leading term of $p(y)$ is $2(a-1)y^{a^2}.$ Writing the $i$th coefficient of $p(y)$ as $c_i$, we thus have $c_{a^2}=2(a-1)$, and that $c_0=p(0)=-(-2)^a$.

For the remaining terms, we use the binomial expansions
\[(y^{a+1}+y^a-2)^a=\sum_{i=0}^a\sum_{j=0}^i\binom{a}{i}\binom{i}{j}(-2)^{a-i}y^{ai+j},\]
\[-y^a(y^a+y^{a-1}-2)^a=-\sum_{i=0}^a\sum_{j=0}^i\binom{a}{i}\binom{i}{j}(-2)^{a-i}y^{ai+a-i+j},\]
\[-2(y^a+y^{a-1}-2)^a=-\sum_{i=0}^a\sum_{j=0}^i\binom{a}{i}\binom{i}{j}(-2)^{a-i+1}y^{ai-i+j}.\]
By matching up the powers of $y$, we can determine the coefficients $c_m$ of $p$. We find for $0\leq m <a-1$ that
\[
c_{am}=(-2)^{a-m}\left(2\binom{a}{m-1}-\binom{a}{m}\right),
\]
along with the special case $m=a-1$ where $c_{a(a-1)}=-2(a-1)^2$.
When $1\leq r<a$, we have
\[c_{am+r}=(-2)^{a-m}\left[\binom{a}{m}\binom{m}{r}-\binom{a}{m}\binom{m}{r-a+m}+\binom{a}{m+1}\binom{m+1}{r-a+m+1}\right].\]

We may simplify the formula for $c_{am}$ when $m\neq 0$ by writing
\[\binom{a}{m}=\frac{a-m+1}{m}\binom{a}{m-1},\]
so
\[c_{am}=(-2)^{a-m}\left(2-\frac{a-m+1}{m}\right)\binom{a}{m-1}=(-2)^{a-m}\frac{3m-a-1}{m}\binom{a}{m-1}.\]
Thus, the sign of $c_{am}$ follows the sign of $(-1)^{a-m}(3m-a-1)$.

We now show that for each $m$ satisfying $0\leq m\leq a-1$, as $r$ increases from $1$ to $a-1$ the bracketed factor of $c_{am+r}$ always starts positive (or zero) and switches to negative at most once. 

 But as
\[\binom{a}{m+1}=\frac{a-m}{m+1}\binom{a}{m},\]
we observe that the sign of the bracket follows the sign of
\[B=B(r,m,a)=\binom{m}{r}-\binom{m}{r-a+m}+\frac{a-m}{m+1}\binom{m+1}{r-a+m+1},\]
so that
\[ c_{am+r}=(-2)^{a-m}\binom{a}{m}B. \]

First observe that since $a\geq m$, $B$ can only be negative if the central term is nonzero. This means $a-m\leq r\leq a$. We now focus on the difference
\[ D=-\binom{m}{r-a+m}+\frac{a-m}{m+1}\binom{m+1}{r-a+m+1}. \]
Since both binomial coefficients are nonzero when $a-m\leq r\leq a$, we have
\[ D=\left(\frac{a-m}{r-a+m+1}-1\right)\binom{m}{r-a+m} = \frac{2a-2m-r-1}{r-a+m+1}\binom{m}{r-a+m}. \]
The sign of $D$ follows the sign of the fraction. Since in our case the denominator is positive, $D$ negative precisely when $r>2a-2m-1.$

Finally, we determine when $B$ is negative by checking when we have
\[\binom{m}{r}<\frac{r-2a+2m+1}{r-a+m+1}\binom{m}{r-a+m}.\]
Since the right side is positive, the inequality certainly holds when the left side is zero, that is, when $r>m$. Now suppose the left side is positive, so $r\leq m$. Then we have
\begin{align*}
\frac{m!}{r!(m-r)!}&<\frac{r-2a+2m+1}{r-a+m+1}\frac{m!}{(r-a+m)!(a-r)!},\\
\frac{(a-r)!}{(m-r)!}&< \frac{r!(r-2a+2m+1)}{(r-a+m+1)!},
\end{align*}
from which it follows that
\[(a-r)(a-r-1)\cdots(m-r+1)< r(r-1)\cdots(r-a+m+2)(r-2a+2m+1).\]
Observe that each factor in the above inequality is nonnegative, and as $r$ increases, the left side decreases and the right side increases. Thus the direction of the inequality switches at most once, depending on the direction of the inequality when $r=m$.
Therefore, the sign of $B$ changes at most once, and this change can be detected when $B$ is negative for $r=a-1$. Explicitly, we find that
\[
B=B(a-1, m-1, a)=-\binom{m-1}{m-2}+
\frac{a-m+1}{m}\binom{m}{m-1}=a-2m+2<0
\]
when
\[m>\frac{a+2}{2}.\]

To summarize, we have shown that for fixed $m$ satisfying $0\leq m<a$, as $r$ increases from 1 to $a-1$, the bracketed factor of $c_{am+r}$ starts positive (or zero for $m=0$) and switches to negative at most once. 

It remains to compare the signs of $c_{am-1}, c_{am},$ and $c_{am+1}$ for each $m$ between 1 and $a-1$. We have that $c_{am}$ has the same sign as $(-1)^{a-m}$ if $m>(a+1)/3$, is zero if $m=(a+1)/3$, and has the opposite sign if $m<(a+1)/3$. For $r=1$ and $m=0$ we have $c_1=0$, and for $m>0$ the sign of $c_{am+1}$ is $(-1)^{a-m}$. Finally for $c_{am-1}$, we have that the sign is $(-1)^{a-m+1}$ for $a/2<m-1$, zero for $a/2=m-1$, and opposite for $a/2>m-1$.

We can now tally the sign changes of $c_i$. Observe that there is a one-to-one correspondence between coefficients $c_{am}, 0\leq m\leq a$ and the sign changes. It may help to refer to the figure below showing the case $a=20$ to see this.

\begin{figure}[h]
\begin{tikzpicture}
\draw[step=0.5cm, black, very thin] (0, -0.5) grid (0.5,10);
\foreach \i in {-0.25,0.75,...,9.75}{
    \draw (0.25, \i) node{$+$};
}
\foreach \i in {0.25,1.25,...,9.25}{
    \draw (0.25, \i) node{$-$};
}
\draw (0.5, 10.25) node{$(-1)^{a-m}$};
\draw[step=0.5cm,black,very thin] (1,0) grid (11,10);
\draw[black, very thin] (1,0)--(1,10);
\draw[black, very thin] (1,-0.5) rectangle (1.5,0);

\draw[black, thick](1.5, 6.5) -- +(90:3.5) node [midway, sloped, above]{negative} -- (1,10) -- (1, 6.5) -- (1.5, 6.5);

\draw[black, thick](1,6.5) -- (1, 6) -- (1.5, 6) -- (1.5, 6.5);
\draw (1.25, 6.25) node{$0$};

\draw[black, thick](1.5, 0) -- +(90:6) node [midway, sloped, above]{positive} -- (1,6) -- (1, 0) -- (1.5, 0);

\draw[black, thick](1,0) -- (1, -0.5) -- (1.5, -0.5) -- (1.5, 0);
\draw (1.25, -0.25) node{$+$};

\draw[black, thick](1.5, 10) -- (10.5,10) -- 
(10.5, 9.5) -- (10, 9.5) -- (10, 9) -- (9.5, 9) -- (9.5, 8.5) -- (9, 8.5) -- (9, 8) -- (8.5, 8) -- (8.5, 7.5) -- (8, 7.5) -- (8, 7) -- (7.5, 7) -- (7.5, 6.5) -- (7, 6.5) -- (7, 6) -- (6.5, 6) -- (6.5, 5.5) -- (5.5, 5.5)--
(5.5, 6) -- (5, 6) -- (5, 6.5) -- (4.5, 6.5) -- (4.5, 7) -- (4, 7) -- (4, 7.5) -- (3.5, 7.5) -- (3.5, 8) -- (3, 8) -- (3, 8.5) -- (2.5, 8.5) -- (2.5, 9) -- (2, 9) -- (2, 9.5) -- (1.5, 9.5);
\draw (6,8) node{$0$};

\draw[black, thick] (6.5, 0) -- (6.5, 2.5) -- (7, 2.5) -- (7, 3) -- (8, 3) -- (8, 3.5) -- (9, 3.5) -- (9, 4) -- (10, 4) -- (10, 4.5) -- (11, 4.5);
\draw[black, thick] (8.5,3.5) -- (8.5, 4) -- (9, 4);
\draw (8.75,3.75) node{$0$};
\draw[black, thick] (9.5, 4) -- (9.5, 4.5) -- (10, 4.5);
\draw (9.75, 4.25) node{$0$};
\draw[black, thick] (10.5, 4.5) -- (10.5, 5) -- (11, 5);
\draw (10.75, 4.75) node{$0$};

\draw[black, thick] (1.5, 0) -- (11, 0) -- (11, 10) -- (10.5, 10);

\draw (4, 4) node {positive};

\draw (9,2) node {negative};
\end{tikzpicture}
\centering
\caption{The signs of the coefficients of $p(y)$ in the case $a=20$.}
\end{figure}

In this figure we arrange the coefficients of $p(y)$ in a square array, starting with $c_0$ in the upper left corner, and then $c_r$ with $r$ increasing to the right from $r=1$ to $r=a-1=19$. The next row represents coefficients $c_a$ to $c_{2a-1}$, and so on. There is one extra coefficient that lives below the lower left corner for $c_{a^2}$. The sign of $B$ for each coefficient is labeled, and the signs of $(-1)^{a-m}$ are in a column on the left side for reference. Starting from the the leftmost column, we initially have sign changes when we leave the column. Then starting at the zero which occurs when $m=a/2=10$, the sign changes as we cross or enter the column. Then when $m-1\geq(a+1)/3=7$, or $m=8$, the behavior changes again. Now the sign changes start occurring earlier, but can still be made to pair up with entries in the leftmost column. Therefore, there are $a+1$ sign changes, establishing the original claim.

\section*{Acknowledgements}

We thank Dave Platt for sharing the code used in \cite{MorrillPlatt}. We thank Carl Pomerance for his helpful comments and feedback.

\bibliographystyle{amsplain}
%\bibliography{LCM_Bibliography}
\providecommand{\bysame}{\leavevmode\hbox to3em{\hrulefill}\thinspace}
\providecommand{\MR}{\relax\ifhmode\unskip\space\fi MR }
% \MRhref is called by the amsart/book/proc definition of \MR.
\providecommand{\MRhref}[2]{%
  \href{http://www.ams.org/mathscinet-getitem?mr=#1}{#2}
}
\providecommand{\href}[2]{#2}

\end{document}